\documentstyle[12pt]{article} 
\def\ds{\displaystyle}
\def\be{\begin{equation}} 
\def\ee{\end{equation}} 
\newtheorem{th}{Theorem}[section]  
\newtheorem{lm}[th]{Lemma}
\newtheorem{cor}[th]{Corollary}
\newtheorem{pr}[th]{Proposition}

\def\RR{{\hbox{I\kern-.2em\hbox{R}}}}
\def\PP{{\hbox{I\kern-.2em\hbox{P}}}}
\def\EE{{\hbox{I\kern-.2em\hbox{E}}}}
\def\ZZ{{\hbox{Z\kern-.4em\hbox{Z}}}}
\def\bx{\mbox{\boldmath $x$}}
\def\bz{\mbox{\boldmath $z$}}
\def\one{{\bf 1}}

\def\R{{\bf R}}
\def\N{{\bf N}}

\def\wh{\widehat}
\def\wt{\widetilde}
\def\Zp{{\bf Z}_+}
\def\qed{\hfill\hbox{\rule{6pt}{6pt}}}

\def\lg{\langle}
\def\rg{\rangle}

\def\Hc{\check{H}}
\def\Hh{\widehat{H}}
\def\Cc{\check{C}}
\def\Ch{\widehat{C}}

\def\cE{{\cal E}}
\def\cF{{\cal F}}

\def\cM{{\cal M}}
\def\cN{{\cal N}}
\def\cP{{\cal P}}

\def\cT{{\cal T}}

\def\cL{{\cal L}}

\begin{document} 
\setlength{\topmargin}{-4mm}
%
%
\begin{center}
{\LARGE \bf 
The coagulation-fragmentation hierarchy 
with homogeneous rates  
and underlying stochastic dynamics 
}
\footnote{
Research supported in part by JSPS, 
Grant-in-Aid for Scientific Research (C) No. 26400142 
}
\end{center} 

%
%


\begin{center}
Kenji Handa 
\end{center} 

\begin{center}
{ 
Department of Mathematics \\
Saga University \\ 
Saga 840-8502 \\ 
Japan
} \\  
e-mail: 
handak@cc.saga-u.ac.jp \\ 
FAX: +81-952-28-8501
\end{center} 
\begin{center}
{\em Dedicated to Professor Hiroshi Sugita  
on the occasion of his 60th birthday}
\end{center}

%
%

\begin{center} 
\begin{minipage}[t]{12cm} 
\small 
A hierarchical system of equations is introduced 
to describe dynamics of `sizes' of infinite 
clusters which coagulate and fragmentate 
with homogeneous rates of certain form. 
We prove that this system of equations is solved weakly  
by correlation measures for stochastic dynamics of  
interval partitions evolving according to 
some split-merge transformations. 
Regarding those processes, a sufficient condition 
for a distribution to be reversible is given.  
Also, an asymptotic result for properly rescaled processes 
is shown to obtain a solution to a nonlinear equation 
called the coagulation-fragmentation equation. 
\end{minipage}
\end{center}

%
%

\section{Introduction} 
\setcounter{equation}{0} 
The phenomena of coagulation and fragmentation 
are studied in various contexts of natural sciences. 
Mathematically, they are considered to be `dual' to 
each other at least in some naive sense 
or to be simply the time-reversal of each other. 
Hence, one naturally expects that 
the coagulation-fragmentation dynamics 
may lead to a nontrivial equilibrium 
in the course of time provided the occurrence of 
coagulation and fragmentation is prescribed 
in an appropriate manner. 
Many authors have examined 
such situations through a nonlinear equation 
called often the coagulation-fragmentation equation. 
It takes the form 
\begin{eqnarray} 
\frac{\partial}{\partial t}c(t,x)
& = & 
\frac{1}{2}\int_0^x\left[K(y,x-y)c(t,y)c(t,x-y)
-F(y,x-y)c(t,x)\right]dy               \nonumber  \\ 
& & 
-\int_0^{\infty}\left[K(x,y)c(t,x)c(t,y)
-F(x,y)c(t,x+y)\right]dy,     \label{1.1} 
\end{eqnarray} 
where $t, x>0$, and 
the functions $K$ and $F$ are supposed to be given, 
nonnegative, symmetric and depending on the mechanisms 
of coagulation and fragmentation, respectively. 
In the literature $c(t,x)$ represents 
the `density' of clusters of size $x$ 
(or particles with mass $x$) at time $t$ 
and the equation (\ref{1.1}) is derived 
heuristically by some physical arguments 
or rigorously for some restricted cases. 
(Among results of the latter kind 
for both nonzero $K$ and $F$, 
we refer \cite{EW00}.) 
However, (\ref{1.1}) is not complete 
for the full description of coagulation-fragmentation phenomena 
since it usually emerges after 
certain contraction procedure such as `propagation of chaos' 
or under intuitive assumptions of asymptotic 
independence among distributions of clusters.

In this paper we study a hierarchical system of equations, 
for a special case of which 
we establish a direct connection with 
an infinite-dimensional stochastic dynamics 
incorporating coagulation and fragmentation. 
For each $k\in \N:=\{1,2,\ldots\}$, 
the $k$th equation of the hierarchy reads 
\begin{eqnarray} 
\lefteqn{
\frac{\partial}{\partial t}c_k(t,z_1,\ldots,z_k)} 
                                         \label{1.2}  \\ 
&=& 
\frac{1}{2}\sum_{l=1}^k\int_0^{z_l}K(y,z_l-y)
c_{k+1}(t,z_1,\ldots,z_{l-1},y,z_l-y,z_{l+1},\ldots,z_k)dy 
                                       \nonumber  \\ 
& & 
-\frac{1}{2}\sum_{l=1}^k\int_0^{z_l}F(y,z_l-y)dy \ 
c_{k}(t,z_1,\ldots,z_k) 
                                       \nonumber  \\ 
& & 
-\sum_{l=1}^k\int_0^{\infty}K(z_l,y)
c_{k+1}(t,z_1,\ldots,z_l,y,z_{l+1},\ldots,z_k)dy    
                                       \nonumber  \\ 
& & 
+\sum_{l=1}^k\int_0^{\infty}F(z_l,y)
c_{k}(t,z_1,\ldots,z_{l-1},z_l+y,z_{l+1},\ldots,z_k)dy  
                                       \nonumber  \\ 
& &  
-\one_{\{k\ge 2\}}\sum_{l<m}^{k}K(z_l,z_m)
c_{k}(t,z_1,\ldots,z_k) 
                                       \nonumber  \\ 
& &  
+\one_{\{k\ge 2\}}\sum_{l<m}^{k}F(z_l,z_m)
c_{k-1}(t,z_1,\ldots,z_{l-1},z_l+z_m,z_{l+1},\ldots, 
z_{m-1},z_{m+1},\ldots,z_k),            \nonumber 
\end{eqnarray} 
where $\one_E$ stands in general for the indicator function of 
a set $E$ and the sum $\sum_{l<m}^{k}$ is taken over 
pairs $(l,m)$ of integers such that $1\le l<m \le k$. 
If the last two terms on the right side of (\ref{1.2}) 
were absent, it is readily checked that 
the system of equations is satisfied by the direct products 
$c^{\otimes k}(t,z_1,\ldots,z_k):=c(t,z_1)\cdots c(t,z_k)$ 
of a solution to (\ref{1.1}). The equations 
(\ref{1.2}) are considered to be much more informative 
in the sense that interactions among an 
arbitrary number of clusters are took into account. 

In fact, a finite-system version of (\ref{1.2}) 
has been discussed for a pure coagulation model 
(i.e. the case $F\equiv 0$) by Escobedo and Pezzotti \cite{EP}. 
Their derivation of (\ref{1.2}) with $F\equiv 0$ 
starts from a finite set of evolution equations satisfied 
by the so-called mass probability functions 
associated with a stochastic coagulation model 
known as the Marcus-Lushnikov process. 
As pointed out in \cite{EP} the situation is similar 
to the derivation of the BBGKY hierarchy in classic 
kinetic theory although the underlying microscopic dynamics 
for the BBGKY hierarchy is not stochastic but deterministic. 
There is an extensive literature discussing 
a stochastic dynamics which serves 
as a basis of an infinite system called the Boltzmann hierarchy, 
which is a thermodynamic limit of the BBGKY hierarchy. 
(See the monograph of Petrina \cite{Petrina} and the references therein.)  
It should be mentioned also that a number of articles 
have discussed stochastic interacting systems 
of finite particles to derive kinetic equations, 
a special case of which is (\ref{1.1}),
in the limit as the number of particles tends to infinity. 
(See e.g. a paper by Eibeck and Wagner \cite{EW}  
and the references therein. Also, 
for a systematic treatment in a general framework  
related to such issues, see 
Kolokoltsov's monograph \cite{Kol}.) 
We intend to explore 
the `coagulation-fragmentation hierarchy' (\ref{1.2}) 
by dealing with stochastic infinite systems directly 
and derive (\ref{1.1}) as a macroscopic equation 
for them through the limit under proper rescaling. 
Such a limit theorem is regarded as 
the law of large numbers for measure-valued processes 
and related to the propagation of chaos. 
(See {\it Remarks} at the end of \S 4.1 below.) 

In the case where the mechanisms of coagulation and fragmentation  
together enjoy the detailed balance condition, i.e., 
\be 
K(x,y)M(x)M(y)=F(x,y)M(x+y)              \label{1.3}
\ee  
for some function $M$, equilibrium behaviors of 
the solution $c(t,x)$ to (\ref{1.1}) with respect to 
a stationary solution of the form 
$x\mapsto M(x)e^{-bx}$ have been studied by 
many authors. In particular, for the equation with 
$K$ and $F$ being positive constants,  
Aizenman and Bak \cite{AB} carried out 
detailed analysis such as a uniform rate for 
strong convergence to equilibrium. 
(See also Stewart and Dubovski \cite{SD}.) 
Lauren{\c c}ot and Mischler \cite{LM} studied 
that convergence 
under certain assumptions for $K, F$ and $M$ 
and suitable conditions on the initial state. 
Such results include particularly an $H$-theorem, 
namely the existence of a Lyapunov functional 
of entropy type for the solution. 

In what follows, 
we shall restrict the discussion to the case where 
\be 
K(x,y)=xy \Hh(x,y), \quad 
F(x,y)= (x+y) \Hc(x,y)                     \label{1.4} 
\ee 
for some homogeneous functions $\Hh$ and $\Hc$ of 
common degree $\lambda\ge 0$, namely,
\be 
\Hh(ax,ay)=a^{\lambda}\Hh(x,y),  \quad 
\Hc(ax,ay)=a^{\lambda}\Hc(x,y)  \quad (a, x , y> 0).     \label{1.5} 
\ee 
To avoid trivialities, we suppose also 
that $\Hh$ and $\Hc$ are not identically zero. Therefore, 
both $K$ and $F$ are necessarily unbounded but 
of polynomial growth at most. 
In case $\theta\Hh=\Hc$ for a constant 
$\theta>0$, (\ref{1.3}) holds for $M(x)=\theta/x$.  

The coagulation and fragmentation phenomena 
have been discussed also in the probability literature. 
See e.g. Bertoin's monograph \cite{B} for systematic accounts of 
stochastic models and random operations 
describing either phenomenon. 
The choice (\ref{1.4}) is mainly motivated by 
a coagulation-fragmentation process 
studied by Mayer-Wolf et al \cite{MWZZ} and Pitman \cite{P}. 
These papers concern the case where 
$\Hh$ and $\Hc$ are constants. 
Having the infinite-dimensional simplex 
\[ 
\Omega_1
=\{\bx=(x_i)_{i=1}^{\infty}:~ x_1\ge x_2 \ge \ldots \ge 0, 
\ \sum_{i}x_i=1 \} 
\] 
as its state space, the process keeps 
the total sum 1 of cluster sizes fixed.  
The special case $\Hh=\Hc \equiv \mbox{const.}$ 
corresponds to the Markov process 
explored in \cite{T1} and \cite{T2}, 
which is associated with `the simplest split-merge operator' 
originally introduced by A. Vershik in the context of 
analysis of the infinite dimensional symmetric group. 
This model was studied extensively in \cite{DMWZZ} 
in a deep and explicit connection with 
a discrete analogue generated by the random transposition, 
for which one may refer to \cite{Sch} for instance. 
For that discrete model, the  
coagulation and fragmentation rates (\ref{1.4}) with 
both $\Hh$ and $\Hc$ being constants,   
naturally emerge as transition probabilities. 
(See (2.2) in \cite{DMWZZ}.)  
In these works it was shown that the celebrated 
Poisson-Dirichlet distribution with parameter $\theta$ 
is a reversible distribution of the process, and 
much efforts were made to prove 
the uniqueness of a stationary distribution. 
In particular, Diaconis et al \cite{DMWZZ} 
succeeded in proving it for $\theta=1$ by giving 
an effective coupling result with the discrete 
coagulation-fragmentation processes. 
(See also Theorem 1.2 in \cite{Sch} and 
Theorem 7.1 in \cite{Gold}.) 
We also refer the reader to \cite{Be} 
for another result of interest 
on a unique stationary distribution for 
the model evolving with a different class of 
coagulation-fragmentations. 
The coagulation and fragmentation 
we will be concerned with are only binary ones. 
Cepeda \cite{C2} constructed stochastic models 
incorporating coagulation and multi-fragmentation 
on a larger state space than that of our models 
(, i.e., $\Omega$ defined below). See Introduction 
and the references in \cite{C2} for previous 
works and development in the study of 
related stochastic models.

By virtue of the homogeneity assumption on $\Hh$ and $\Hc$   
we can consider the generalized process 
associated with (\ref{1.4}) not only on $\Omega_1$ 
but also on the infinite-dimensional cone 
\[ 
\Omega 
=\{\bz=(z_i)_{i=1}^{\infty}:~ z_1\ge z_2 \ge \ldots \ge 0, 
\ 0<\sum_{i}z_i<\infty \}.  
\] 
In fact, the major arguments in this paper 
exploit some ingredients from theory of point processes. 
The idea is that each $\bz=(z_i)\in\Omega$ 
can be identified with the locally finite point-configuration 
\[ 
\xi=\sum_{i}\one_{\{z_i>0\}}\delta_{z_i}
\] 
on the interval $(0,\infty)$, where 
$\delta_{z_i}$ is the delta distribution concentrated at $z_i$. 
Indeed, in one of our main results, 
the notion of correlation measures 
will make us possible to reveal 
an exact connection between 
hierarchical equations (\ref{1.2}) and 
the coagulation-fragmentation process with rates (\ref{1.4}). 
As another result based on the point process calculus 
we will present a class of 
coagulation-fragmentation processes having 
the Poisson-Dirichlet distributions or 
certain variants (including the laws of gamma 
point processes) as their reversible distributions, 
clarifying what mathematical structures 
are responsible for this result. 
That structure will be described in terms of 
correlation functions 
together with Palm distributions, 
certain conditional laws for the point process. 
(See (\ref{2.14}) and (\ref{2.15}) below.) 
We mention also that 
the reversibility will play some key roles  
in discussing the existence of strong solutions to 
(\ref{1.2}). As for the original equation (\ref{1.1}), 
introducing rescaled models 
which depend on the scaling parameter $N$, 
one can discuss its derivation from 
the associated measure-valued processes 
in which each point is assigned mass $1/N$. 
Such a result is formulated as a limit theorem  
for Markov processes as $N\to\infty$ 
and one of the key steps is 
to show the tightness of their laws, 
which is far from trivial 
because there is less restriction on 
grows order of $K$ and $F$. 
As will turn out later our setting of 
the degrees of homogeneity plays  
an essential role to overcome difficulties of this sort. 

The organization of this paper is as follows. 
In Section 2, we introduce 
the coagulation-fragmentation process 
associated with rates (\ref{1.4})  
and give an equivalent description of the model 
in terms of the corresponding point process. 
Section 3 discusses 
a weak version of (\ref{1.2}), 
which will turn out to be satisfied by the 
correlation measures of 
our coagulation-fragmentation process. 
After some preliminary arguments are made 
for rescaled models in Section 4, 
a solution to (\ref{1.1}) will be obtained 
from properly rescaled empirical measures 
in Section 5. 

\section{The coagulation-fragmentation process 
associated with split-merge transformations}
\setcounter{equation}{0} 
\subsection{Definition of the models} 
As mentioned in Introduction 
the rates $K$ and $F$ are supposed to 
be of the form (\ref{1.4}) with $\Hh$ and $\Hc$ being 
homogeneous functions of degree $\lambda\ge0$  throughout. 
Notice that this is equivalent to the condition that 
$K$ and $F$ are homogeneous functions  
of degree $\lambda+2$ and $\lambda+1$, 
respectively, for some $\lambda\ge 0$. 
As far as coagulation rates are concerned, 
the homogeneity, though mathematically a strong condition, 
is satisfied typically by examples of kernels 
used in the physical literature as seen in Table 1 in \cite{A}  
(although our framework excludes any of such examples).  
See also \cite{FL} and \cite{C}, 
which discuss the equation with homogeneous(-like) $K$. 
In the rest, the following two conditions are 
also imposed without mentioning:  \\ 
(H1) $\Hh$ is a symmetric, nonnegative 
measurable function on $(0,\infty)^2$ such that 
\[ 
\Ch:=\sup\{\Hh(u,1-u) :~0<u<1 \}\in(0,\infty).  
\] 
(H2) $\Hc$ is a symmetric, nonnegative 
measurable function on $(0,\infty)^2$ such that 
\[ 
\Cc:= \int_0^1 \Hc(u,1-u)du \in(0,\infty).  
\] 
In general, a homogeneous function $H$ 
on $(0,\infty)^2$ is determined by 
its degree $\lambda$ and the function 
$h(u):=H(u,1-u)$ on $(0,1)$ through the relation 
$H(x,y)=(x+y)^{\lambda}h(\frac{x}{x+y})$. 
As for the fragmentation rate, 
the homogeneity (\ref{1.5}) combined 
with (\ref{1.4}) implies that the overall rate of 
fragmentation of an $x$-sized cluster is 
necessarily given by the power-law form:   
\[ 
\frac{1}{2}\int_0^xF(y,x-y)dy 
= \frac{1}{2}x^2\int_0^1\Hc(ux,(1-u)x)du
= \frac{\Cc}{2}x^{2+\lambda}. 
\] 
Such a situation is featured by the 
coagulation-fragmentation equation 
studied in \cite{BLL}, whose conditions 
for coagulation rates are also well adapted to our setting. 
(See {\it Remark} at the end of Section 5 
for related discussions.) 
Also, \cite{VZ} examined the interplay between 
degrees of coagulation and fragmentation 
in the context of stability analysis. 
At the beginning of Section 4 we will 
mention another role of such interplay 
between them in the study of rescaled processes. 
\medskip  \\ 
\noindent  
{\it Examples.}~ 
(i) Consider $H(x,y)=(xy)^{a}(x^b+y^b)^c$,  
for which $\lambda=2a+bc$. 
Without loss of generality, we can suppose that $b\ge 0$. Then 
$H$ satisfies (H1) (resp. (H2)) if $a\ge 0$ (resp. $a>-1$). \\ 
(ii) Let  $b>0$ and define ${H}(x,y)=(xy)^{a}|x^b-y^b|^c$, 
for which $\lambda=2a+bc$. 
${H}$ satisfies (H1) (resp. (H2)) 
if $a\ge 0$ and $c\ge 0$ (resp. $a>-1$ and $bc>-1$). 
Indeed, it is readily observed that 
${H}(u,1-u)/u^a \to 1 (u \downarrow 0)$ 
and 
${H}(u,1-u)/|2u-1|^{bc} \to b2^{-(2a+bc)} 
(u \to 1/2)$.  \\ 
(iii) Given $\lambda\ge 0$, let 
${H}(x,y)=(x^a+y^a)(x^{\lambda-a}+y^{\lambda-a})$. 
It follows that (H1) (resp. (H2)) is satisfied 
if $0\le a\le \lambda$ (resp. $-1<a<\lambda+1$).        \\ 
(iv) An example of discontinuous homogeneous function 
of degree $\lambda$ is 
$H(x,y)=x^{\lambda}\one_{\{x\ge ay\}}
+y^{\lambda}\one_{\{y\ge ax\}}$, where $a>0$. 
Another one is $H(x,y)=(xy)^{\lambda/2}
\one_{[a,1/a]}(x/y)$ with $0<a<1$. 
For each example and any $\lambda\ge 0$, 
both  (H1)  and (H2) hold. 

\medskip 

In order to define our coagulation-fragmentation process 
as a continuous-time Markov process on $\Omega$, 
suitable modifications to the formulation 
in \cite{MWZZ} are made in the following manner. 
Let $\Omega$ be equipped with the product topology 
and $B(\Omega)$ (resp. $B(\Omega_1)$) be the Banach space of 
bounded Borel functions on $\Omega$ (resp. $\Omega_1$) 
with the sup norm $\Vert \cdot \Vert_{\infty}$.  
For $\bz=(z_i)\in \Omega$ put $|\bz|=\sum z_i$. 
A useful inequality is $\sum z_i^{1+a}\le |\bz|^{1+a}$ 
for any $a\ge 0$, which is implied by 
$\sum (z_i/|\bz|)^{1+a}\le \sum (z_i/|\bz|) =1$. 
Define a bounded linear operator $\wt{L}$ on $B(\Omega)$ by 
\begin{eqnarray} 
\wt{L}\Phi(\bz)
& = & 
\frac{1}{2|\bz|^{2+\lambda}}\sum_{i\ne j}K(z_i,z_j)
\left(\Phi(M_{ij}\bz)-\Phi(\bz)\right) \nonumber \\ 
& & +\frac{1}{2|\bz|^{2+\lambda}}\sum_i
\int_0^{z_i}dy F(y,z_i-y)\left(\Phi(S_i^{(y)}\bz)-\Phi(\bz)\right),        \label{2.1} 
\end{eqnarray} 
where $M_{ij}\bz\in\Omega$ 
(resp. $S_i^{(y)}\bz\in\Omega$) 
is obtained from a sequence $\bz=(z_k)$ by merging  
$z_i$ and $z_j$ into $z_i+z_j$ 
(resp. by splitting $z_i$ into $y$ and $z_i-y$) 
and then by reordering. Noting that 
$M_{ij}\bz=\bz$ (resp. $S_i^{(y)}\bz=\bz$) 
whenever $z_iz_j=0$ (resp. $z_i=0$), we regard 
the sum $\sum_{i\ne j}$ (resp. $\sum_i$) in (\ref{2.1}) 
as the sum taken over 
$i\ne j$ (resp $i$) such that $z_iz_j>0$ (resp. $z_i>0$). 
We adopt such convention 
for the same kind of expression throughout. 
The boundedness of $L$ is seen 
easily from (H1) and (H2) in view of 
alternative expression for (\ref{2.1}) 
\begin{eqnarray} 
\wt{L}\Phi(\bz)
& = & \frac{1}{2}\sum_{i\ne j}
\frac{z_i}{|\bz|}\cdot\frac{z_j}{|\bz|} 
\left(\frac{z_i+z_j}{|\bz|}\right)^{\lambda} 
\Hh\left(\frac{z_i}{z_i+z_j},\frac{z_j}{z_i+z_j}\right)
\left(\Phi(M_{ij}\bz)-\Phi(\bz)\right) \nonumber \\ 
& & +\frac{1}{2}\sum_i 
\left(\frac{z_i}{|\bz|}\right)^{2} 
\left(\frac{z_i}{|\bz|}\right)^{\lambda} 
\int_0^{1}du \Hc(u,1-u)\left(\Phi(S_i^{u}\bz)-\Phi(\bz)\right),        \label{2.2} 
\end{eqnarray} 
where $S_i^{u}\bz:=S_i^{(uz_i)}\bz$. 
Indeed, it follows that 
$\Vert \wt{L}\Phi\Vert_{\infty}\le(\Ch\vee\Cc)
\Vert \Phi\Vert_{\infty}$. 
Here and in what follows we use the notation 
$a \vee b := \max\{a,b\}$ and $a \wedge b := \min\{a,b\}$.  
The standard argument (e.g., \S 2 of Chapter 4 in \cite{EK}) 
shows that $ \wt{L}$ generates 
a continuous-time Markov jump process 
$\{\wt{Z}(t)=(\wt{Z}_i(t))_{i=1}^{\infty}:~t\ge 0\}$, say, on $\Omega$. 
It is obvious that $|\wt{Z}(t)|=|\wt{Z}(0)|$ 
for all $t\ge 0$ a.s. 
Similarly, $L_1$, the restriction of $\wt{L}$ on $B(\Omega_1)$, namely, 
\begin{eqnarray} 
L_1\Phi(\bx) 
& = & 
\frac{1}{2}\sum_{i\ne j}K(x_i,x_j) 
\left(\Phi(M_{ij}\bx)-\Phi(\bx)\right) \nonumber \\ 
& & +\frac{1}{2}\sum_i 
\int_0^{x_i}dy F(y,x_i-y)\left(\Phi(S_i^{(y)}\bx)-\Phi(\bx)\right) \nonumber \\ 
& = & 
\frac{1}{2}\sum_{i\ne j} x_ix_j \Hh(x_i,x_j) 
\left(\Phi(M_{ij}\bx)-\Phi(\bx)\right) \nonumber \\ 
& & +\frac{1}{2}\sum_i x_i^{2+\lambda} 
\int_0^{1}du \Hc(u,1-u)\left(\Phi(S_i^{u}\bx)-\Phi(\bx)\right)        \label{2.3} 
\end{eqnarray} 
generates a continuous-time Markov process 
$\{X(t)=(X_i(t))_{i=1}^{\infty}:~t\ge 0\}$ on $\Omega_1$. 
In the case where $\Hh\equiv 1$ and $\lambda=0$, 
(\ref{2.3}) can be thought of as the generator of  
continuous-time version of a Markov chain studied in 
\cite{MWZZ}. Moreover, the operator 
(\ref{2.3}) is a special case (more specifically, 
the binary fragmentation case) of 
the generator considered in \cite{C2}, although 
in order for the model to be defined also on $\Omega$ 
we need homogeneity of $K$, whereas in \cite{C2} 
certain continuity of $K$ is imposed. 

For each $a>0$, define the dilation map 
$D_a:\Omega \to \Omega $ by $D_a(\bz)=a\bz:=(az_i)$. 
The relationship between $Z(t)$ and $X(t)$ 
described in the following lemma is fundamental. 
\begin{lm}
(i) Suppose that a process $\{\wt{Z}(t):~t\ge 0\}$ 
generated by $\wt{L}$ is given. Then  
the $\Omega_1$-valued process 
$\{X(t)=(X_i(t))_{i=1}^{\infty}:~t\ge 0\}$ 
defined by 
$X_i(t)=\wt{Z}_i(t)/|\wt{Z}(t)|=\wt{Z}_i(t)/|\wt{Z}(0)|$ 
is a process generated by $L_1$. \\
(ii) Suppose that a $(0,\infty)$-valued random variable $V$ 
and a process $\{(X_i(t))_{i=1}^{\infty}:~t\ge 0\}$ 
generated by $L_1$ are mutually independent. Then  
the $\Omega$-valued process 
$\{\wt{Z}(t)=(\wt{Z}_i(t))_{i=1}^{\infty}:~t\ge 0\}$ defined 
by $\wt{Z}_i(t)=VX_i(t)$ is a process generated by $\wt{L}$. 
\end{lm}
{\it Proof.}~ 
(i) Take $\Phi\in B(\Omega)$ arbitrarily  
and denote by $E[\ \cdot \ |Z(0)=\bz]$ 
the expectation with respect to the process 
generated by $L$ starting from $\bz\in\Omega$. 
Since $v^{-1}(M_{ij}\bz)=M_{ij}(v^{-1}\bz)$ and 
$|\bz|^{-1}(S_{i}^{(uz_i)}\bz)=S_{i}^u(|\bz|^{-1}\bz)$,  
we see from (\ref{2.2}) 
\[
\wt{L}(\Phi\circ D_{1/v})(\bz)=(L_1\Phi)(v^{-1}\bz), 
\quad \bz\in \Omega, v=|\bz|. 
\] 
Hence, for any $t>0$, by Fubini's theorem 
\begin{eqnarray*} 
\lefteqn{E\left[\Phi(X(t))\right]
-E\left[\Phi(X(0))\right]}               \\ 
& = &  
E\left[\Phi(|\wt{Z}(0)|^{-1}\wt{Z}(t))\right]
-E\left[\Phi(|\wt{Z}(0)|^{-1}\wt{Z}(0))\right]   \\ 
& = & 
\int_{\Omega}P(\wt{Z}(0)\in d\bz)
\left\{E\left[(\Phi\circ D_{1/|\bz|})
(\wt{Z}(t))|\wt{Z}(0)=\bz\right]
-(\Phi\circ D_{1/|\bz|})(\bz)\right\}     \\ 
& = & 
\int_{\Omega}P(\wt{Z}(0)\in d\bz)
\int_0^tdsE\left[\wt{L}(\Phi\circ D_{1/|\bz|})(\wt{Z}(s))| 
\wt{Z}(0)=\bz\right]   \\ 
& = & 
\int_0^tds\int_{\Omega}P(\wt{Z}(0)\in d\bz)
E\left[L_1\Phi(|\bz|^{-1}\wt{Z}(s))|\wt{Z}(0)=\bz\right]  \\ 
& = & 
\int_0^tds E\left[L_1\Phi(|\wt{Z}(0)|^{-1}\wt{Z}(s))\right] 
\ = \ 
\int_0^tds E\left[L_1\Phi(X(s))\right].   
\end{eqnarray*} 
This proves the first assertion. \\ 
(ii) Based on the relation 
$L_1(\Phi\circ D_v)(\bx)=\wt{L}\Phi(v\bx)$ 
for $\bx\in \Omega_1$ and $v>0$, the proof of the second assertion is very similar to that for (i). So we omit it. \qed 

\medskip 

We call $V$ in Lemma 2.1 a lifting variable. 
Roughly speaking, lifting a process on $\Omega_1$ generated by $L_1$ 
yields a process on  $\Omega$ generated by $\wt{L}$. 
We need to consider an unbounded operator 
${L}\Phi(\bz)=|\bz|^{2+\lambda} \wt{L}\Phi(\bz)$ or 
\begin{eqnarray}  
{L}\Phi(\bz)
& = & 
\frac{1}{2}\sum_{i\ne j}K(z_i,z_j)
\left(\Phi(M_{ij}\bz)-\Phi(\bz)\right) 
                 \nonumber              \\ 
&  &  
+\frac{1}{2}\sum_i\int_0^{z_i}dy F(y,z_i-y) 
\left(\Phi(S_i^{(y)}\bz)-\Phi(\bz)\right).             \label{2.4} 
\end{eqnarray} 
The corresponding process 
$\{Z(t):~t\ge 0\}$ on  $\Omega$ can be obtained from 
a process $\{\wt{Z}(t):~t\ge 0\}$ generated by $\wt{L}$ 
with the same initial law by a random time-change 
\[ 
{Z}(t):=\wt{Z}(|\wt{Z}(0)|^{2+\lambda}t). 
\] 
This can be shown by general theory of Markov processes, 
e.g., Theorem 1.3 in Chapter 6 of \cite{EK}, or 
more directly, by making the following observation: 
for any $\Phi\in B(\Omega)$ such that 
$L\Phi\in B(\Omega)$, 
by the optional sampling theorem  
\begin{eqnarray*} 
\lefteqn{
\Phi(Z(t))-\Phi(Z(0))
-\int_0^tdu L\Phi({Z}(u)) }                 \\ 
& = & 
\Phi(\wt{Z}(|\wt{Z}(0)|^{2+\lambda}t))-\Phi(\wt{Z}(0)) 
-\int_0^tdu|\wt{Z}(0)|^{2+\lambda}
\wt{L}\Phi(\wt{Z}(|\wt{Z}(0)|^{2+\lambda}u))  \\ 
& = & 
\Phi(\wt{Z}(|\wt{Z}(0)|^{2+\lambda}t))-\Phi(\wt{Z}(0)) 
- \int_0^{|\wt{Z}(0)|^{2+\lambda}t}ds\wt{L}\Phi(\wt{Z}(s))
\end{eqnarray*} 
is a martingale. 

\subsection{Reformulation in terms of point processes} 
We proceed to reformulate 
the above-mentioned processes as Markov processes 
taking values in a space of point-configurations on $(0,\infty)$. 
To discuss in the setting of point processes, 
we need the following notation. 
Set $\Zp=\{0,1,2,\ldots\}$ and let $\cN$ be 
the set of $\Zp$-valued Radon measures  
on $(0,\infty)$. Each element $\eta$ of $\cN$ 
is regarded as a counting measure associated with 
a locally finite point-configuration on $(0,\infty)$ with multiplicity. 
We equip $\cN$ with the vague topology and 
use the notation $|\eta|:=\int_{(0,\infty)}v\eta(dv)$ 
for $\eta\in\cN$ and 
$\lg f,\nu\rg:=\int_{(0,\infty)}f(v)\nu(dv)$ 
for a measure $\nu$ and a Borel function 
$f$ on $(0,\infty)$. 
Denote by $B_+(S)$ the set of nonnegative 
bounded Borel measurable functions 
on a topological space $S$. For simplicity, 
we set $B_+=B_+((0,\infty))$ and 
use the notation $B_+^k$ instead of $B_+((0,\infty)^k)$ 
for $k=2,3,\ldots$.  
As mentioned roughly in Introduction, 
the subsequent argument is based on 
the one-to-one correspondence between 
$\bz=(z_i)\in \Omega$ and 
\[ 
\Xi(\bz):=\sum_i\one_{\{z_i>0\}}\delta_{z_i}
\in \{\eta\in\cN:~|\eta|<\infty\}=:\cN_1. 
\] 
Clearly, if $\eta=\Xi(\bz)$, then $|\eta|=|\bz|$ and 
$\eta([\epsilon,\infty))\le |\bz|/\epsilon$ 
for $\epsilon>0$. Note that the map 
$\Xi:\Omega\to \cN_1$ 
is bi-measurable. It follows that 
\be 
\Xi(M_{ij}\bz)-\Xi(\bz)
=\delta_{z_i+z_j}-\delta_{z_i}-\delta_{z_j} 
\quad \mbox{if} \quad z_i,z_j>0,  
                                        \label{2.5} 
\ee 
and 
\be 
\Xi(S_{i}^{(y)}\bz)-\Xi(\bz)
=\delta_{y}+\delta_{z_i-y}-\delta_{z_i} 
\quad \mbox{if} \quad z_i>y>0.  
                                        \label{2.6} 
\ee 
Thus, employing $\Xi(\bz)$ rather than $\bz$ itself 
enables us to avoid an unnecessary complication arising 
from reordering of the sequence. 

Besides, owing to the map $\Xi$,  
the arguments below 
make use of some effective tools 
in theory of point processes, 
such as correlation measures and 
(reduced) Palm distributions. 
(See e.g. \S 13.1 of \cite{DV2} for general accounts.) 
We shall give their definitions 
for a locally finite point process 
$\xi=\sum\delta_{Z_i}$ on $(0,\infty)$. 
In what follows, the domain of integration will be 
suppresed as long as it is $(0,\infty)^k$ for some $k\in\N$, 
which should be clear from context. 
The first correlation measure $q_1$ is the mean measure of $\xi$, 
and for $k=2,3,\ldots$ the $k$th correlation measure $q_k$ 
is the mean measure of the modified product counting measure  
\be 
\xi^{[k]}:=\sum_{i_1,\ldots,i_k(\ne)}
\delta_{(Z_{i_1},\ldots,Z_{i_k})}   \label{2.7}  
\ee
on $(0,\infty)^k$, where $\sum_{i_1,\ldots,i_k(\ne)}$ indicates that 
the sum is taken over $k$-tuplets $(i_1,\ldots,i_k)$ 
such that $i_l\ne i_m$ whenever $l\ne m$.  
The entire system $\{q_1,q_2,\ldots\}$ of 
correlation measures determines uniquely the law of $\xi$. 
(The identity (\ref{2.13}) below is 
the structure underlying this fact.) 
The density of $q_k$ is called 
the $k$th correlation function of $\xi$ if it exists. 
Furthermore, letting 
$\{P_{z_1,\ldots,z_k}:~z_1,\ldots,z_k\in (0,\infty)\}$ 
be a family of Borel probability measures on $\cN$, 
we call $P_{z_1,\ldots,z_k}$ 
the $k$th-order reduced Palm distribution of $\xi$ 
at $(z_1,\ldots,z_k)$ if 
for any $G\in B_+((0,\infty)^k\times \cN)$ 
\begin{eqnarray} 
\lefteqn{E\left[\sum_{i_1,\ldots,i_k(\ne)}
G(Z_{i_1},\ldots,Z_{i_k},\xi)\right]
= E\left[\int \xi^{[k]}(dz_1\cdots dz_k)
G(z_1,\ldots,z_k, \xi)\right]}    \nonumber  \\ 
& = & 
\int q_k(dz_1\cdots dz_k) E_{z_1,\ldots,z_k}
\left[G(z_1,\ldots,z_k,\eta+\delta_{z_1}+\cdots+\delta_{z_k})\right].  
                                        \label{2.8}  
\end{eqnarray} 
Here and throughout, $E_{z_1,\ldots,z_k}$ denotes the expectation 
with respect to $P_{z_1,\ldots,z_k}$, so that 
\[ 
E_{z_1,\ldots,z_k}\left[\Psi(\eta)\right]
=
\int_{\cN}\Psi(\eta)P_{z_1,\ldots,z_k}(d\eta), 
\quad \Psi\in B_+(\cN). 
\] 
(Rigorously speaking, certain measurability 
of $P_{z_1,\ldots,z_k}$ in $(z_1,\ldots,z_k)$ is 
required just as in the definition of 
regular conditional distributions. 
See Proposition 13.1.IV of \cite{DV2} for the case $k=1$.) 
An intuitive interpretation for $P_{z_1,\ldots,z_k}$ 
is the conditional law of 
$\xi-\delta_{z_1}-\cdots-\delta_{z_k}$ 
given $\xi(\{z_1\})\cdots\xi(\{z_k\})>0$. 
We call an equality of the type (\ref{2.8}) 
the $k$th-order Plam formula for $\xi$. 

To state the main result of this section, 
let us recall some known results on 
Poisson-Dirichlet point process, 
which is by definition the point-configuration associated with 
a random element of $\Omega_1$ 
distributed according to the Poisson-Dirichlet distribution. 
For each $\theta>0$, 
let ${\rm PD}(\theta)$ denote the Poisson-Dirichlet 
distribution with parameter $\theta$. 
(See e.g. \cite{King}, \cite{B}, \cite{Feng} 
for the definition.) Suppose that $X=(X_i)$ is 
a random element of $\Omega_1$ 
whose law is ${\rm PD}(\theta)$ 
and consider the associated point process 
$\xi^{(\theta)}=\sum\delta_{X_i}$ on the interval $(0,1)$, 
which we call simply the ${\rm PD}(\theta)$ process. 
It was shown in \cite{W} that 
the $k$th correlation function of $\xi^{(\theta)}$ 
takes the form 
\be 
(x_1,\ldots,x_{k}) \mapsto 
\frac{\theta^k}{x_1\cdots x_k} 
\left(1-
\sum_{j=1}^{k}x_{j}\right)^{\theta-1}
\one_{\Delta_k^{\circ}}(x_1,\ldots,x_{k}),  
                                          \label{2.9}  
\ee 
where 
\[ 
\Delta_k^{\circ}=\left\{(x_1,\ldots,x_{k}): 
x_1,\ldots,x_k>0, \ 
x_1+\cdots+x_k<1 \right\}.    
\] 
By a special case of Corollary 1 of \cite{Le} 
the $k$th-order reduced Palm distribution 
of $\xi^{(\theta)}$ 
at $(x_1,\ldots,x_{k})\in\Delta_k^{\circ}$ coincides with 
the law of 
$\xi^{(\theta)}\circ D_{1-x_1-\cdots-x_{k}}^{-1}$. 
Moreover, such a self-similar property of Palm distributions 
characterizes the family $\{{\rm PD}(\theta):~\theta>0\}$ 
as was shown in Theorem 2 (iii) of \cite{Le}. 

\medskip 

The next lemma shows how the aforementioned 
notions are well adapted 
for dealing with our generators. Denote by $B_c$ the totality of 
bounded Borel functions on $(0,\infty)$ with 
compact support and let $B_{+,c}=B_+\cap B_c$. 
The support and the sup norm of a function $f$ are 
denoted by ${\rm supp}(f)$ and $\Vert f\Vert_{\infty}$, respectively. 
We now give a class of functions on $\Omega$ 
for which our generators act in a tractable manner. Set 
\[ 
\widetilde{B}_{+}=\{\phi\in B_{+}:
~\phi-1\in B_c, \Vert \phi-1\Vert_{\infty}<1\}. 
\] 
For each $\phi\in \widetilde{B}_{+}$ 
we can define a function $\Pi_{\phi}$ on $\Omega$ by 
$\Pi_{\phi}(\bz)=\prod_{i:z_i>0}\phi(z_i)$,   
noting that the right side is a finite product. 
By abuse of notation, we also write 
$\Pi_{\phi}(\xi)$ for $\Pi_{\phi}(\bz)$ 
when $\xi=\Xi(\bz)$. Thus 
\be 
\Pi_{\phi}(\xi) 
= \exp\left(\sum_i\one_{\{z_i>0\}}\log \phi(z_i) \right)
= e^{\lg \log \phi,\xi\rg}.                                     \label{2.10} 
\ee
An important remark is that the class 
$\{\Pi_{\phi}:~\phi\in\widetilde{B}_{+}\}$ is 
measure-determining on $\Omega$ 
because it includes all functions of the form 
$\bz \mapsto \exp(-\lg f,\Xi(\bz)\rg)$ with $f\in B_{+,c}$. 
Given a function on $\Omega$, we regard it also as 
a function on $\Omega_{\le R}:=\{\bz\in\Omega:~|\bz|\le R\}$ 
for any $R>0$. 
It is clear that $L$ restricted on $B(\Omega_{\le R})$ is bounded.  
\begin{lm}
Let $\phi\in \widetilde{B}_{+}$ and 
set $f=\phi-1$. Then for any $\bz\in\Omega$ 
\be 
0\le \Pi_{\phi}(\bz) 
\le 
(1-\Vert f\Vert_{\infty})^{-\xi({\rm supp}(f))},   \label{2.11} 
\ee 
where $\xi=\Xi(\bz)$. Moreover, 
$\Pi_{\phi}\in B(\Omega_{\le R})$ for any $R>0$ and 
\begin{eqnarray} 
L\Pi_{\phi}(\bz)  
& = &  
\frac{1}{2}\sum_{i\ne j}K(z_i,z_j)
\left[\phi(z_i+z_j)-\phi(z_i)\phi(z_j)\right]
\prod_{k\ne i,j}\phi(z_k)                \label{2.12}      \\ 
&  & 
+\frac{1}{2}\sum_i \int_0^{z_i}dy F(y,z_i-y)
\left[\phi(y)\phi(z_i-y)-\phi(z_i)\right]
\prod_{k\ne i}\phi(z_k)                  \nonumber  \\ 
& = & 
\frac{1}{2}
\int\xi^{[2]}(dv_1dv_2)K(v_1,v_2)
\left[\phi(v_1+v_2)-\phi(v_1)\phi(v_2)\right]
\Pi_{\phi}(\xi-\delta_{v_1}-\delta_{v_2})  \nonumber     \\ 
&  & 
+\frac{1}{2}\int\xi(dv)\int_0^vdy F(y,v-y)
\left[\phi(y)\phi(v-y)-\phi(v)\right]
\Pi_{\phi}(\xi-\delta_{v}).                \nonumber  
\end{eqnarray} 
\end{lm}
{\it Proof.}~
Let $\phi(0)=1$ and  $f(0)=0$ by convention. 
To prove (\ref{2.11}) observe that 
\be 
\prod_{i}\phi(z_i) 
= 
\prod_{i}(1+f(z_i)) 
= 
1+\sum_{k=1}^{\infty}\frac{1}{k!}
\sum_{i_1,\ldots,i_k(\ne)}f(z_{i_1})\cdots f(z_{i_k}).   \label{2.13} 
\ee 
So, using the notation 
${\alpha \choose k}
=\alpha(\alpha-1)\cdots(\alpha-k+1)/k!$, we get 
\begin{eqnarray*} 
0 \le \Pi_{\phi}(\bz) 
& \le & 
1+\sum_{k=1}^{\infty}\frac{\Vert f\Vert_{\infty}^k}{k!}
\sum_{i_1,\ldots,i_k(\ne)}
\one_{\{z_{i_1},\ldots,z_{i_k}\in{\rm supp}(f)\}} \\ 
& = & 
1+\sum_{k=1}^{\infty}\Vert f\Vert_{\infty}^k 
{\xi({\rm supp}(f)) \choose k}  \\ 
& = & 
(1-\Vert f\Vert_{\infty})^{-\xi({\rm supp}(f))}, 
\end{eqnarray*} 
and thus (\ref{2.11}) follows. 
Putting 
$\epsilon=\inf{\rm supp}(f)>0$, we get 
\[ 
\xi({\rm supp}(f))\le \xi([\epsilon,\infty)) 
\le \epsilon^{-1}|\bz|. 
\] 
This combined with (\ref{2.11}) implies that $\Pi_{\phi}$ 
is bounded on $\Omega_{\le R}$. (\ref{2.12}) is verified by 
direct calculations with the help of 
(\ref{2.5}) and (\ref{2.6}). \qed 

\medskip 

\subsection{Reversible cases} 
We demonstrate the power of 
point process calculus by proving 
an extension of the reversibility result 
due to Mayer-Wolf et al \cite {MWZZ}. 
(As for the stationarity result, a proof 
based on the underlying Poisson process 
can be found in \S 7.3 of \cite{Gold}.)  
To this end, we recall that 
the $k$th correlation function $q_{k}$ of the 
${\rm PD}(\theta)$ process $\xi^{(\theta)}$ 
is given in (\ref{2.9}) 
and that the $k$th-order reduced Palm distribution 
$P_{x_1,\ldots,x_k}$ of $\xi^{(\theta)}$ at 
$(x_1,\ldots,x_k)\in \Delta_k^{\circ}$ is 
the law of $\xi^{(\theta)}\circ D_{1-x_1-\cdots-x_k}^{-1}$. 
In particular, for any $(x_1,x_2)\in\Delta_2^{\circ}$ 
\be 
x_1x_2q_2(x_1,x_2)=\theta (x_1+x_2)q_1(x_1+x_2)  \label{2.14}  
\ee 
and 
\be 
P_{x_1,x_2}=P_{x_1+x_2}.                         \label{2.15}  
\ee 
As will be shown in the next theorem, 
these identities are responsible 
for the reversibility of ${\rm PD}(\theta)$  
with respect to processes associated with bounded operators 
on $B(\Omega_1)$ of the form  
\begin{eqnarray*} 
L_1^{(Q,\theta)}\Phi(\bx) 
& =  & 
\frac{1}{2}\sum_{i\ne j} x_ix_j Q(x_i,x_j) 
\left(\Phi(M_{ij}\bx)-\Phi(\bx)\right) \nonumber \\ 
& &  +\frac{\theta}{2}\sum_i x_i^{2} 
\int_0^{1}du Q(ux_i,(1-u)x_i) 
\left(\Phi(S_i^{u}\bx)-\Phi(\bx)\right),  
\end{eqnarray*} 
where $Q$ is any nonzero bounded, 
symmetric nonnegative function on 
$\{(x,y) |~x,y> 0, x+y \le 1\}$. 
(It should be noted that Theorem 12 in \cite{MWZZ}   
proved essentially the symmetry of $L_1^{(Q,\theta)}$ with 
$Q\equiv \mbox{const.}$ with respect to 
${\rm PD}(\theta)$.) More generally, we consider 
\begin{eqnarray} 
L_1^{\sharp}\Phi(\bx) 
& = & 
\frac{1}{2}\sum_{i\ne j}K_1(x_i,x_j) 
\left(\Phi(M_{ij}\bx)-\Phi(\bx)\right) \nonumber \\ 
& & +\frac{1}{2}\sum_i x_i
\int_0^{1}du F_1(ux_i,(1-u)x_i) 
\left(\Phi(S_i^{u}\bx)-\Phi(\bx)\right)  
                             \label{2.16}  
\end{eqnarray} 
with $K_1$ and $F_1$ being symmetric nonnegative 
functions on $\{(x,y) |~x,y>0, x+y \le 1\}$ such that 
$(x_i)\mapsto 
\sum_{i\ne j}K_1(x_i,x_j)\one_{\{x_ix_j>0\}}$ and 
$(x_i)\mapsto \sum x_i \int_0^{1}du F_1(ux_i,(1-u)x_i)$ 
are bounded functions on $\Omega_1$. We may and do 
suppose that $K_1(x,y)=0$ whenever $xy=0$. 
\begin{th}
(i) Let  $X=(X_i)_{i=1}^{\infty}$ be a random element 
of $\Omega_1$ and suppose that 
the first and second correlation functions 
$q_1$ and $q_2$ of $\xi:=\Xi(X)$ exist and satisfy 
\be 
K_1(x_1,x_2)q_2(x_1,x_2)P_{x_1,x_2} 
=F_1(x_1,x_2)q_1(x_1+x_2) P_{x_1+x_2}, 
                                        \label{2.17} 
\ee 
$\mbox{a.e.-}(x_1,x_2)\in \Delta_2^{o}$. 
Here, $P_{x_1,x_2}$ and $P_{x_1+x_2}$ are 
the reduced Palm distributions of $\xi$ and 
the above equality is understood as the one 
between two measures on $\cN$. Then 
the process generated by $L_1^{\sharp}$ is 
reversible with respect to the law of $X$. 
\\   
(ii) ${\rm PD}(\theta)$ is a reversible distribution 
of the process generated by $L_1^{(Q,\theta)}$. 
\end{th}
{\it Proof.}~ 
(i) Since $L_1^{\sharp}$ is bounded, 
we only have to check the symmetry 
\[ 
E\left[\Phi(X)L_1^{\sharp}\Psi(X)\right]
=E\left[\Psi(X)L_1^{\sharp}\Phi(X)\right]
\] 
for any $\Phi,\Psi\in B(\Omega_1)$. 
We will prove stronger equalities 
\begin{eqnarray*}
\lefteqn{\lg\lg\Phi,\Psi\rg\rg_{\rm coag} 
\ := \  
E\left[\Phi(X)\sum_{i\ne j}K_1(X_i,X_j)\Psi(M_{ij}X)\right]} 
                                                         \\    
& = &  
E\left[\sum_{i}X_i\int_0^1duF_1(uX_i,(1-u)X_i) 
\Phi(S_{i}^uX)\Psi(X)\right] 
\ =: \ 
{}_{\rm frag}\lg\lg\Phi,\Psi\rg\rg,       
\end{eqnarray*} 
which are regarded as a sort of 
coagulation-fragmentation duality. 
Furthermore, it suffices to verify for functions 
$\Phi=\Pi_{\phi}$ and $\Psi=\Pi_{\psi}$ 
with $\phi,\psi\in \widetilde{B}_+$. Thanks to the 
first-order and second-order Palm formulae for  $\xi=\sum\one_{\{X_i>0\}}\delta_{X_i}$, 
similar calculations to (\ref{2.12}) show that 
\begin{eqnarray*} 
\lefteqn{
\lg\lg\Pi_{\phi},\Pi_{\psi}\rg\rg_{\rm coag}}       \\  
& = & 
E\left[\sum_{i\ne j}K_1(X_i,X_j)\phi(X_i)\phi(X_j)
\psi(X_i+X_j)\prod_{k\ne i,j}(\phi(X_k)\psi(X_k))\right]   \\ 
& = & 
E\left[\int\xi^{[2]}(dx_1dx_2)K_1(x_1,x_2)\phi(x_1)\phi(x_2)
\psi(x_1+x_2)\Pi_{\phi\psi}(\xi-\delta_{x_1}-\delta_{x_2})\right]   \\ 
& = & 
\int_{\Delta_{2}^{\circ}}q_2(x_1,x_2)
K_1(x_1,x_2)\phi(x_1)\phi(x_2)\psi(x_1+x_2)
E_{x_1,x_2}\left[\Pi_{\phi\psi}(\eta)\right]dx_1dx_2 
\end{eqnarray*} 
and 
\begin{eqnarray*} 
\lefteqn{
{}_{\rm frag}\lg\lg\Pi_{\phi},\Pi_{\psi}\rg\rg}       \\  
& = & 
E\left[\sum_{i}X_i\int_0^1du F_1(uX_i,(1-u)X_i)
\phi(uX_i)\phi((1-u)X_i)\psi(X_i)\prod_{k\ne i}(\phi(X_k)\psi(X_k))\right] \\ 
& = & 
E\left[\int\xi(dv)v\int_0^1du F_1(uv,(1-u)v)
\phi(uv)\phi((1-u)v)\psi(v)
\Pi_{\phi\psi}(\xi-\delta_{v})\right] \\ 
& = & 
\int_0^1q_1(v)v\int_0^1F_1(uv,(1-u)v)
\phi(uv)\phi((1-u)v)\psi(v)
E_{v}\left[\Pi_{\phi\psi}(\eta)\right]dudv. 
\end{eqnarray*} 
By virtue of (\ref{2.17}) we obtain the desired equality 
$\lg\lg\Phi,\Psi\rg\rg_{\rm coag}
={}_{\rm frag}\lg\lg\Phi,\Psi\rg\rg$ 
after the change of variables $uv=:x_1,(1-u)v=:x_2$ \\ 
(ii) This assertion is immediate 
by noting that (\ref{2.17}) with 
\[ 
K_1(x_1,x_2)=x_1x_2Q(x_1,x_2) \quad \mbox{and} 
\quad 
F_1(x_1,x_2)=\theta(x_1+x_2)Q(x_1,x_2) 
\] 
is valid for the ${\rm PD}(\theta)$ process  
because of (\ref{2.14}) and (\ref{2.15}).  
The proof of Theorem 2.3 is complete. 
\qed 

\medskip 

\noindent 
{\it Remarks.}~
(i) It would be interesting to investigate the class 
of (nonnegative unbounded) functions $Q$ for which 
the operator $L_1^{(Q,\theta)}$ defines  
a Markov process on $\Omega_1$. 
For example, if taking $Q(x,y)=(xy)^{-1}$ 
is allowed in that sense, the operator 
\begin{eqnarray*} 
L_1^{(Q,\theta)}\Phi(\bx) 
= 
\frac{1}{2}\sum_{i\ne j} 
\left(\Phi(M_{ij}\bx)-\Phi(\bx)\right) 
  +\frac{\theta}{2}\sum_i \int_0^{1}\frac{du}{u(1-u)} 
\left(\Phi(S_i^{u}\bx)-\Phi(\bx)\right)  
\end{eqnarray*}  
could deserve further exploration 
in connection with e.g. `asymptotic frequency' 
of some exchangeable fragmentation-coalescence 
process studied in \cite{Be}. It is pointed out that 
at least the unbounded coagulation operator 
in the above can be treated 
within the fame work of \cite{C2}.  \\ 
(ii) Alternative direction of generalization of 
the processes reversible with respect to 
${\rm PD}(\theta)$ 
is based on an obvious generalization of 
(\ref{2.14}) and (\ref{2.15}), i.e., 
\[ 
x_1\cdots x_{k+1}q_{k+1}(x_1,\ldots,x_{k+1}) 
=\theta^k (x_1+\cdots +x_{k+1})
q_1(x_1+\cdots +x_{k+1}) 
\] 
and 
\[ 
P_{x_1,\ldots,x_{k+1}}=P_{x_1+\cdots +x_{k+1}}, 
\] 
in which $k\in\N$ is arbitrary and 
$(x_1,\ldots,x_{k+1})\in \Delta_{k+1}^{o}$. 
The corresponding 
process on $\Omega_1$ 
incorporates multiple-coagulation 
and multiple-fragmentation. (cf. 
the transition kernel (\ref{2.5}) in \cite{EW} 
or  Example 1.8 in \cite{Kol}. 
See also \cite{C2} for more general scheme 
for the multiple-fragmentation.) 
One of the simplest examples of the generator 
of such a process is 
\begin{eqnarray*} 
\lefteqn{
L_{1,k}^{(\theta)}\Pi_{\phi}(\bx)}            \\  
& := &  
\sum_{i_1,\ldots,i_{k+1}(\ne)} x_{i_1}\cdots x_{i_{k+1}} 
\left[\phi(x_{i_1}+\cdots+x_{i_{k+1}})
-\phi(x_{i_1})\cdots\phi(x_{i_{k+1}})\right]
\prod_{j\ne i_1,\ldots,i_{k+1}}\phi(x_j)        \\ 
&  & 
+\theta^k\sum_i x_i^{k+1}\int_{\Delta_k}du_1\cdots du_k 
\left[\phi(u_1x_i)\cdots\phi(u_{k}x_i)
\phi((1-|u|)x_i)-\phi(x_i)\right]
\prod_{j\ne i}\phi(x_j)                      \\ 
& = & 
\int\xi^{[k+1]}(dv_1\cdots dv_{k+1})v_1\cdots v_{k+1} \\ 
& & 
\hspace*{10mm} \times
\left[\phi(v_1+\cdots+v_{k+1})
-\phi(v_1)\cdots\phi(v_{k+1})\right]          
\Pi_{\phi}(\xi-\delta_{v_1}-\cdots-\delta_{v_{k+1}})  \\ 
&  & 
+\theta^k\int\xi(dv)v^{k+1}
\int_{\Delta_k}du_1\cdots du_k              \\ 
& & 
\hspace*{10mm} \times 
\left[\phi(u_1v)\cdots\phi(u_{k}v)
\phi((1-|u|)v)-\phi(v)\right]
\Pi_{\phi}(\xi-\delta_{v}),          
\end{eqnarray*} 
where $|u|=u_1+\cdots +u_{k}$ and $\xi=\Xi(\bx)$. 
Clearly $L_{1,1}^{(\theta)}=2L_{1}^{(1,\theta)}$. 
The calculations in the proof of Theorem 2.3 
is easily modified to prove that 
${\rm PD}(\theta)$ is still a reversible distribution 
of the process generated by $L_{1,k}^{(\theta)}$. 
The details are left to the reader.  \\ 
(iii) The spectral gap of a suitable extension 
$\overline{L_1^{(Q,\theta)}}$of $L_1^{(Q,\theta)}$ vanishes. 
Indeed, letting $\Psi_{\delta}(\bz)=\sum_iz_i^{\delta}$ for 
$\bz=(z_i)\in \Omega$ and $\delta>0$, 
we see, with the help of Lemma 6.4 in \cite{H09}, 
that $\Psi_{\delta}$ is square integrable 
with respect to ${\rm PD}(\theta)$ and that 
its variance ${\rm var}(\Psi_{\delta})$ is given by 
\[ 
{\rm var}(\Psi_{\delta})
\ = \   
\frac{\Gamma(\theta+1)\Gamma(2\delta)}
{\Gamma(\theta+2\delta)}   
+\Gamma(\delta)^2 
\left(\frac{\theta\Gamma(\theta+1)}{\Gamma(\theta+2\delta)}
-\frac{\Gamma(\theta+1)^2}{\Gamma(\theta+\delta)^2}\right) 
\ =: \   
\chi_1(\delta)+\chi_2(\delta),     
\] 
where $\Gamma(\cdot)$ is Gamma function. 
As $\delta \downarrow 0$, 
${\rm var}(\Psi_{\delta})\to \infty$ since 
$\chi_1(\delta)\to \infty$ and 
\begin{eqnarray*} 
\chi_2(\delta)& =&  
\frac{\theta\Gamma(\theta+1)\Gamma(\delta+1)^2}
{\Gamma(\theta+2\delta)\Gamma(\theta+\delta)^2}   
\cdot 
\frac{\Gamma(\theta+\delta)^2
-\Gamma(\theta)\Gamma(\theta+2\delta)}{\delta^2} \\ 
& \to&  
\frac{\theta^2}{\Gamma(\theta)^2}
\left(\Gamma'(\theta)^2-\Gamma(\theta)
\Gamma''(\theta)\right) 
\quad  (\mbox{by l'Hospital's rule}). 
\end{eqnarray*} 
As for Dirichlet form 
$\cE(\Psi_{\delta}): = 
E\left[\Psi_{\delta}(X)(-\overline{L_1^{(Q,\theta)}})
\Psi_{\delta}(X)\right]$   
in which $X=(X_i)_{i=1}^{\infty}$ is 
${\rm PD}(\theta)$-distributed, 
by the dominated convergence theorem 
\begin{eqnarray*} 
\cE(\Psi_{\delta})  
& = & 
\frac{1}{2}E\left[\overline{L_1^{(Q,\theta)}}(\Psi_{\delta}^2)(X)
-2\Psi_{\delta}(X)\overline{L_1^{(Q,\theta)}}\Psi_{\delta}(X)\right]       \\ 
& =&  
\frac{1}{4} E\left[\sum_{i\ne j}X_iX_jQ(X_i,X_j)
\{\Psi_{\delta}(M_{ij}X)-\Psi_{\delta}(X)\}^2\right] \\ 
&  & 
+ \frac{\theta}{4}  
E\left[\sum_{i}X_i^2\int_0^1du Q(uX_i,(1-u)X_i) 
\{\Psi_{\delta}(S_{i}^uX)-\Psi_{\delta}(X)\}^2\right] \\ 
& =&  
\frac{1}{4} E\left[\sum_{i\ne j}X_iX_jQ(X_i,X_j)
\left\{(X_i+X_j)^{\delta}-X_i^{\delta}-X_j^{\delta}\right\}^2\right] \\ 
&  & 
+ \frac{\theta}{4} 
E\left[\sum_{i}X_i^2\int_0^1du Q(uX_i,(1-u)X_i)  
\left\{(uX_i)^{\delta}+((1-u)X_i)^{\delta}
-X_i^{\delta}\right\}^2\right]                          \\ 
& \to & 
\frac{1}{4} E\left[\sum_{i\ne j}X_iX_jQ(X_i,X_j)\right]
+\frac{\theta}{4} E\left[\sum_{i}X_i^2\int_0^1du 
Q(uX_i,(1-u)X_i)  \right]<\infty    
\end{eqnarray*} 
as $\delta \downarrow 0$. 
(In fact, the second equality in the above needs 
justification. This can be done by considering 
bounded functions 
$\Psi_{\delta, \epsilon}(\bx) 
=\sum x_i^{\delta}\one_{\{x_i\ge \epsilon\}}$ 
on $\Omega_1$, taking the limit as 
$\epsilon \downarrow 0$ and 
applying Lebesgue's convergence theorem.) 
Consequently, 
$\cE(\Psi_{\delta})/{\rm var}(\Psi_{\delta})\to 0$, 
and hence the exponential convergence to equilibrium 
does not hold for the process generated by 
$L_1^{(Q,\theta)}$. It seems, however, 
that the exact speed of convergence is unknown 
even in the case $Q\equiv \mbox{const}$. 

\medskip 

The equality (\ref{2.17}) is thought of as a 
probabilistic counterpart of the detailed balance 
condition (\ref{1.3}). It should be noted that 
(\ref{2.17}) is equivalent to the validity of 
two equalities 
$K_1(x_1,x_2)q_2(x_1,x_2)=F_1(x_1+x_2)q_1(x_1+x_2)$ 
and  $P_{x_1,x_2}=P_{x_1+x_2}$, 
a.e.-$(x_1,x_2)\in \Delta_2^{o}$ 
because of the triviality that 
$P_{x_1,x_2}(\cN)=P_{x_1+x_2}(\cN)=1$. 
The reader may wonder 
whether there is any distribution other 
than Poisson-Dirichlet distributions 
which enjoys the relation (\ref{2.17}) 
for some explicit $K_1$ and $F_1$. 
The following examples are intended to 
give answers to that question by discussing 
certain deformations of ${\rm PD}(\theta)$. 

\medskip 

\noindent 
{\it Examples.}~ 
Let $(X_i)_{i=1}^{\infty}$ be ${\rm PD}(\theta)$-distributed 
and $\phi$ be a nonnegative measurable function 
on $(0,1)$. Assume that  $0<a:=E\left[\prod_{i}\phi(X_i)\right]<\infty$ and 
define a probability measure $\wt{P}$ on $\Omega_1$ by 
\[ 
\wt{P}(\bullet)=a^{-1} 
E\left[\prod_{i}\phi(X_i):~(X_i)_{i=1}^{\infty}\in\bullet \right]. 
\] 
It is not difficult to show that 
the first and second correlation functions 
$\wt{q}_1$ and $\wt{q}_2$ and 
the first-order and second-order reduced 
Palm distributions of $\sum \delta_{X_i}$   
under $\wt{P}$ are given in terms of 
those of the ${\rm PD}(\theta)$ process 
(namely, $q_1, q_2, P_v$ and $P_{x_1,x_2}$) by 
\[ 
\wt{q}_1(v)=a^{-1}\phi(v)q_1(v) 
E_v\left[\prod_{i}\phi(X_i)\right], 
\] 
\[ 
\wt{q}_2(x_1.x_2)
=a^{-1}\phi(x_1)\phi(x_2)q_2(x_1,x_2) 
E_{x_1,x_2}\left[\prod_{i}\phi(X_i)\right], 
\] 
\[ 
\wt{P}_v(\bullet)
=E_v\left[\prod_{i}\phi(X_i):~ 
\sum_i\delta_{X_i}\in \bullet \right] 
\left(E_v\left[\prod_{i}\phi(X_i)\right]\right)^{-1} 
\] 
and 
\[ 
\wt{P}_{x_1,x_2}(\bullet)
=E_{x_1,x_2}\left[\prod_{i}\phi(X_i):~ 
\sum_i\delta_{X_i}\in \bullet \right] 
\left(E_{x_1,x_2}\left[\prod_{i}\phi(X_i)\right]\right)^{-1},  
\] 
respectively. Notice that the above formula for 
$\wt{P}_v$ (resp. $\wt{P}_{x_1,x_2}$) 
is valid in $\wt{q}_1(v)dv$-a.e. 
(resp. $\wt{q}_1(x_1,x_2)dx_1dx_2$-a.e.) sense, 
so that the denominator can be 
considered to be positive. 
Combining with (\ref{2.14}) and (\ref{2.15}), we obtain  
\[
x_1x_2\phi(x_1+x_2)\wt{q}_2(x_1,x_2)
=\theta (x_1+x_2) \phi(x_1)\phi(x_2)\wt{q}_1(x_1+x_2) 
\] 
and $\wt{P}_{x_1,x_2}=\wt{P}_{x_1+x_2}$. Here are 
two examples of $K_1$ and $F_1$. \\ 
(i) The above two identities show that 
(\ref{2.17}) is satisfied by $\sum \delta_{X_i}$ 
under $\wt{P}$ when we choose 
\[ 
K_1(x_1,x_2)=x_1x_2\phi(x_1+x_2) \quad \mbox{and}  \quad 
F_1(x_1,x_2)=\theta (x_1+x_2) \phi(x_1)\phi(x_2). 
\]  
In this case, (\ref{2.16}) reads 
\begin{eqnarray*} 
L_1^{\sharp}\Phi(\bx) 
& = & 
\frac{1}{2}\sum_{i\ne j}x_ix_j \phi(x_i+x_j) 
\left(\Phi(M_{ij}\bx)-\Phi(\bx)\right) \nonumber \\ 
& & +\frac{\theta}{2}\sum_i x_i^2
\int_0^{1}du \phi(ux_i)\phi((1-u)x_i)
\left(\Phi(S_i^{u}\bx)-\Phi(\bx)\right), 
\end{eqnarray*} 
which defines a bounded operator 
whenever $\phi$ is bounded. For example, 
fixing $s\in (0,1)$ arbitrarily and 
setting $\phi(u)=\one_{(0,s]}(u)$, 
we see easily that $a=P(X_1\le s)$ and 
the associated reversible distribution $\wt{P}$ 
is identified with the conditional law of 
$(X_i)_{i=1}^{\infty}$ given that $X_1\le s$. \\ 
(ii) For another choice 
\[ 
K_1(x_1,x_2)=x_1x_2 \quad \mbox{and}  \quad 
F_1(x_1,x_2)=\theta (x_1+x_2) \phi(x_1)\phi(x_2)/\phi(x_1+x_2) 
\]  
(\ref{2.16}) becomes 
\begin{eqnarray*} 
L_1^{\sharp}\Phi(\bx) 
& = & 
\frac{1}{2}\sum_{i\ne j}x_ix_j 
\left(\Phi(M_{ij}\bx)-\Phi(\bx)\right) \nonumber \\ 
& & +\frac{\theta}{2}\sum_i \frac{x_i^2}{\phi(x_i)}
\int_0^{1}du \phi(ux_i)\phi((1-u)x_i)
\left(\Phi(S_i^{u}\bx)-\Phi(\bx)\right). 
\end{eqnarray*} 
(Compare with the transition kernel 
studied in \cite{MWZZ}.) 
This operator is bounded for any 
uniformly positive, bounded function $\phi$ on $(0,1)$. 
To check, we take $\phi(u)=\exp(bu)$ with 
$b$ being an arbitrary real number  
and then verify that $\wt{P}={\rm PD}(\theta)$ 
and $L_1^{\sharp}=L_1^{(1,\theta)}$. 
The distribution $\wt{P}$ for $\phi(u)=\exp(bu^2)$ 
has been discussed as the equilibrium measure of 
a certain model in population genetics. 
(See \cite{H05} and the references  therein.)

\medskip 

To explore an analogue of Theorem 2.3 
for processes on $\Omega$, we discuss 
the process generated by $L$ in (\ref{2.4}). 
The special case $\theta\Hh\equiv \Hc$ has 
the generator 
\begin{eqnarray} 
L^{(H,\theta)}\Phi(\bz) 
& := & 
\frac{1}{2}\sum_{i\ne j} z_iz_j H(z_i,z_j) 
\left(\Phi(M_{ij}\bz)-\Phi(\bz)\right) \nonumber \\ 
& & +\frac{\theta}{2}\sum_i z_i^{2} 
\int_0^{1}du H(uz_i,(1-u)z_i) 
\left(\Phi(S_i^{u}\bz)-\Phi(\bz)\right),  \label{2.18} 
\end{eqnarray} 
where $H$ is a symmetric, nonnegative 
homogeneous function $H$ 
of degree $\lambda\ge 0$ satisfying (H1) and $\theta>0$. By conditioning (or cut-off) argument, we get 
\begin{th}
(i) Let  $Z=(Z_i)_{i=1}^{\infty}$ be a random element 
of $\Omega$ and suppose that 
the first and second correlation functions 
$r_1$ and $r_2$ of $\xi:=\Xi(Z)$ exist and satisfy 
\be 
K(z_1,z_2)r_2(z_1,z_2)P_{z_1,z_2} 
=F(z_1,z_2)r_1(z_1+z_2) P_{z_1+z_2}, 
                                        \label{2.19} 
\ee 
$\mbox{a.e.-}(z_1,z_2)\in (0,\infty)^2$. 
Here, $P_{z_1,z_2}$ and $P_{z_1+z_2}$ are 
the reduced Palm distributions of $\xi$ and 
the above equality is understood as the one 
between two measures on $\cN$. Then 
the process generated by $L$ is 
reversible with respect to the law of $Z$. 
\\ 
(ii) Let $L^{(H,\theta)}$ be as in (\ref{2.18}). 
Suppose that 
$(X_i)_{i=1}^{\infty}$ is $\rm{PD}(\theta)$-distributed. 
Then, for any $(0,\infty)$-valued random variable $V$ 
independent of $X$, 
the law of an $\Omega$-valued random element 
$(VX_i)_{i=1}^{\infty}$ is a reversible distribution 
of the process generated by $L^{(H,\theta)}$. 
\end{th}
{\it Proof.}~ 
(i) Let $R>0$ be such that $P(|Z|\le R)>0$. 
First, consider the process generated by $L$ 
with initial distribution 
$P^{(R)}(\bullet):=P(Z\in \bullet|~ |Z|\le R)$. 
Then it is clear that such a process lies in 
$\Omega_{\le R}$, and hence 
its generator $L$ is essentially bounded. 
So, just as in the proof of Theorem 2.3 
the proof of the reversibility with respect to 
$P^{(R)}$ can reduce to verifying that 
the equalities corresponding to (\ref{2.19}) hold 
for $\xi=\Xi(Z)$ under the conditional law $P^{(R)}$. 
It is not difficult to show that, 
under $P^{(R)}$, $\xi=\Xi(Z)$ has 
the first and second correlation functions 
\[ 
r_1^{(R)}(z):=r_1(z) \one_{\{z\le R\}}
\frac{P_{z}(|\eta|\le R-z)}{P(|Z|\le R)}, 
\] 
\[ 
r_2^{(R)}(z_1,z_2):=r_2(z_1,z_2) \one_{\{z_1+z_2\le R\}}
\frac{P_{z_1,z_2}(|\eta|\le R-(z_1+z_2))}{P(|Z|\le R)} 
\] 
and the first order and second order 
reduced Palm distributions 
\[
P^{(R)}_z(\bullet)=P_z(\bullet~|~|\eta|\le R-z), 
\quad P^{(R)}_{z_1,z_2}(\bullet)
=P_{z_1,z_2}(\bullet~|~|\eta|\le R-(z_1+z_2)). 
\] 
These formulas combined with (\ref{2.19}) yield 
\[  
K(z_1,z_2)r^{(R)}_2(z_1,z_2)P^{(R)}_{z_1,z_2} 
=F(z_1,z_2)r^{(R)}_1(z_1+z_2) P^{(R)}_{z_1+z_2}, 
\] 
which is sufficient to imply the reversibility of 
the process $\{Z(t):~t\ge 0\}$ 
generated by $L$ with initial distribution $P^{(R)}$ 
for the aforementioned reason. 
Before taking the limit as $R\to\infty$, 
we interpret the reversibility obtained so far 
in terms of conditional expectations 
as follows: for any $n\in\N$, 
$0<t_1<\cdots<t_n<T$ and 
$\Phi_1,\ldots, \Phi_n\in B(\Omega)$ 
\begin{eqnarray*} 
\lefteqn{
E\left[\Phi_1(Z(t_1))\cdots \Phi_n(Z(t_n))
|~|Z(0)|\le R\right]}                                    \\  
& = & 
E\left[\Phi_1(Z(T-t_1))\cdots \Phi_n(Z(T-t_n))
|~|Z(0)|\le R\right]. 
\end{eqnarray*} 
By letting $R\to \infty$ the required reversibility 
has been proved. \\
(ii) 
Consider the lifted process 
$\{\wt{Z}(t)=(VX_i(t))_{i=1}^{\infty}:~t\ge 0\}$, 
where $\{X(t)=(X_i(t))_{i=1}^{\infty}:~t\ge 0\}$ 
is generated by $L_1^{(H,\theta)}$,  
independent of $V$ and 
such that $X(0)=(X_i)_{i=1}^{\infty}$. 
By Lemma 2.1 (ii) $\{\wt{Z}(t)\}$ is generated by 
$\wt{L}$ with $\Hh=H$ and $\Hc=\theta H$. Accordingly 
\[ 
Z(t):=\wt{Z}(|\wt{Z}(0)|^{2+\lambda}t)
=VX(V^{2+\lambda}t)
\] 
is a process generated by $L^{(H,\theta)}$ 
and clearly $Z(0)=(VX_i)_{i=1}^{\infty}$. 
Letting $\rho$ denote the law of $V$, 
we see from the reversibility of $\{X(t)\}$ 
proved in Theorem 2.3 (ii) that for any $n\in\N$, 
$0<t_1<\cdots<t_n<T$ and 
$\Phi_1,\ldots, \Phi_n\in B(\Omega)$
\begin{eqnarray*} 
\lefteqn{
E\left[\Phi_1(Z(t_1))\cdots \Phi_n(Z(t_n))\right]}  \\ 
& = & 
\int\rho(dv)E\left[\Phi_1(vX(v^{2+\lambda}t_1))
\cdots \Phi_n(vX(v^{2+\lambda}t_n))\right]      \\ 
& =&  
\int\rho(dv)E\left[\Phi_1(vX(v^{2+\lambda}T-v^{2+\lambda}t_1))
\cdots\Phi_n(vX(v^{2+\lambda}T-v^{2+\lambda}t_n))\right]     \\
& =&  
E\left[\Phi_1(Z(T-t_1))\cdots \Phi_n(Z(T-t_n))\right]. 
\end{eqnarray*} 
This proves that the law of $Z(0)=(VX_i)_{i=1}^{\infty}$ 
is a reversible distribution of $\{Z(t):~t\ge 0\}$, 
a process generated by $L^{(H,\theta)}$. 
\qed 

\medskip 

In fact, alternative proof of Theorem 2.4 (ii) 
exists and is based on (\ref{2.19}) 
together with the following static result 
on the correlation measures 
and the reduced Palm distributions of 
`the lifted $\rm{PD}(\theta)$ process' 
$\sum\delta_{VX_i}$. 
\begin{lm}
Let $(X_i)_{i=1}^{\infty}$ be 
$\rm{PD}(\theta)$-distributed 
and suppose that a $(0,\infty)$-valued 
random variable $V$ independent of $(X_i)_{i=1}^{\infty}$ is  given. 
Then, for each $k\in\N$, 
the $k$th correlation function $r_{k}$ on $(0,\infty)^k$ of $\sum\delta_{VX_i}$is given by  
\be 
r_{k}(z_1,\ldots,z_k)
=\frac{\theta^k}{z_1\cdots z_k}\int_{(|z|,\infty)}
\rho(dv)\left(1-\frac{|z|}{v}\right)^{\theta-1},  \label{2.20}
\ee 
where $|z|=z_1+\cdots+z_k$ and 
$\rho$ is the law of $V$. 
Moreover, for any $z=(z_1,\ldots,z_k)\in(0,\infty)^k$ 
such that $P(V>|z|)>0$, the expectation $E_{z_1,\ldots,z_k}$ 
with resect to the $k$th-order reduced Palm distribution 
of $\sum\delta_{VX_i}$ at $z=(z_1,\ldots,z_k)$ 
is characterized by the formula 
\be 
E_{z_1,\ldots,z_k}\left[\prod_i\phi(VX_i)\right] 
=
\frac{\ds{\int_{(|z|,\infty)}\rho(dv)
\left(1-\frac{|z|}{v}\right)^{\theta-1}
E\left[\prod_i\phi((v-|z|)X_i)\right]}}
{\ds{\int_{(|z|,\infty)}\rho(dv)\left(1-\frac{|z|}{v}\right)^{\theta-1}}}
,                                                           \label{2.21}
\ee 
in which $\phi\in\wt{B}_+$ is arbitrary.   \par 
In the special case where $V$ has the gamma density 
\be 
\rho_{\theta,b}(v):=\Gamma(\theta)^{-1}
b^{\theta}v^{\theta-1}e^{-bv}\one_{(0,\infty)}(v) \label{2.22}
\ee  
with $b>0$, $\sum \delta_{VX_i}$ 
is a Poisson point process on $(0,\infty)$ 
with mean measure $\theta y^{-1}e^{-by}dy$. 
\end{lm}
{\it Proof.}~Let $f\in B_+^k$ be arbitrary. 
By the assumed independence and 
the Palm formula for the $\rm{PD}(\theta)$ process 
$\sum \delta_{X_i}$ 
\begin{eqnarray*} 
\lefteqn{
E\left[\sum_{i_1,\ldots,i_k(\ne)}f(VX_{i_1},\ldots,VX_{i_k})
\prod_{j\neq i_1,\ldots,i_k} \phi(VX_j) \right]}  \\ 
& = & 
\int\rho(dv)E\left[\sum_{i_1,\ldots,i_k(\ne)}
f(vX_{i_1},\ldots,vX_{i_k})
\prod_{j\neq i_1,\ldots,i_k} \phi(vX_j) \right]     \\ 
& =&  
\int\rho(dv)\int_{\Delta_k^{\circ}} f(vx_1,\ldots,vx_k)
\frac{\theta^k(1-|x|)^{\theta-1}}{x_1\cdots x_k}dx_1\cdots dx_k 
E\left[\prod_{j} \phi(v(1-|x|)X_j) \right]    \\
& =&  
\int_{(0,\infty)^k}dz_1\cdots dz_k f(z_1,\ldots,z_k)  \\ 
& &  
\times \frac{\theta^k}{z_1\cdots z_k} 
\int_{(|z|,\infty)}\rho(dv) 
\left(1-\frac{|z|}{v}\right)^{\theta-1}
E\left[\prod_{j} \phi((v-|z|)X_j) \right]. 
\end{eqnarray*} 
Taking $\phi\equiv 1$ yields (\ref{2.20}). 
Therefore,  
the above equalities suffice to imply (\ref{2.21}). 
 The last assertion follows from 
\[ 
\frac{\theta^k}{z_1\cdots z_k}\int_{(|z|,\infty)} 
\rho_{\theta,b}(v)
\left(1-\frac{|z|}{v}\right)^{\theta-1}dv
=\frac{\theta^k}{z_1\cdots z_k}e^{-b|z|}
=\prod_{i=1}^k\left(\frac{\theta}{z_i}e^{-bz_i}\right), 
\]  
which is nothing but the $k$th correlation function of 
a Poisson point process on $(0,\infty)$ with 
mean density $\theta y^{-1}e^{-by}$.  
\qed 

\medskip 

\noindent 
We call the above-mentioned Poisson process 
the gamma point process with parameter $(\theta,b)$. 
It is worth noting that in view of (\ref{2.20}) and 
(\ref{2.21}) the equalities  
\[ 
z_1z_2r_2(v_1,v_2)=\theta (z_1+z_2)r_1(z_1+z_2) 
\quad \mbox{and} \quad  
E_{z_1,z_2}=E_{z_1+z_2}  
\] 
which correspond to (\ref{2.14}) and (\ref{2.15}), 
respectively, hold true for 
any lifted $\rm{PD}(\theta)$ process, 
and as a result it satisfies also (\ref{2.19}) with 
\be 
K(x,y)=xy H(x,y)  \quad 
\mbox{and} \quad 
F(x,y)=\theta(x+y) H(x,y).                   \label{2.23} 
\ee 
In the forthcoming section it will be shown that 
the time-dependent system of correlation measures 
of the process generated by $L$ solves the 
hierarchical equation (\ref{1.2}) weakly. 
In connection with Theorem 2.4 (ii), we remark that 
the system of the correlation functions $\{r_k\}$ 
given by (\ref{2.20}) is verified directly 
to be a stationary solution to (\ref{1.2}) with (\ref{2.23}). 
In these calculations, merely the following structure 
is relevant: 
\[ 
z_1\cdots z_k r_k(z_1,\cdots, z_k)
=\theta^k g(z_1+\cdots+z_k)
\] 
for some function $g$ independent of $k$ 
(although (\ref{2.20}) shows us the exact form of $g$).

\section{Hierarchical equations for correlation measures}
\setcounter{equation}{0} 
The purpose of this section is to derive equations 
(\ref{1.2}) as those describing the time evolution of 
the correlation measures associated with 
coagulation-fragmentation processes 
introduced in the previous section. 
As far as the reversible process generated by 
$L_1^{(Q,\theta)}$ with $Q\equiv\mbox{const.}$ is concerned, 
such an attempt is found essentially in \cite{MWZZ} 
for the purpose of showing the uniqueness of stationary distributions.  
In that paper, however, the stationary hierarchical equations 
(called `the basic relations' on p.19) for the 
correlation functions are incorrect, 
overlooking a term coming from coagulation between 
clusters with specific sizes given 
(i.e., a term involving $p_{k-1}$). 
It is not clear that developing such an approach could 
make one possible to settle the uniqueness issue 
in much more general setting. 
More specifically, in the reversible case described 
in Theorem 2.3 (resp. Theorem 2.4), 
it seems reasonable to expect the existence of 
a functional of distributions on 
$\Omega_1$ (resp. $\Omega$) which decays 
under the time evolution governed by 
$L_1^{\sharp}$ (resp. $L$). (Notice that 
the uniqueness of stationary distributions 
cannot hold for the processes on $\Omega$.)   

Intending only to derive an infinite system of equations 
describing fully the time evolution of our model, 
we begin by introducing a weak version of (\ref{1.2}), 
namely, the equation obtained by operating on  
a test function by its formal adjoint. 
Since $K$ and $F$ are unbounded, 
suitable integrability conditions must be required for solutions. 
Let $B_{+,c}^k$ denote 
the totality of functions in $B_{+}^k$ with compact support. 
Also, the abbreviated notation $\bz_k=(z_1,\ldots,z_k)$ and 
$d\bz_k=dz_1\cdots dz_k$ are used in the integral expressions.  A family of measures $\{c_k(t,d\bz_k):~t\ge 0, k\in\N\}$ 
is said to be admissible if the following two conditions are fulfilled: \\ 
(A1) Each $c_k(t,\cdot)$ is a locally finite measure on $(0,\infty)^k$.  \\ 
(A2) For any $f \in B_{+,c}^{k}$, $g \in B_{+,c}^{k-1}$
and $l\in\{1,\ldots,k\}$, the three functions below are 
locally integrable on $[0,\infty)$: 
\begin{eqnarray*} 
{\rm (i)}~ t & \mapsto &  \int f(\bz_k)c_k(t,d\bz_k),              \\ 
{\rm (ii)}~ t & \mapsto & 
\int z_l z_{l+1} f(z_1,\ldots,z_{l-1},z_{l}+z_{l+1},z_{l+2}\ldots,z_{k+1}) 
c_{k+1}(t,d\bz_{k+1}),    
                                                      \\ 
{\rm (iii)}~ t & \mapsto & 
\left\{ 
\begin{array}{ll} 
\ds{\int z_{1}(z_{1}^{1+\lambda}\vee 1) c_{1}(t,dz_{1})} & (k=1),       \\ 
\ds{\int z_{l}(z_{l}^{1+\lambda}\vee 1)  
g(z_1,\ldots,z_{l-1},z_{l+1},\ldots,z_{k}) c_{k}(t,d\bz_{k})} & (k\ge 2).  
\end{array} 
\right. 
\end{eqnarray*} 
Given an admissible $\{c_k(t,d\bz_k):~t\ge 0, k\in\N\}$, 
we call it a weak solution of the hierarchical 
coagulation-fragmentation equation (\ref{1.2}) 
with kernels $K$ and $F$ given by (\ref{1.4}) 
if for any $t>0$, $k\in\N$ and $f\in B_{+,c}^k$ 
\begin{eqnarray} 
\lefteqn{
\int f(\bz_k)c_k(t,d\bz_k)
-\int f(\bz_k)c_k(0,d\bz_k)} 
                                                                     \label{3.1}  \\ 
&=& 
\frac{1}{2}\sum_{l=1}^k\int_0^tds \int K(z_l,z_{l+1})
f({\rm Coag}_{l,l+1}\bz_{k+1})c_{k+1}(s,d\bz_{k+1}) 
                                       \nonumber  \\ 
& & 
-\frac{1}{2}\sum_{l=1}^k\int_0^tds 
\int\int_0^{z_l}dy F(y,z_l-y) f(\bz_k)c_{k}(s,d\bz_k) 
                                       \nonumber  \\ 
& & 
-\sum_{l=1}^k\int_0^tds \int K(z_l,z_{l+1})
f(z_1,\ldots,z_l,z_{l+2},\ldots,z_{k+1})
c_{k+1}(s,d\bz_{k+1})
                                       \nonumber  \\ 
& & 
+\sum_{l=1}^k\int_0^tds \int \int_0^{z_l}dyF(y,z_l-y) 
f(z_1,\ldots,z_{l-1},y,z_{l+1},\ldots,z_{k})c_{k}(s,d\bz_{k})  
                                       \nonumber  \\ 
& &  
-\one_{\{k\ge 2\}}\sum_{l<m}^{k}\int_0^tds \int K(z_l,z_m)
f(\bz_k)c_{k}(s,d\bz_{k}) 
                                       \nonumber  \\ 
& &  
+\one_{\{k\ge 2\}}\sum_{l<m}^{k}\int_0^tds 
\int\int_0^{z_l}dyF(y,z_l-y)
f({\rm Frag}_{l,m}^{(y)}\bz_{k-1}) 
c_{k-1}(s,d\bz_{k-1}),            \nonumber 
\end{eqnarray} 
where 
\[ 
{\rm Coag}_{l,l+1}\bz_{k+1}
=(z_1,\ldots,z_{l-1},z_l+z_{l+1},z_{l+2},\ldots,z_{k+1})  
\] 
and 
\[ 
{\rm Frag}_{l,m}^{(y)}\bz_{k-1}
=(z_1,\ldots,z_{l-1},y,z_{l+1},\ldots, 
z_{m-1},z_l-y,z_{m},\ldots,z_{k-1}). 
\] 
By (\ref{1.4}), the equation 
(\ref{3.1}) actually takes a more specific form 
\begin{eqnarray} 
\lefteqn{
\int f(\bz_k)c_k(t,d\bz_k)
-\int f(\bz_k)c_k(0,d\bz_k)} 
                                                                    \label{3.2}  \\ 
&=& 
\frac{1}{2}\sum_{l=1}^k\int_0^tds \int z_l z_{l+1}
\Hh(z_l,z_{l+1})f({\rm Coag}_{l,l+1}\bz_{k+1})
c_{k+1}(s,d\bz_{k+1}) 
                                       \nonumber  \\ 
& & 
-\frac{1}{2}\sum_{l=1}^k\int_0^tds \int z_l 
\int_0^{z_l}dy\Hc(y,z_l-y)f(\bz_k)c_{k}(s,d\bz_k) 
                                       \nonumber  \\ 
& & 
-\sum_{l=1}^k\int_0^tds \int z_l z_{l+1}
\Hh(z_l,z_{l+1})f(z_1,\ldots,z_l,z_{l+2},\ldots,z_{k+1})
c_{k+1}(s,d\bz_{k+1})
                                       \nonumber  \\ 
& & 
+\sum_{l=1}^k\int_0^tds \int z_l \int_0^{z_l}dy
\Hc(y,z_l-y)f(z_1,\ldots,z_{l-1},y,z_{l+1},\ldots,z_k)
c_{k}(s,d\bz_{k})  
                                       \nonumber  \\ 
& &  
-\one_{\{k\ge 2\}}\sum_{l<m}^{k}\int_0^tds \int z_l z_m
\Hh(z_l,z_m)f(\bz_k)c_{k}(s,d\bz_{k}) 
                                       \nonumber  \\ 
& &  
+\one_{\{k\ge 2\}}\sum_{l<m}^{k}\int_0^t ds  
\int z_l  \int_0^{z_l} dy \Hc(y,z_l-y)
f({\rm Frag}_{l,m}^{(y)}\bz_{k-1}) c_{k-1}(s,d\bz_{k-1})    \nonumber  \\ 
& =: & 
I_1-I_2-I_3+I_4-I_5+I_6.                                              \label{3.3} 
\end{eqnarray} 
Obviously, (A1) ensures that two integrals on 
the left side of (\ref{3.1}) is finite. 
Considering the terms on the right side, 
we prepare the following bounds:  
by homogeneity (\ref{1.5}), (H1) and (H2) together 
\begin{eqnarray}
\Hh(x,y) 
& = & 
\Hh\left((x+y)\frac{x}{x+y},(x+y)\frac{x}{x+y}\right)   
 \le \Ch(x+y)^{\lambda}                                         \label{3.4}    \\ 
& \le & \Ch(1+x)^{\lambda} (y^{\lambda}\vee 1)  
\ \le \  \Ch(1+x)^{\lambda} (y^{1+\lambda}\vee 1)      \label{3.5}
\end{eqnarray} 
and 
\be 
\int_0^{x}dy \Hc(y,x-y) 
=x \int_0^1du \Hc(ux,(1-u)x) = \Cc x^{1+\lambda}.              \label{3.6}
\ee 
\begin{lm}
Assume that $\{c_k(t,d\bz_k):~t\ge 0, k\in\N\}$ is admissible. 
Then every term on the right side of (\ref{3.1}) 
(or equivalently of (\ref{3.2})) is finite. 
\end{lm}
{\it Proof.}~ 
We discuss $I_i$'s on (\ref{3.3}) 
instead of the terms on the right side of (\ref{3.1}). 
It follows from (A2) (i) 
that $I_2$ and $I_5$ are finite. 
Also, $I_1$ is finite because of (\ref{3.4}) and (A2) (ii). 
$I_3$ converges by (\ref{3.5}) together with (A2) (iii), 
and similarly the finiteness of $I_4$ is 
due to (\ref{3.6}) and (A2) (iii) 
with 
\[ 
g(z_{1},\ldots,z_{k-1}) =
\sup_{y>0} f(z_{1},\ldots,z_{{l-1}},y,z_{{l}},\ldots,z_{k-1}).  
\] 
Lastly, again by (A2) (i), $I_6$ is finite 
since each function 
\[ 
h_{l,m}(z_{1},\ldots,z_{k-1}):=
z_l \! \int_0^{z_l} \! dy \Hc(y,z_l-y)
f({\rm Frag}_{l,m}^{(y)}\bz_{k-1}) 
\] 
is an element of $B_{+,c}^{k-1}$ for $k\ge 2$.  
Indeed, taking $\epsilon>0$ so that 
$f(z_{1},\ldots,z_{k})=0$ whenever $0<z_l<\epsilon$,  
we see that $h_{l,m}(z_1,\ldots,z_{k-1})=0$ 
for any $z_l\in(0,\epsilon)$, and 
analogously, taking $R>0$ so that 
$f(z_{1},\ldots,z_{k})=0$ whenever 
$z_l>R$ or $z_m>R$, we see that  
$h_{l,m}(z_1,\ldots,z_{k-1})=0$ for any $z_l>2R$. 
The proof of Lemma 3.1 is complete. \qed 

\medskip 

The main result of this section yields stochastic construction 
of a solution to (\ref{1.2}) with kernels we are concerned with. 
\begin{th}
(i) Let $\{Z(t):~t\ge 0\}$ be the process 
generated by $L$ and 
suppose that $E[|Z(0)|^k]<\infty$ for all $k\in\N$.   
For each  $t\ge 0$ and $k\in \N$ 
denote by $r_k(t,d\bz_k)$ 
the $k$th correlation measure of $\Xi(Z(t))$.  
Then 
$\{r_k(t,d\bz_k):~t\ge 0, k\in\N\}$ is admissible and 
solves weakly the hierarchical 
coagulation-fragmentation equation (\ref{1.2}) with 
kernels $K$ and $F$ given by (\ref{1.4}). \\ 
(ii) Let $\wh{Q}$ and $\check{Q}$ be 
symmetric, nonnegative bounded functions on 
$\{(x,y) |~x,y> 0, x+y \le 1\}$ and set 
\[ 
K_1(x,y)=xy\wh{Q}(x,y) \quad  \mbox{and} \quad  
F_1(x,y)=(x+y)\check Q(x,y). 
\]  
Let $\{X(t):~t\ge 0\}$ be the process 
generated by $L_1^{\sharp}$ in (\ref{2.16}). 
For each  $t\ge 0$ and $k\in \N$ 
denote by $q_k(t,d\bx_k)$ 
the $k$th correlation measure of $\Xi(X(t))$.  
Then 
$\{q_k(t,d\bz_k):~t\ge 0, k\in\N\}$ is admissible 
and solves weakly the hierarchical 
coagulation-fragmentation equation (\ref{1.2}) with 
kernels $K_1$ and $F_1$.   
\end{th}
Recalling the definition of correlation measures 
(cf. (\ref{2.7})), 
the proof of this theorem is basically done by 
calculating carefully $L\Phi(\bz)$ or  
$L_1^{\sharp}\Phi(\bz)$ for 
\be  
\Phi(\bz)=\sum_{i_1,\ldots,i_k(\ne)}f(z_{i_1},\ldots,z_{i_k}),  
                                                 \label{3.7}
\ee 
where $f\in B_{+,c}^k$ is arbitrary. 
Here, we understand that 
$f(z_{i_1},\ldots,z_{i_k})=0$ when 
$z_{i_1} \cdots z_{i_k}=0$, so that 
\[ 
\sum_{i_1,\ldots,i_k(\ne)}f(z_{i_1},\ldots,z_{i_k})
= 
\int f(y_1,\ldots,y_k)\xi^{[k]}(dy_1\cdots dy_k), 
\] 
where $\xi=\Xi(\bz)$. 
Although such function's $\Phi$ on $\Omega$ 
may be unbounded,  we can control its growth order 
as will be seen in the next lemma. For each $a>0$ 
denote by $\cF_a$ the class of measurable functions $\Psi$ 
on $\Omega$ such that, for some constant $C<\infty$,   
$|\Psi(\bz)|\le C|\bz|^a$ for all $\bz\in\Omega$. 
\begin{lm}
Let $f\in B_{+,c}^k$ and $\Phi$ be as in (\ref{3.7}). 
Then $\Phi\in\cF_k$. 
\end{lm}
{\it Proof.}~ 
Define 
$\epsilon =\inf\{\min\{z_1,\ldots,z_k\}:(z_1,\ldots,z_k)
\in {\rm supp}(f)\}$, which is strictly positive 
because ${\rm supp}(f)$ is assumed to be 
a compact subset of  $(0,\infty)^k$. Therefore 
\[ 
0 \le \Phi(z) \le \Vert f\Vert_{\infty}
\sum_{i_1,\ldots,i_k} 
\frac{z_{i_1}}{\epsilon}\cdots \frac{z_{i_k}}{\epsilon}
\one_{\{(z_{i_1},\ldots,z_{i_k})\in {\rm supp}(f)\}} 
\le \Vert f\Vert_{\infty}\frac{|\bz|^{k}}{\epsilon^k}. 
\]  
This proves Lemma 3.3. \qed

\medskip 

\noindent 
At the core of our proof of Theorem 3.2 is 
\begin{lm}
Let $f\in B_{+,c}^k$ and $\Phi$ be as in (\ref{3.7}). 
Then for each $\bz=(z_i)_{i=1}^{\infty}\in \Omega$  
\begin{eqnarray} 
\lefteqn{L\Phi(\bz)} 
                                                       \label{3.8} \\ 
& = & 
\frac{1}{2}
\sum_{l=1}^k \sum_{i_1,\ldots,i_{k+1}(\ne)}K(z_{i_l},z_{i_{l+1}})
f({\rm Coag}_{l,l+1}(z_{i_1},\ldots,z_{i_{k+1}})) 
                                                    \nonumber \\ 
&  & 
-\frac{1}{2}\sum_{l=1}^k \sum_{i_1,\ldots,i_{k}(\ne)}
\int_0^{z_{i_l}}dyF(y,z_{i_l}-y)f(z_{i_1},\ldots,z_{i_k}).   \nonumber  \\ 
&  & 
-\sum_{l=1}^k \sum_{i_1,\ldots,i_{k+1}(\ne)}K(z_{i_l},z_{i_{l+1}})
f(z_{i_1},\ldots,z_{i_l},z_{i_{l+2}},\ldots,z_{i_{k+1}})    \nonumber \\
&  & 
+\sum_{l=1}^k \sum_{i_1,\ldots,i_{k}(\ne)} 
\int_0^{z_{i_l}}dyF(y,z_{i_l}-y) f(z_{i_1},\ldots,z_{i_{l-1}},y,z_{i_{l+1}},\ldots,z_{i_{k}})
                                       \nonumber \\ 
& &   
-\one_{\{k\ge 2\}}\sum_{l<m}^{k} 
\sum_{i_1,\ldots,i_{k}(\ne)}K(z_{i_l},z_{i_m})
f(z_{i_1},\ldots,z_{i_{k}})                         \nonumber  \\ 
&   & 
+\one_{\{k\ge 2\}}\sum_{l<m}^{k} 
\sum_{i_1,\ldots,i_{k-1}(\ne)} \int_0^{z_{i_l}}dyF(y,z_{i_l}-y) 
f({\rm Frag}_{l,m}^{(y)}(z_{i_1},\ldots,z_{i_{k-1}}))
                                       \nonumber \\ 
& =: & 
\frac{1}{2}\sum_{l=1}^k\Psi_l^{(1)}(\bz) 
-\frac{1}{2}\sum_{l=1}^k\Psi_l^{(2)}(\bz)
-\sum_{l=1}^k \Psi_l^{(3)}(\bz)            \nonumber  \\
& & 
+\sum_{l=1}^k \Psi_l^{(4)}(\bz) 
-\one_{\{k\ge 2\}}\sum_{l<m}^{k} \Psi_{l,m}^{(5)}(\bz)
+\one_{\{k\ge 2\}}\sum_{l<m}^{k} \Psi_{l,m}^{(6)}(\bz). \nonumber 
\end{eqnarray} 
Moreover, $\Psi_l^{(1)}\in \cF_{k+1}$, 
$\Psi_l^{(2)}\in \cF_{k}$, $\Psi_l^{(3)}\in \cF_{k+1+\lambda}$,  
$\Psi_l^{(4)}\in \cF_{k+1+\lambda}$, 
$\one_{\{k\ge 2\}}\Psi_{l,m}^{(5)}\in \cF_{k}$ 
and $\one_{\{k\ge 2\}}\Psi_{l,m}^{(6)}\in \cF_{k+2}$.  
Thus $L\Phi$ belongs to 
the linear span $\overline{\cF}$ of $\bigcup_{a>0}\cF_a$. 
\end{lm}
{\it Proof.}~
The main proof of (\ref{3.8}) consists of 
almost algebraic calculations (which are completely 
independent of  other arguments)  
and so it is postponed. 
We here prove only the assertions for 
$\Psi_l^{(i)} (i\in\{1,2,3,4\})$ 
and $\Psi_{l,m}^{(j)} (j\in\{5,6\})$.  
The arguments below are based on similar observations to 
those in the proof of Lemma 3.1. 
From Lemma 3.3, it is evident that 
$\Psi_l^{(2)}$ and $\Psi_{l,m}^{(5)}$ are elements of $\cF_k$. 
Relying on the fact that the function $h_{l,m}$ 
in the proof of Lemma 3.1 belongs to $B_{+,c}^{k-1}$, 
one can verify similarly that 
$\Psi_{l,m}^{(6)}\in\cF_{k-1}$ for $k\ge 2$. 
As for $\Psi_l^{(1)}$, taking $\epsilon>0$ and $R>\epsilon$ 
such that $[\epsilon, R]^k\supset {\rm supp}(f)$ 
and using (\ref{3.4}), 
we observe as in the proof of Lemma 3.3 that 
for $l\in\{1,\ldots,k\}$ 
\begin{eqnarray*} 
|\Psi_{l}^{(1)}(\bz)|
&\le &
\Ch \sum_{i_1,\ldots,i_{k+1}(\ne)}
z_{i_l}z_{i_{l+1}}(z_{i_l}+z_{i_{l+1}})^{\lambda}
f({\rm Coag}_{l,l+1}(z_{i_1},\ldots,z_{i_{k+1}}))                  \\ 
& \le  &
\Ch \Vert f\Vert_{\infty} 
\sum_{i_1,\ldots,i_{k+1}} z_{i_l}z_{i_{l+1}} R^{\lambda}
\frac{z_{i_1}}{\epsilon}\cdots \frac{z_{i_{l-1}}}{\epsilon}   
\frac{z_{i_{l+2}}}{\epsilon}\cdots \frac{z_{i_{k+1}}}{\epsilon}       \\  
& = & 
\Ch R^{\lambda} \Vert f\Vert_{\infty}
\frac{|\bz|^{k+1}}{\epsilon^{k-1}}, 
\end{eqnarray*} 
by which $\Psi_{l}^{(1)}\in\cF_{k+1}$. Analogously 
\begin{eqnarray*} 
|\Psi_{l}^{(3)}(\bz)|
&\le &
\Ch \sum_{i_1,\ldots,i_{k+1}(\ne)}
z_{i_l}z_{i_{l+1}}(z_{i_l}+z_{i_{l+1}})^{\lambda} 
f(z_{i_1},\ldots,z_{i_l},z_{i_{l+2}},\ldots,z_{i_{k+1}})                  \\ 
& \le  &
\Ch \sum_{i_1,\ldots,i_{k+1}(\ne)} R z_{i_{l+1}} |z|^{\lambda}
f(z_{i_1},\ldots,z_{i_l},z_{i_{l+2}},\ldots,z_{i_{k+1}})        \\  
& \le  & 
\Ch R |\bz|^{1+\lambda} 
\sum_{i_1,\ldots,i_{k}(\ne)} f(z_{i_1},\ldots,z_{i_{k}}), 
\end{eqnarray*} 
which combined with Lemma 3.3 implies 
that $\Psi_l^{(3)}\in\cF_{k+1+\lambda}$. 
Lastly, we shall show that $\Psi_l^{(4)}\in\cF_{k+1+\lambda}$. 
For $k=1$ we have by (\ref{3.6}) 
\[ 
|\Psi_1^{(4)}(\bz)| 
\le \Cc \sum_i z_i^{2+\lambda}  \Vert f\Vert_{\infty} 
\le \Cc |\bz|^{2+\lambda}  \Vert f\Vert_{\infty} 
\] 
and hence $\Psi_1^{(4)}\in\cF_{2+\lambda}$. 
For $k\ge 2$, by noting that 
\[ 
{\overline f}_l(z_1,\ldots,z_{k-1})
:= \sup_{y>0} f(z_{1},\ldots,z_{{l-1}},y,z_{{l}},\ldots,z_{{k-1}})
\] 
belongs to $B_{+,c}^{k-1}$ and observing that 
\begin{eqnarray*} 
|\Psi_l^{(4)}(\bz)| 
& \le &  
\Cc\sum_{i_1,\ldots,i_{k}(\ne)} z_{i_l}^{2+\lambda} 
{\overline f}_l(z_{i_1},\ldots,z_{i_{l-1}},z_{i_{l+1}},\ldots,z_{i_{k}})   \\ 
& \le & 
\Cc  |\bz|^{2+\lambda}\sum_{i_1,\ldots,i_{k-1}(\ne)}
{\overline f}_l(z_{i_1},\ldots,z_{i_{k-1}}), 
\end{eqnarray*} 
we deduce from Lemma 3.3 
that $\Psi_l^{(4)}\in\cF_{k+1+\lambda}$. 
Consequently, we have shown 
that $L\Phi\in\overline{\cF}$.  
\qed 

\medskip 

\noindent  
{\it Proof of Theorem 3.2 (i).}~
Let $f\in B_{+,c}^k$ be arbitrary. 
First, we must show the admissibility of $r_k(t,d\bz_k)$. 
The conditions (A1) and (A2) (i) for $r_k(t,d\bz_k)$ 
are easily seen to hold by combining Lemma 3.3 with 
\[ 
E[|Z(t)|^a]=E[|Z(0)|^a]=:m_a<\infty 
\] 
for all $a\ge 1$. 
(A2) (ii) can be verified by a similar bound 
to that for $|\Psi_l^{(1)}(\bz)|$ in the proof of Lemma 3.4. 
Moreover, (A2) (iii) follows from observations that 
\begin{eqnarray*} 
\lefteqn{
\int z_{1} (z_{1}^{1+\lambda}\vee 1) r_{1}(t,dz_{1})  
\le E\left[\sum_i(Z_{i}(t)+Z_{i}(t)^{2+\lambda})\right]}        \\ 
& \le & 
E\left[|Z(t)|+|Z(t)|^{2+\lambda}\right]  
\ = \ m_{1}+m_{2+\lambda} 
\end{eqnarray*} 
and that for $k\ge 2$ and $l\in\{1,\ldots,k\}$ 
\begin{eqnarray*} 
\lefteqn{\int z_{l} (z_{l}^{1+\lambda}\vee 1)   
g(z_1,\ldots,z_{l-1},z_{l+1},\ldots,z_{k}) r_{k}(t,d\bz_{k})}  \\ 
& \le & 
C\int (z_{l}+z_{l}^{2+\lambda})  
z_1\cdots z_{l-1} z_{l+1} \cdots z_{k} r_{k}(t,d\bz_{k})    \\ 
& \le & 
C E\left[(|Z(t)|+|Z(t)|^{2+\lambda})|Z(t)|^{k-1}\right] 
\ = \ C \left(m_k+m_{1+\lambda+k}\right),  
\end{eqnarray*} 
where $C$ is a finite constant. 
Second, we claim that, for $\Phi$ given by (\ref{3.7}) and $t>0$ 
\be 
\int f(\bz_k)r_k(t,d\bz_k)
-\int f(\bz_k)r_k(0,d\bz_k) 
=\int_0^tds E[L\Phi(Z(s))].  
                                                 \label{3.9}
\ee 
Define, for each $R>0$, 
$\Phi^{(R)}(\bz)=\one_{\{|\bz|\le R\}}\Phi(\bz)$. 
Then $\Phi^{(R)}\in B(\Omega)$ by Lemma 3.3 and 
clearly 
$L\Phi^{(R)}(\bz)=\one_{\{|\bz|\le R\}}L\Phi(\bz)$. 
Furthermore, by virtue of Lemma 3.4, 
$\one_{\{|\bz|\le R\}}L\Phi(\bz)$ is bounded 
and so is $L\Phi^{(R)}$.  
These observations together imply that 
\begin{eqnarray*} 
E\left[\Phi^{(R)}(Z(t))\right]
-E\left[\Phi^{(R)}(Z(0))\right] 
&=&
\int_0^tds E\left[L\Phi^{(R)}(Z(s))\right]                         \\ 
&=&
\int_0^tds E\left[\one_{\{|Z(s)|\le R\}}L\Phi(Z(s))\right]. 
\end{eqnarray*} 
Noting that every moment of 
$|Z(t)|=|Z(0)|$ is finite 
by the assumption and that $L\Phi\in\overline{\cF}$, 
we get (\ref{3.9}) by taking the limit as $R\to\infty$ 
with the help of Lebesgue's convergence theorem. 

Integrating the right side of (\ref{3.8}) 
with respect to the law of $Z(s)$ 
and then plugging the resulting expression for the expectation 
$E[L\Phi(Z(s))]$ into (\ref{3.9}) yield 
(\ref{3.1}) with $r_k(s,d\bz_k)$ in place of $c_k(s,d\bz_k)$. 
We thus obtained the required equations 
for \{$r_k(t,d\bz_k):~t\ge 0, k\in \N\}$ 
and the proof of Theorem 3.2 (i) is complete,  
provided that the identity (\ref{3.8}) is shown by 
calculations which are self-contained.    

\medskip 

\noindent 
{\it Proof of (\ref{3.8}).}~ 
Recalling the definition (\ref{2.4}) of $L$, 
we now calculate $L\Phi(\bz)$ for $\Phi$ 
of the form (\ref{3.7}). 
Observe that, for any $i\ne j$ such that $z_iz_j>0$,  
the `coagulation difference' $\Phi(M_{ij}\bz)-\Phi(\bz)$ 
equals 
\begin{eqnarray*} 
\lefteqn{\sum_{l=1}^k \sum_{i_1,\ldots,i_{k-1}(\ne)}^{\{i,j\}^c}
f(z_{i_1},\ldots,z_{i_{l-1}},z_i+z_j,z_{i_l},\ldots,z_{i_{k-1}})}    \\ 
& - & 
\sum_{l=1}^k \sum_{i_1,\ldots,i_{k-1}(\ne)}^{\{i,j\}^c}
f(z_{i_1},\ldots,z_{i_{l-1}},z_i,z_{i_l},\ldots,z_{i_{k-1}})    \nonumber \\
& -  &
\sum_{l=1}^k \sum_{i_1,\ldots,i_{k-1}(\ne)}^{\{i,j\}^c}
f(z_{i_1},\ldots,z_{i_{l-1}},z_j,z_{i_l},\ldots,z_{i_{k-1}})
                                       \nonumber \\ 
& - & 
\one_{\{k\ge 2\}}\sum_{l<m}^{k} 
\sum_{i_1,\ldots,i_{k-2}(\ne)}^{\{i,j\}^c}
f(z_{i_1},\ldots,z_{i_{l-1}},z_i,z_{i_l},\ldots,z_{i_{m-2}},z_j,
z_{i_{m-1}},\ldots,z_{i_{k-2}})
                                       \nonumber \\ 
& -&   
\one_{\{k\ge 2\}}\sum_{l<m}^{k} 
\sum_{i_1,\ldots,i_{k-2}(\ne)}^{\{i,j\}^c}
f(z_{i_1},\ldots,z_{i_{l-1}},z_j,z_{i_l},\ldots,z_{i_{m-2}},z_i,
z_{i_{m-1}},\ldots,z_{i_{k-2}})                       \nonumber  \\ 
=: \ \ 
\lefteqn{\Sigma_{ij}^{(1)}(\bz)
-\Sigma_{ij}^{(2)}(\bz)   
-\Sigma_{ij}^{(3)}(\bz) 
-\Sigma_{ij}^{(4)}(\bz)  
-\Sigma_{ij}^{(5)}(\bz),}                 \nonumber 
\end{eqnarray*} 
where $\sum_{i_1,\ldots,i_{k-1}(\ne)}^{\{i,j\}^c}$ stands for the 
sum taken over $(k-1)$-tuples $(i_1,\ldots,i_{k-1})$ 
of positive integers such that $i_1,\ldots,i_{k-1}\in\{i,j\}^c$ 
are mutually distinct. 
Noting that $\Sigma_{ij}^{(2)}(\bz)=\Sigma_{ji}^{(3)}(\bz)$ and 
$\Sigma_{ij}^{(4)}(\bz)=\Sigma_{ji}^{(5)}(\bz)$, we get 
\begin{eqnarray} 
\lefteqn{\frac{1}{2}
\sum_{i\ne j}K(z_i,z_j)(\Phi(M_{ij}\bz)-\Phi(\bz))} \label{3.10} \\  
& = & 
\frac{1}{2}
\sum_{l=1}^k \sum_{i_1,\ldots,i_{k+1}(\ne)}K(z_{i_l},z_{i_{l+1}}) 
f(z_{i_1},\ldots,z_{i_{l-1}},z_{i_l}+z_{i_{l+1}},z_{i_{l+2}},\ldots,z_{i_{k+1}}) \nonumber \\ 
&  & 
-\sum_{l=1}^k \sum_{i_1,\ldots,i_{k-1}(\ne)}K(z_{i_l},z_{i_{l+1}}) 
f(z_{i_1},\ldots,z_{i_{l-1}},z_{i_l},z_{i_{l+2}},\ldots,z_{i_{k+1}})    \nonumber \\
& &   
-\one_{\{k\ge 2\}}\sum_{l<m}^{k} 
\sum_{i_1,\ldots,i_{k}(\ne)}K(z_{i_l},z_{i_m}) 
f(z_{i_1},\ldots,z_{i_{k}}).                                     \nonumber 
\end{eqnarray} 

Similarly, for each $i\in\N$ with $z_i>0$ 
and any $y\in(0,z_i)$, 
the `fragmentation difference' $\Phi(S_{i}^{(y)}\bz)-\Phi(\bz)$ 
is expressed as 
\begin{eqnarray*} 
\lefteqn{\sum_{l=1}^k \sum_{i_1,\ldots,i_{k-1}(\ne)}^{\{i\}^c}
f(z_{i_1},\ldots,z_{i_{l-1}},y,z_{i_l},\ldots,z_{i_{k-1}})}   \\ 
& + & 
\sum_{l=1}^k \sum_{i_1,\ldots,i_{k-1}(\ne)}^{\{i\}^c}
f(z_{i_1},\ldots,z_{i_{l-1}},z_i-y,z_{i_l},\ldots,z_{i_{k-1}})
                                       \nonumber \\ 
& +  & 
\one_{\{k\ge 2\}}\sum_{l<m}^{k} 
\sum_{i_1,\ldots,i_{k-2}(\ne)}^{\{i\}^c}
f({\rm Frag}_{l,m}^{(y)}(z_{i_1},\ldots,z_{i_{l-1}},z_i,z_{i_l},\ldots,z_{i_{k-2}}))
                                       \nonumber \\ 
& + & 
\one_{\{k\ge 2\}}\sum_{l<m}^{k} 
\sum_{i_1,\ldots,i_{k-2}(\ne)}^{\{i\}^c}
f({\rm Frag}_{l,m}^{(z_i-y)}(z_{i_1},\ldots,z_{i_{l-1}},z_i,z_{i_l},\ldots,z_{i_{k-2}})) 
                                       \nonumber \\ 
& - & 
\sum_{l=1}^k \sum_{i_1,\ldots,i_{k-1}(\ne)}^{\{i\}^c}
f(z_{i_1},\ldots,z_{i_{l-1}},z_i,z_{i_l},\ldots,z_{i_{k-1}})  \nonumber \\ 
=: \ \ 
\lefteqn{\Sigma_{i}^{(6)}(y,\bz)
+\Sigma_{i}^{(7)}(y,\bz)   
+\Sigma_{i}^{(8)}(y,\bz) 
+\Sigma_{i}^{(9)}(y,\bz)  
-\Sigma_{i}^{(10)}(\bz)}                   \nonumber 
\end{eqnarray*} 
because for $1\le l<m\le k$  
\begin{eqnarray*} 
\lefteqn{{\rm Frag}_{l,m}^{(y)}(z_{i_1},\ldots,z_{i_{l-1}},
z_i,z_{i_l},\ldots,z_{i_{k-2}})}                    \\ 
&= &
(z_{i_1},\ldots,z_{i_{l-1}},y,z_{i_l},\ldots,z_{i_{m-2}},
z_i-y,z_{i_{m-1}},\ldots,z_{i_{k-2}})
\end{eqnarray*} 
and 
\begin{eqnarray*} 
\lefteqn{{\rm Frag}_{l,m}^{(z_i-y)}(z_{i_1},\ldots,z_{i_{l-1}},
z_i,z_{i_l},\ldots,z_{i_{k-2}})}                   \\ 
& = & 
(z_{i_1},\ldots,z_{i_{l-1}},z_i-y,z_{i_l},\ldots,z_{i_{m-2}},
y,z_{i_{m-1}},\ldots,z_{i_{k-2}}). 
\end{eqnarray*} 
Here, $\sum_{i_1,\ldots,i_{k-1}(\ne)}^{\{i\}^c}$ 
indicates the sum taken over $(k-1)$-tuples $(i_1,\ldots,i_{k-1})$ 
of distinct positive integers which are different from $i$. 
By noting two identities 
$\Sigma_{i}^{(6)}(y,\bz)=\Sigma_{i}^{(7)}(z_i-y,\bz)$ 
and $\Sigma_{i}^{(8)}(y,\bz)=\Sigma_{i}^{(9)}(z_i-y,\bz)$ 
\begin{eqnarray} 
\lefteqn{\frac{1}{2}
\sum_iz_i\int_0^{z_i}dy F(y,z_i-y)(\Phi(S_i^{(y)}\bz)-\Phi(\bz))}   \label{3.11} \\ 
& = & 
\sum_{l=1}^k \sum_{i_1,\ldots,i_{k}(\ne)}z_{i_l}
\int_0^{z_{i_l}}dy F(y, z_{i_l}-y) 
f(z_{i_1},\ldots,z_{i_{l-1}},y,z_{i_{l+1}},\ldots,z_{i_{k}})
                                       \nonumber \\ 
&   & 
+\one_{\{k\ge 2\}}\sum_{l<m}^{k} 
\sum_{i_1,\ldots,i_{k-1}(\ne)}z_{i_l}
\int_0^{z_{i_l}}dy F(y, z_{i_l}-y) 
f({\rm Frag}_{l,m}^{(y)}(z_{i_1},\ldots,z_{i_{k-1}}))
                                       \nonumber \\ 
&  & 
-\frac{1}{2}\sum_{l=1}^k \sum_{i_1,\ldots,i_{k}(\ne)}
z_{i_l}\int_0^{z_{i_l}}dy F(y, z_{i_l}-y) f(z_{i_1},\ldots,z_{i_k}).   \nonumber  
\end{eqnarray} 
Consequently, (\ref{3.8}) is deduced from (\ref{3.10}) and (\ref{3.11}). 
\qed  

\medskip 

\noindent 
{\it Proof of Theorem 3.2 (ii).}~
The proof is essentially the same as the proof of (i). 
We just give some comments. 
The admissibility is shown similarly by $|X(0)|=1$. 
Concerning the analogue of Lemma 3.4, 
(\ref{3.8}) with $K_1$ and $F_1$ 
in place of $K$ and $F$, respectively, holds true  
since the proof of (\ref{3.8}) itself is almost algebraic. 
Moreover, the assertions corresponding to 
the second half of Lemma 3.4 (i.e., the assertions for 
$\Psi_l^{(i)}$'s and $\Psi_{l,m}^{(i)}$'s) 
are also valid for $\lambda=0$ 
as is easily seen from the proof. 
We complete the proof of Theorem 3.2       \qed 

\medskip 

So far we have discussed only weak solutions. 
The final result in this section shows, 
under a certain condition on the initial distribution, 
the existence of 
a strong solution to (\ref{1.2}) with 
$\Hc=\theta\Hh$ for some constant $\theta>0$,  
which ensures reversibility of the underlying process 
as was shown in Theorem 2.4 (ii). 
We will require also for the initial distribution 
to be regarded as a `perturbation' 
from some reversible distribution, 
namely the law of a lifted ${\rm PD}(\theta)$ process.   
Recall that a lifted $\rm{PD}(\theta)$ process is 
of the form $\sum\delta_{VX_i}$, where 
a $(0,\infty)$-valued random variable $V$ 
and a $\rm{PD}(\theta)$-distributed random element 
$(X_i)_{i=1}^{\infty}$ of $\Omega_1$ are mutually independent. 
\begin{pr}
Let $L^{(H,\theta)}$ be as in (\ref{2.18}) and 
$\{Z(t)=(Z_i(t))_{i=1}^{\infty}: t\ge 0\}$ 
be a process generated by $L^{(H,\theta)}$. 
Set $\xi(t)=\Xi(Z(t))$ and 
denote by $\{\cT_t\}_{t\ge 0}$ the semigroup 
associated with $\{\xi(t): t\ge 0\}$. 
Suppose that the law of $\xi(0)$ is 
absolutely continuous with respect to 
the law of a lifted $\rm{PD}(\theta)$ process 
$\sum \delta_{VX_i}$. \\ 
(i) For any $k\in\N$ and $t\ge 0$, 
the $k$th correlation function of $\xi(t)$ is given by  
\be 
r_k(t,\bz_k):=
r_{k}(z_1,\ldots,z_k) 
\int_{\cN}P_{z_1,\ldots, z_k}(d\eta)(\cT_t\Psi^*)
(\eta+\delta_{z_1}+\cdots+\delta_{z_k}),        \label{3.12}
\ee 
where $r_{k}(z_1,\ldots,z_k)$ is 
the $k$th correlation function (\ref{2.20}) 
of $\sum \delta_{VX_i}$,  
$P_{z_1,\ldots,z_k}$ is the 
$k$th-order reduced Palm distribution of $\sum \delta_{VX_i}$ 
at $\bz_k=(z_1,\ldots,z_k)$ characterized by (\ref{2.21}) 
and $\Psi^*$ is the density of 
the law of $\xi(0)$ with respect to 
the law of $\sum \delta_{VX_i}$. 
(It is understood that $r_k(t,\bz_k)=0$ 
whenever $P(V>{|\bz_k|})=0$.)  \\ 
(ii) Suppose additionally that 
$E[|Z(0)|^k]<\infty$ for all $k\in\N$.  
Then the family of nonnegative measurable functions 
$\{r_k(t,\bz_k):~t\ge 0, k\in\N\}$ given in (i) solves 
the equation (\ref{1.2}) with 
$K(x,y)=xy H(x,y)$ and $F(x,y)=\theta(x+y) H(x,y)$ 
in the following sense: for any $k\in N$ and $t\ge 0$ 
\be 
r_k(t,\bz_k)-r_k(0,\bz_k)-
\int_0^t \cL_k(s,\bz_k)ds=0,  \ \mbox{a.e.} \ 
\bz_k\in(0,\infty)^k,                            \label{3.13} 
\ee 
where $\cL_k(s,\bz_k)$ is the right side of (\ref{1.2}) 
with $c_{k+1}$,  $c_{k}$, $c_{k-1}$ and $t$ replaced by 
$r_{k+1}$,  $r_{k}$, $r_{k-1}$ and $s$, respectively. 
\end{pr}
{\it Proof.}~(i) Let $f\in B_+^k$ be arbitrary 
and $\Phi$ be as in (\ref{3.7}). Thus, 
by abuse of notation as in the previous section 
\[ 
\Phi(\xi)
=\sum_{i_1,\ldots,i_{k}(\ne)}f(z_{i_1},\ldots,z_{i_k}) 
\] 
for $(z_i)_{i=1}^{\infty}\in\Omega$ and 
$\xi=\Xi((z_i)_{i=1}^{\infty})=\sum\one_{\{z_i>0\}} \delta_{z_i}$. 
By the assumption of absolute continuity 
together with the reversibility implied by Theorem 2.4 (ii) 
\begin{eqnarray*} 
\lefteqn{
E\left[\sum_{i_1,\ldots,i_k(\ne)}
f(Z_{i_1}(t),\ldots,Z_{i_k}(t))\right]  
\ = \   
E\left[\Phi(\xi(t))\right]
\ = \
E\left[(\cT_t\Phi)(\xi(0))\right]}      \\ 
& = & 
\int_{\cN} (\cT_t\Phi)(\eta)P\left(\xi(0)\in d\eta\right)   
\ = \ 
\int_{\cN} (\cT_t\Phi)(\eta) \Psi^*(\eta) 
P\left(\sum_i \delta_{VX_i}\in d\eta\right) \\ 
& = & 
\int_{\cN} \Phi(\eta) (\cT_t\Psi^*)(\eta) 
P\left(\sum_i \delta_{VX_i}\in d\eta\right) \\ 
& = & 
\int f(\bz_k)
r_{k}(\bz_k)d\bz_k
\int_{\cN} P_{z_1,\ldots,z_k}(d\eta)
 (\cT_t\Psi^*)(\eta+\delta_{z_1}+\cdots+\delta_{z_k}),  
\end{eqnarray*} 
where the last equality is deduced from 
the Palm formula (\ref{2.8}) combined with Lemma 2.5. 
This proves (\ref{3.12}).  \\ 
(ii) As an immediate consequence of Theorem 3.2 (i) 
we have  
\[ 
\int f(\bz_k)\left[r_k(t,\bz_k)-r_k(0,\bz_k)
-\int_0^t \cL_k(s,\bz_k)ds\right]d\bz_k=0  
\] 
for all $f\in B_{+,c}^k$. 
Replace $f(\bz_k)$ by $f(\bz_k)z_1\cdots z_k$ to get 
\be 
\int f(\bz_k)\left[r_k(t,\bz_k)-r_k(0,\bz_k)
-\int_0^t \cL_k(s,\bz_k)ds\right]z_1\cdots z_kd\bz_k=0. 
                                                      \label{3.14} 
\ee 
It is easily verified from the assumption on the moments 
of $|Z(0)|$ that the signed measure 
\[   
\left[r_k(t,\bz_k)-r_k(0,\bz_k)
-\int_0^t \cL_k(s,\bz_k)ds\right]z_1\cdots z_kd\bz_k 
\]  
is expressed as a linear combination of (at most) 
8 finite measures on $(0,\infty)^k$.  Therefore, 
(\ref{3.14}) implies that it must vanish 
and accordingly  (\ref{3.13}) holds.         \qed 

\medskip 

\noindent 
{\it Example.}~ 
In the case where the lifted process 
$\sum \delta_{VX_i}$ is a gamma point process with parameter 
$(\theta,b)$ (see at the end of Section 2), 
the absolute continuity assumption in Proposition 3.5 is 
satisfied e.g. when $\xi(0)$ is a Poisson process 
on $(0,\infty)$ with mean density 
of the form $e^{h(y)}\theta e^{-by}y^{-1}$ and 
$\int|e^{h(y)}-1|e^{-by}y^{-1}dy<\infty$. 
In that case, the density $\Psi^*$ mentioned 
in Proposition 3.5 (i) is given by 
\[ 
\Psi^*(\eta)=
\exp\left[\lg h,\eta\rg-\theta\int(e^{h(y)}-1)
e^{-by}y^{-1}dy\right]. 
\] 
(See e.g. Lemma 2.4 of \cite{B}.) Also, 
since the reduced Palm distributions 
of any Poisson process are identical with its law, 
(\ref{3.12}) becomes 
\[ 
r_k(t,\bz_k)=
\frac{\theta^k}{z_1 \cdots z_k}
e^{-b(z_1+\cdots+z_k)}  
E\left[ (\cT_t\Psi^*)
(\eta+\delta_{z_1}+\cdots+\delta_{z_k})\right],
\] 
where $\eta$ is 
a gamma point process with parameter $(\theta,b)$.

\section{Preliminary results for rescaled processes}
\setcounter{equation}{0} 
\subsection{Models with a scaling parameter} 
Both this section and the subsequent section 
are devoted to a derivation of the equation (\ref{1.1}) 
from properly rescaled coagulation-fragmentation processes. 
In principle, the procedure is similar to that in \cite{EW00} 
although that paper assumes 
the conditions, among others, of the form 
\[ 
K(x,y)= o(x)o(y) \quad  \mbox{and}  \quad 
\int_0^xF(y,x-y)dy=o(x) \quad 
\mbox{as}  \quad x,y\to \infty 
\] 
for the rates. 
In our situation, these conditions are never met 
since by (\ref{1.4}) and (\ref{1.5}) 
\[ 
K(x,y)=xy(x+y)^{\lambda} 
\Hh\left(\frac{x}{x+y},\frac{y}{x+y}\right) 
\quad  \mbox{and}  \quad 
\int_0^xF(y,x-y)dy=\Cc x^{2+\lambda}. 
\] 
However, it will turn out that our setting on the degrees of 
homogeneity of $K$ and $F$ provides with us certain 
effective ingredients to overcome difficulties due to such growth orders. More specifically, it will play 
an essential role in the proof of Proposition 4.5 below. 
In this connection we mention that 
\cite{VZ} discussed the relation 
between occurrence of `steady-state solutions' 
for coagulation-fragmentation equations and 
the degrees of homogeneity. According to 
the authors' criterion based on analysis 
of the moments in several basic examples 
our setting on the degrees is in the region 
corresponding to systems for which 
steady states occur. 

Let us specify the model we are concerned with 
in the rest of the paper. Following \cite{EW00}, 
we introduce a scaling parameter $N=1,2,\ldots$ 
and modify the generator $L$ in (\ref{2.4}) 
by replacing $K$ by $K/N$. 
To be more precise, define 
\begin{eqnarray} 
L^N\Phi(\bz) 
& = & 
\frac{1}{2N}\sum_{i\ne j} z_iz_j \Hh(z_i,z_j) 
\left(\Phi(M_{ij}\bz)-\Phi(\bz)\right) \nonumber \\ 
& & 
+\frac{1}{2}\sum_i z_i 
\int_0^{z_i}dy \Hc(y,z_i-y)
\left(\Phi(S_i^{(y)}\bz)-\Phi(\bz)\right)        \label{4.1} 
\end{eqnarray} 
and denote by 
$\{Z^N(t)=(Z_i^N(t))_{i=1}^{\infty}:~t\ge 0\}$ 
the process generated by $L^N$. 
The rescaled process we will study actually is 
\[ 
\xi^N(t):=\frac{1}{N}\Xi(Z^N(t))
=\frac{1}{N}\sum_{i}\one_{(0,\infty)}(Z_i^N(t))
\delta_{Z_i^N(t)}. 
\] 
This process is regarded as a process taking values in 
$\cM$, the totality of locally finite measures on $(0,\infty)$. 
$\cM$ is equipped with the vague topology. 
For $a\ge 0$, let $\psi_a$ denote the power function 
$\psi_a(y)=y^a$, so that $\lg\psi_a, \zeta\rg$ stands for 
the `$a$th moment' of $\zeta\in \cM$. 
We consider $c_0\in\cM$ such that for some $\delta>1$  
\be 
0< \langle \psi_1, c_0 \rangle <\infty \quad 
\mbox{and} \quad 
\langle \psi_{2+\lambda+\delta}, c_0 \rangle <\infty.    \label{4.2}  
\ee 
Note that (\ref{4.2}) implies that 
$\lg\psi_a, c_0\rg<\infty$ for any $a\in[1,2+\lambda+\delta]$.  
Concerning initial distributions of $\{Z^N(t):~t\ge 0\}$, 
we suppose that  
\be 
E\left[\left|Z^N(0)\right|^{2+\lambda+\delta}\right]<\infty 
\quad \mbox{for each $N=1,2,\ldots$}.                 \label{4.3} 
\ee 
It shall be required also 
that $\xi^N(0)$ converges to $c_0$ 
in distribution as $N\to\infty$. 
A typical case where such a convergence holds is 
given in the following lemma, which is stated 
in a general setting.  
In the rest, ${\rm Po}(m)$ stands for the law of 
a Poisson point process on $(0,\infty)$ 
with mean measure $m\in\cM$. 
\begin{lm}
Let $\zeta \in \cM$ be arbitrary. 
Assume that $\sum \delta_{Y_i^N}$ is 
${\rm Po}(N\zeta)$-distributed for each $N=1,2,\ldots$.  
Then $\eta^N:=N^{-1}\sum \delta_{Y_i^N}$ 
converges to $\zeta$ in distribution as $N\to\infty$. 
Moreover, if $\lg\psi_1, \zeta\rg<\infty$ and 
$\lg\psi_a, \zeta\rg<\infty$ for some $a>1$, then 
\be 
E\left[\lg\psi_1, \eta^N\rg^a\right] 
=E\left[\left(\frac{\sum Y_i^N}{N}\right)^a\right] 
\to \lg\psi_1, \zeta\rg^a 
\quad \mbox{as} \quad N\to\infty.      \label{4.4} 
\ee 
In particular, $\sup_N
E\left[\left(N^{-1}\sum Y_i^N\right)^a\right] <\infty$. 
\end{lm}
Since the proof of Lemma 4.1 is rather lengthy and 
not relevant to other parts of this paper, 
the proof will be given in Appendix. 
It is worth noting here that requiring 
the law of $\sum \delta_{Y_i^N}$ to be Poisson 
automatically implies that every correlation measure 
of it is the direct product of the mean measure. 

By looking at the limit points of 
$\{\xi^N(t):~t\ge 0\}$ as $N\to\infty$ 
we intend to derive a weak solution to (\ref{1.1}) 
with $K$ and $F$ given by (\ref{1.4}), namely 
\begin{eqnarray} 
\lefteqn{
\frac{\partial}{\partial t}c(t,x)}                   \label{4.5}  \\ 
& = & 
\frac{1}{2}\int_0^x\left[y(x-y)\Hh(y,x-y)c(t,y)c(t,x-y)
-x \Hc(y,x-y)c(t,x)\right]dy               \nonumber  \\ 
& & 
-\int_0^{\infty}\left[xy\Hh(x,y)c(t,x)c(t,y)
-(x+y)\Hc(x,y)c(t,x+y)\right]dy.          \nonumber  
\end{eqnarray} 
One may realize that in the case 
where both $\Hh$ and $\Hc$ are constants, 
$a$ and $b$, say, respectively, 
(\ref{4.5}) can be transformed by considering 
the `size-biased version' $c^{\star}(t,x):=xc(t,x)$ into 
\begin{eqnarray*} 
x^{-1}\frac{\partial}{\partial t}c^{\star}(t,x) 
& = & 
\frac{1}{2}\int_0^x\left[a c^{\star}(t,y)c^{\star}(t,x-y)
-b c^{\star}(t,x)\right]dy                          \\ 
& & 
-\int_0^{\infty}\left[a c^{\star}(t,x)c^{\star}(t,y)
-b c^{\star}(t,x+y)\right]dy. 
\end{eqnarray*} 
This equation is very similar to (\ref{1.1}) 
with $K$ and $F$ being constants, the solution of 
which has been studied extensively in 
\cite{AB} and \cite{SD}. However, 
it is not clear whether or not there is any direct  
connection between solutions of these two equations.  
Turning to (\ref{4.5}),  the weak form 
with test functions $f\in B_{c}$ 
and initial measure $c_0$ reads 
\begin{eqnarray*} 
\lefteqn{
\int f(x)c(t,dx) -\int f(x)c_0(dx)}  \\ 
&=& 
\frac{1}{2}\int_0^tds 
\int c(s,dx)c(s,dy)xy \Hh(x,y) (\Box f)(x,y)      \\ 
& & 
-\frac{1}{2}\int_0^tds 
\int c(s,dx) x\int_0^{x}dy\Hc(y,x-y) (\Box f)(y,x-y),       
\end{eqnarray*} 
where $(\Box f)(x,y)=f(x+y)-f(x)-f(y)$.   
A rough idea for its derivation can be described as follows. 
For $f\in B_+\cup B_c$ and $\bz\in\Omega$, set 
$\Phi_f(\bz)=\sum f(z_i)$, 
adopting the convention that $f(0)=0$.  
By calculating $L^N\Phi_f$ (cf. (\ref{3.8}) with $k=1$), 
it is observed that for any $f\in B_{c}$ 
\[ 
M_f^N(t):=
\lg f, \xi^N(t) \rg -\lg f, \xi^N(0) \rg 
-\frac{1}{N}\int_0^t L^N\Phi_f(Z^N(s))ds  
\] 
is a martingale and 
\begin{eqnarray} 
\lefteqn{
\lg f,\xi^N(t)\rg -\lg f,\xi^N(0)\rg-M_f^N(t)}  \label{4.6}  \\ 
&=& 
\frac{1}{2}\int_0^tds 
\int \xi^N(s)^{[2]}(dxdy)xy \Hh(x,y) (\Box f)(x,y)       \nonumber  \\ 
& & 
-\frac{1}{2}\int_0^tds 
\int \xi^N(s)(dx)x\int_0^{x}dy \Hc(y,x-y) (\Box f)(y,x-y), \nonumber 
\end{eqnarray} 
where 
\[ 
\xi^N(s)^{[2]}
=\frac{1}{N^2}\sum_{i\ne j}
\one_{(0,\infty)}(Z_i^N(s)Z_j^N(s)) 
\delta_{(Z_i^N(s),Z_j^N(s))}. 
\] 
(By virtue of Lemma 3.4 with $k=1$, the 
integrability of $M^N_f(t)$ is ensured by (\ref{4.3}).)  
As far as the limit as $N\to\infty$ is concerned, 
$\xi^N(s)^{[2]}$ in the right side of (\ref{4.6})  
can be replaced by $\xi^N(s)^{\otimes 2}$ 
under a suitable assumption on the convergence of 
$\xi^N(0)$. Indeed, letting $R>0$ be such that 
${\rm supp}(f) \subset [0,R]$, we have by (\ref{3.4}) 
\begin{eqnarray*} 
\lefteqn{ 
\left|\int \left(\xi^N(s)^{[2]}
-\xi^N(s)^{\otimes 2}\right)(dxdy) 
xy \Hh(x,y) (\Box f)(x,y) \right|}               \\ 
& = & 
\left|\frac{1}{N^2}
\sum_{i} Z_i^N(s)^2 \Hh(Z_i^N(s),Z_i^N(s))  \left\{f(2Z_i^N(s))-2f(Z_i^N(s))\right\}\right| \\ 
& \le & 
\frac{\Ch}{N^2}
\sum_{i} Z_i^N(s)^2 (2Z_i^N(s))^{\lambda}
\left|f(2Z_i^N(s))-2f(Z_i^N(s))\right|              \\ 
& \le & 
\frac{\Ch R(2R)^{\lambda}}{N^2} 
\sum_{i} Z_i^N(s) \cdot 3\Vert f \Vert_{\infty} \\ 
& = & \frac{3\Ch2^{\lambda}R^{1+\lambda}}{N} 
\Vert f \Vert_{\infty} \lg \psi_1,\xi^N(0)\rg. 
\end{eqnarray*} 
Thus, at least at formal level, 
the derivation of a weak solution to (\ref{4.5}) 
from $\xi^N(t)$ would reduce to proving that 
$M_f^N(t)$ vanishes 
in a suitable sense as $N\to \infty$.  \par 
To this end, we shall calculate the 
quadratic variation $\lg M_f^N\rg (t)$. 
However, the square integrability of $M_f^N(t)$ 
is nontrivial under the condition (\ref{4.3}). 
We will guarantee this in the next lemma 
by the cutoff argument as in the proof of Theorem 3.2. 
Also, the following inequality will be used to bound 
the quadratic variation: for any $s,t\ge 0$ and $a\ge 0$ 
\be   
(s+t)^a \le C_{1,a}(s^a+t^a),              \label{4.7}             
\ee 
where $C_{1,a}=2^{a-1}\vee 1$. 
This inequality for $a>1$ is deduced from 
H\"older's inequality and the one for $0<a\le 1$ is 
implied by the identity 
\be  
(s+t)^a-s^a-t^a
=a(a-1)\int_0^sdu\int_0^tdv(u+v)^{a-2}.      \label{4.8}    
\ee 
\begin{lm}
Suppose that 
$E\left[\left|Z^N(0)\right|^{2+\lambda}\right]<\infty$. 
For each $f\in B_{c}$,  
$\{M_f^N (t):~t\ge 0\}$ is a square integrable 
martingale with quadratic variation 
\[ 
\lg M_f^N\rg (t) 
=
\int_0^t \Gamma^N_f(Z^N(s))ds,   
\] 
where  $\Gamma^N_f: \Omega\to \R_+$ is defined to be 
\begin{eqnarray} 
\Gamma^N_f(\bz)    
& = & 
\frac{1}{N^2}\left\{L^N((\Phi_f)^2)(\bz)
-2\Phi_f(\bz)L^N\Phi_f(\bz)\right\}  \nonumber     \\ 
& = & 
\frac{1}{2N^3}\sum_{i\ne j} z_iz_j \Hh(z_i,z_j) 
\left\{f(z_i+z_j)-f(z_i)-f(z_j)\right\}^2            \label{4.9}    \\         
& & 
+\frac{1}{2N^2}\sum_i z_i^{2} 
\int_0^{1}du \Hc(uz_i,(1-u)z_i)
\left\{f(uz_i)+f((1-u)z_i)-f(z_i)\right\}^2      \nonumber      \\ 
& =: & 
\frac{1}{2}\Sigma_K^N(\bz)
+\frac{1}{2}\Sigma_F^N(\bz).          \nonumber 
\end{eqnarray} 
Moreover, the following estimates hold:  
\be 
\Sigma_K^N(\bz)
\le 18\Ch C_{1,\lambda}
\frac{\Vert f \Vert_{\infty}^2}{N}\cdot\frac{|\bz|}{N}
\sum_i\frac{z_i^{1+\lambda}}{N}
= 18\Ch C_{1,\lambda}\frac{\Vert f \Vert_{\infty}^2}{N}
\langle \psi_1,\xi^N \rangle    
\langle \psi_{1+\lambda},\xi^N \rangle         \label{4.10}  
\ee 
and 
\be 
\Sigma_F^N(\bz)
\le 9\Cc\frac{\Vert f \Vert_{\infty}^2}{N} 
\sum_i\frac{z_i^{2+\lambda}}{N}  
= 9\Cc\frac{\Vert f \Vert_{\infty}^2}{N}  
\langle \psi_{2+\lambda},\xi^N \rangle,         \label{4.11} 
\ee 
where $\xi^N=N^{-1}\Xi(\bz)
=N^{-1}\sum \one_{\{z_i>0\}}\delta_{z_i}$.  
\end{lm}
{\it Proof.}~ 
For any $R>0$, let 
$\Phi_f^{(R)}(\bz)=\Phi_f (\bz)\one_{\{|\bz|\le R\}}$ 
and observe that 
\begin{eqnarray*} 
M_f^{N,R}(t) 
& := & 
\frac{1}{N}\Phi_f^{(R)}(Z^N(t))
-\frac{1}{N}\Phi_f^{(R)}(Z^N(0))
-\frac{1}{N}\int_0^t L^N\Phi_f^{(R)}(Z^N(s))ds   \\ 
& = & 
M_f^{N}(t) \one_{\{|Z^N(0)|\le R\}} 
\end{eqnarray*} 
is a bounded martingale with quadratic variation 
\begin{eqnarray*} 
\lg M_f^{N,R} \rg (t) 
& = & 
\frac{1}{N^2} \int_0^t 
\left[L^N((\Phi_f^{(R)})^2)(Z^{N}(s))
-2\Phi_f^{(R)}(Z^{N}(s)) L^N\Phi_f^{(R)}(Z^{N}(s)) \right]ds \\ 
& = & 
\int_0^t \Gamma_f^N(Z^N(s))ds 
\one_{\{|Z^N(0)|\le R\}}. 
\end{eqnarray*} 
The expression (\ref{4.9}) for $\Gamma_f^N(\bz)$ 
is deduced from (\ref{4.1}). By (\ref{3.4}) and (\ref{4.7}) 
\begin{eqnarray*} 
\Sigma_K^N(\bz)           
& \le  & 
9\Ch C_{1,\lambda}\Vert f \Vert_{\infty}^2\sum_{i\ne j}
\frac{z_iz_j }{N^3} (z_i^{\lambda}+z_j^{\lambda})    \\ 
& = & 
18\Ch C_{1,\lambda}\Vert f \Vert_{\infty}^2\sum_{i\ne j}
\frac{z_i^{1+\lambda}z_j }{N^3}                 \\ 
& \le  & 
18\Ch C_{1,\lambda}\frac{\Vert f \Vert_{\infty}^2}{N}
\cdot\frac{|\bz|}{N}\sum_i\frac{z_i^{1+\lambda}}{N}.  
\end{eqnarray*} 
This proves (\ref{4.10}), whereas 
(\ref{4.11}) is immediate from (\ref{3.6}). 
The assumption together with these two estimates 
ensure the integrability of 
$\int_0^t\Gamma_f^N(Z^N(s))ds$. Therefore, 
by the monotone convergence theorem 
\begin{eqnarray*} 
E\left[(M_f^{N}(t))^2 \right] 
& = & 
\lim_{R\to\infty}E\left[(M_f^{N,R}(t))^2 \right]      \\ 
& = & 
\lim_{R\to\infty}
E\left[\int_0^t \Gamma_f^N(Z^N(s))ds 
\one_{\{|z|\le R\}}\right]                         \\ 
& =  & 
E\left[\int_0^t \Gamma_f^N(Z^N(s))ds \right]<\infty.    
\end{eqnarray*} 
Once the square integrability of $M_f^{N}(t)$ 
is in hand, 
one can show further by a similar argument that 
$(M_f^{N}(t))^2-\int_0^t \Gamma_f^N(Z^N(s))ds$ is 
a martingale. The proof of Lemma 4.2 is complete. 
\qed 

\medskip 

\noindent 
{\it Remarks.}~ (i) A heuristic derivation of (\ref{1.1}) 
based on the hierarchical structure discussed 
in the previous section is available: 
due to the rescaling $K \mapsto K/N$ and 
$\xi \mapsto \xi^N :=\xi /N$, the equation 
solved weakly by the family 
$\{r_k^N(t,d\bz_k):~t\ge0, k\in\N\}$ 
of correlation measures of $\xi^N(t)$ becomes 
\begin{eqnarray*} 
\lefteqn{
\frac{\partial}{\partial t}r_k^N(t,z_1,\ldots,z_k)} \\ 
&=& 
\frac{1}{2}\sum_{l=1}^k\int_0^{z_l}K(y,z_l-y)
r^N_{k+1}(t,z_1,\ldots,z_{l-1},y,z_l-y,z_{l+1},\ldots,z_k)dy \\ 
& & 
-\frac{1}{2}\sum_{l=1}^k\int_0^{z_l}F(y,z_l-y)dy \ 
r^N_{k}(t,z_1,\ldots,z_k)                           \\ 
& & 
-\sum_{l=1}^k\int_0^{\infty}K(z_l,y)
r^N_{k+1}(t,z_1,\ldots,z_l,y,z_{l+1},\ldots,z_k)dy    \\ 
& & 
+\sum_{l=1}^k\int_0^{\infty}F(z_l,y)
r^N_{k}(t,z_1,\ldots,z_{l-1},z_l+y,z_{l+1},\ldots,z_k)dy    \\ 
& &  
-\frac{\one_{\{k\ge 2\}}}{N}
\sum_{l<m}^{k}K(z_l,z_m)r^N_{k}(t,z_1,\ldots,z_k)  \\ 
& &  
+\frac{\one_{\{k\ge 2\}}}{N}\sum_{l<m}^{k}F(z_l,z_m)
r^N_{k-1}(t,z_1,\ldots,z_{l-1},z_l+z_m,z_{l+1},\ldots, 
z_{m-1},z_{m+1},\ldots,z_k). 
\end{eqnarray*} 
So, the last two terms would be expected to 
vanish in the limit as $N\to \infty$ and 
the limits $r_k$ of $r_k^N$, if they exist, 
would solve the same equations 
as the ones satisfied by the direct products 
$c^{\otimes k}(t,z_1,\ldots,z_k)
=c(t,z_1)\cdots c(t,z_k)$ of a solution to (\ref{1.1}). 
This procedure has been accomplished in 
\cite{EP} for a pure coagulation model.  \\ 
(ii) We also give a remark 
on the asymptotic equivalence 
between the moment measures and 
the correlation measures of the rescaled process 
$\xi^N(t)$. (cf. Lemma 1.16 in \cite{Kol}. 
The reader is cautioned that our terminology 
`moment measure' is in conflict with that of \cite{Kol}.) 
For each $k\in\{2,3,\ldots\}$, 
under the assumption that 
$\sup_N E[(|Z^N(0)|/N)^{k-1}]<\infty$, 
it holds that for any $f\in B_{+,c}^k$ 
\[ 
\left|E\left[\frac{1}{N^k}\sum_{i_1,\ldots,i_k}
f(Z_{i_1}^N(t),\ldots,Z_{i_k}^N(t))\right] 
-E\left[\frac{1}{N^k}\sum_{i_1,\ldots,i_k(\ne)}
f(Z_{i_1}^N(t),\ldots,Z_{i_k}^N(t))\right]\right|
\le \frac{C}{N} 
\] 
for some constant $C$ independent of $N$.  
Indeed, it is sufficient to verify this by assuming that 
$f(z_1,\ldots,z_k)=
z_1\cdots z_k\one_{[\epsilon,R]}(z_1)\cdots 
\one_{[\epsilon,R]}(z_k)$ with $0<\epsilon<R$, 
for which the above difference is dominated by 
a finite sum of expectations of the form 
\[ 
E\left[\frac{1}{N^k}\sum_{i_1,\ldots,i_j} 
Z_{i_1}^N(t)^{n_1}\one_{[\epsilon,R]}(Z_{i_1}^N(t))
\cdots 
Z_{i_j}^N(t)^{n_j}\one_{[\epsilon,R]}(Z_{i_j}^N(t))\right], 
\] 
where $j\in\{1,\ldots,k-1\}$, $n_1,\ldots,n_j\in \N$ 
is such that $n_1+\cdots+n_j=k$ 
and hence $n_{l}\ge 2$ for some $l\in\{1,\ldots,j\}$. 
The desired bound follows by noting that 
\[ 
Z_{i_l}^N(t)^{n_l}\one_{[\epsilon,R]}(Z_{i_l}^N(t)) 
\le 
R Z_{i_l}^N(t)^{n_l-1}\one_{[\epsilon,R]}(Z_{i_l}^N(t)) 
\] 
and observing that the above expectation is 
less than or equal to 
\[ 
E\left[\frac{R}{N^k}|Z^N(t)|^{k-1}\right]
=\frac{R}{N}
E\left[\left(\frac{|Z^N(0)|}{N}\right)^{k-1}\right]. 
\] 

\subsection{Key estimates for the martingale}
Lemma 4.2 implies that for some constant 
$C_2>0$ independent of $N, t$ and $f$  
\begin{eqnarray} 
\lefteqn{
C_2 \Vert f \Vert_{\infty}^{-2} 
 E\left[\lg M_f^N\rg (t)\right]}         \nonumber     \\ 
& \le &  
\frac{1}{N}E\left[\frac{|Z^N(0)|}{N}\int_0^t ds
\sum_i\frac{Z^N_i(s)^{1+\lambda}}{N} \right] 
+ \frac{1}{N}E\left[\int_0^t ds \sum_i
\frac{Z^N_i(s)^{2+\lambda}}{N} \right].       \label{4.12} 
\end{eqnarray} 
In fact, the two expectations on (\ref{4.12}) are 
related to each other in such a way that 
\begin{eqnarray} 
\lefteqn{
E\left[\frac{|Z^N(0)|}{N}\int_0^t ds 
\sum_i\frac{Z^N_i(s)^{1+\lambda}}{N} \right]}  \label{4.13}   \\ 
& \le &  
\left(t
E\left[\left(\frac{|Z^N(0)|}{N}\right)^{2+\lambda}\right]
\right)^{\frac{1}{1+\lambda}}
\left(E\left[\int_0^t ds \sum_i 
\frac{Z^N_i(s)^{2+\lambda}}{N}\right]
\right)^{\frac{\lambda}{1+\lambda}},       \nonumber 
\end{eqnarray} 
which is a special case 
($\alpha=1, \gamma=\lambda$ and $\epsilon=1$)  
of the following lemma. 
\begin{lm}
For arbitrary $\alpha, \gamma \ge 0$ and $\epsilon>0$  
\begin{eqnarray} 
\lefteqn{
E\left[\left(\frac{|Z^N(0)|}{N}\right)^{\alpha}\int_0^t ds
\sum_i\frac{Z^N_i(s)^{1+\gamma}}{N} \right]}    \label{4.14}   \\ 
& \le &  
\left(t
E\left[\left(\frac{|Z^N(0)|}{N}
\right)^{1+\alpha+\frac{\alpha\gamma}{\epsilon}}\right]
\right)^{\frac{\epsilon}{\gamma+\epsilon}}
\left(E\left[\int_0^t ds \sum_i 
\frac{Z^N_i(s)^{1+\gamma+\epsilon}}{N}\right]
\right)^{\frac{\gamma}{\gamma+\epsilon}}.       \nonumber 
\end{eqnarray} 
\end{lm}
{\it Proof.}~
Since (\ref{4.14}) for $\gamma=0$ is obvious, 
we may assume that $\gamma>0$.  
Then put $p=1+\gamma/\epsilon$ 
and $q=1+\epsilon/\gamma$, which are mutually conjugate. 
By virtue of H\"older's inequality with respect to 
the `weight' $E[N^{-1}\int_0^tds \sum_i \cdot ]$ 
\begin{eqnarray*} 
\lefteqn{
E\left[\left(\frac{|Z^N(0)|}{N}\right)^{\alpha}\int_0^t ds
\sum_i\frac{Z^N_i(s)^{1+\gamma}}{N} \right]}     \\ 
& = & 
E\left[\int_0^t \frac{ds}{N}\sum_i
\left(\frac{|Z^N(0)|^{\alpha}}{N^{\alpha}}
Z^N_i(s)^{\frac{1}{p}}\right)
Z^N_i(s)^{1+\gamma-\frac{1}{p}}\right]   \\ 
& \le &  
\left(E\left[\int_0^t \frac{ds}{N}\sum_i
\left(\frac{|Z^N(0)|}{N}\right)^{\alpha p}Z^N_i(s)\right]
\right)^{\frac{1}{p}}
\left(E\left[\int_0^t \frac{ds}{N} \sum_i 
Z^N_i(s)^{(1+\gamma-\frac{1}{p})q}\right]
\right)^{\frac{1}{q}}   \\ 
& = &  
\left(t
E\left[\left(\frac{|Z^N(0)|}{N}\right)^{\alpha p+1}\right]
\right)^{\frac{\epsilon}{\gamma+\epsilon}}
\left(E\left[\int_0^t ds \sum_i 
\frac{Z^N_i(s)^{1+\gamma+\epsilon}}{N}\right]
\right)^{\frac{\gamma}{\gamma+\epsilon}}.      
\end{eqnarray*} 
Here, the final equality is due to 
$(1+\gamma-\frac{1}{p})q=1+\gamma+\epsilon$. 
(\ref{4.14}) has been obtained  
since $\alpha p+1=1+\alpha+\alpha\gamma/\epsilon$. 
\qed 

\medskip 

The expression in the right side of (\ref{4.13}) 
motivates us to define 
\[ 
\overline{m}_{a}= 
\sup_N E\left[\left(\frac{|Z^N(0)|}{N}\right)^{a}\right] 
=\sup_N E\left[\lg \psi_1,\xi^N(0)\rg^a \right] 
\] 
for $a\ge 0$. We will discuss under the condition that 
\be 
\overline{m}_{2+\lambda+\delta}< \infty.     \label{4.15}  
\ee 
This condition is stronger than (\ref{4.3}) and 
valid in the case of Poisson processes 
considered in Lemma 4.1. So, we shall focus 
attention on estimation of the second term of (\ref{4.12}). 
Before doing it, we prepare an elementary but 
technically important lemma. 
\begin{lm}
Let $A,B>0$, a real number $C$ and $q>1$ be given. 
If $A$ satisfies a `self-dominated inequality' 
$A \le B A^{1/q}+ C$, then 
\[ 
A 
\le B^{\frac{q}{q-1}}+\frac{q}{q-1}(C \vee 0). 
\]  
\end{lm}
{\it Proof.}~ Since the assertion is obvious 
when $C\le 0$, we only have to consider the case $C>0$.  
Denote by $\wt{A}$ the right side of the required inequality. 
The equation $x=B x^{1/q}+C$ for $x>0$ 
has a unique root and 
$x>B x^{1/q}+C$ for sufficiently large $x$. 
Therefore, it is enough to show that 
$\wt{A} > B \wt{A}^{1/q}+C.$ 
Observe that for $s,t>0$ 
\[ 
(s+t)^{\frac{1}{q}}-s^{\frac{1}{q}} 
\ = \ \frac{1}{q} \int_0^tdu(s+u)^{\frac{1}{q}-1}
\ < \ \frac{1}{q}  \int_0^tdu s^{\frac{1}{q}-1} 
\ = \  \frac{1}{q}s^{\frac{1}{q}-1}t 
\] 
and thus $(s+t)^{1/q}<s^{1/q} + s^{(1/q)-1}t/q$. 
Plugging $s=B^{\frac{q}{q-1}}$ 
and $t=\frac{q}{q-1}C$ gives 
\begin{eqnarray*} 
B\wt{A}^{1/q}+C 
& < & 
B\left(B^{\frac{1}{q-1}}
+B^{-1}\frac{q}{q-1}C\cdot\frac{1}{q}\right)+C   \\ 
& = &  
B^{\frac{q}{q-1}}+\frac{q}{q-1}C   
\ = \ 
\wt{A}   
\end{eqnarray*} 
as desired. 
\qed 

\medskip 
\noindent 
The next proposition supplies the key to proceeding further.   
\begin{pr}
Assume that (\ref{4.2}) holds for some $\delta>1$. 
Suppose that $\sum \delta_{Z_i^N(0)}$ 
has mean measure $Nc_0$ for each $N=1,2,\ldots$ 
and that (\ref{4.15}) holds. 
Then, for each $t>0$,  there exist constants 
$\overline{C}(t)$, $\overline{C}_0(t)$ and $C^*(t)$ 
independent of $N$ such that 
\be 
E\left[\int_0^t ds 
\langle \psi_{2+\lambda+\delta},\xi^N(s) \rangle \right] 
= 
E\left[\int_0^t ds \sum_i
\frac{Z^N_i(s)^{2+\lambda+\delta}}{N} \right]   
\le \overline{C}(t)                                \label{4.16}  
\ee  
\be 
E\left[\int_0^t ds 
\langle \psi_{2+\lambda},\xi^N(s) \rangle \right] 
= 
E\left[\int_0^t ds \sum_i
\frac{Z^N_i(s)^{2+\lambda}}{N} \right]  
\le \overline{C}_0(t)                                 \label{4.17}  
\ee 
and 
\be  
E\left[\langle \psi_{\delta},\xi^N(t) \rangle \right] 
= 
E\left[\sum_i
\frac{Z^N_i(t)^{\delta}}{N} \right]  \le C^*(t).    \label{4.18}                  
\ee 
\end{pr}
{\it Proof.}~ 
Define bounded functions $\Psi_{\delta}^{(R)}(\bz)= 
\sum_i z_i^{\delta}\one_{\{|\bz|\le R\}}$ on $\Omega$ 
for all $R>0$. By direct calculations 
\begin{eqnarray*} 
L^N\Psi_{\delta}^{(R)}(\bz) 
& = & 
\frac{1}{2N}\sum_{i\ne j}z_iz_j \Hh(z_i,z_j) 
\left\{(z_i+z_j)^{\delta}-z_i^{\delta}-z_j^{\delta}\right\}
\one_{\{|\bz|\le R\}}                         \\ 
&  & 
+\frac{1}{2}\sum_i z_i^{2+\lambda+\delta}
\int_0^{1}du \Hc(u,1-u) 
\left\{u^{\delta}+(1-u)^{\delta}-1\right\}     
\one_{\{|\bz|\le R\}}. 
\end{eqnarray*} 
A crucial point here is that by $\delta>1$ and (\ref{4.8}) 
\[ 
\check{C}_{\delta}:=
\int_0^{1}du \Hc(u,1-u)\left\{1-u^{\delta}-(1-u)^{\delta}\right\}\in (0,\Cc). 
\] 
Taking expectation of the martingale 
\[ 
\Psi_{\delta}^{(R)}(Z^N(t))
-\Psi_{\delta}^{(R)}(Z^N(0))
-\int_0^tds L^N\Psi_{\delta}^{(R)}(Z^N(s)) 
\] 
and then letting $R\to \infty$ yield 
\begin{eqnarray*}  
\lefteqn{
\frac{\check{C}_{\delta}}{2}E\left[\int_0^t ds 
\sum_i Z^N_i(s)^{2+\lambda+\delta} \right]
+E\left[\sum_i Z^N_i(t)^{\delta} \right]} \\ 
& = & 
E\left[\sum_i Z^N_i(0)^{\delta} \right] 
+ \frac{1}{2N}E\left[\int_0^t ds 
\sum_{i \ne j} Z^N_i(s)Z^N_j(s)
\Hh(Z^N_i(s),Z^N_j(s))  \right.         \\ 
& & 
\hspace*{5cm} 
\left. \cdot \left\{(Z^N_i(s)+Z^N_j(s))^{\delta}
-Z^N_i(s)^{\delta}-Z^N_j(s)^{\delta}\right\}  \right]  \\ 
& \le & 
N\lg \psi_{\delta}, c_0\rg 
+ \frac{\Ch}{2N}E\left[\int_0^t ds 
\sum_{i \ne j} Z^N_i(s)Z^N_j(s)
(Z^N_i(s)+Z^N_j(s))^{\lambda+\delta} \right]  \\ 
& \le & 
N\lg \psi_{\delta}, c_0\rg 
+ \frac{\Ch C_{1,\lambda+\delta}}{2N}
E\left[\int_0^t ds 
\sum_{i \ne j} Z^N_i(s)Z^N_j(s)
\left\{Z^N_i(s)^{\lambda+\delta}
+Z^N_j(s)^{\lambda+\delta}\right\} \right]  \\ 
& \le & 
N\lg \psi_{\delta}, c_0\rg 
+ \frac{\Ch C_{1,\lambda+\delta}}{N}
E\left[|Z^N(0)| \int_0^t ds 
\sum_{i} Z^N_i(s)^{1+\lambda+\delta} \right]. 
\end{eqnarray*} 
Note that all the expectations in the above are finite 
by (\ref{4.15}). We have obtained 
\begin{eqnarray}  
\lefteqn{
\frac{\check{C}_{\delta}}{2}E\left[\int_0^t ds 
\sum_i \frac{Z^N_i(s)^{2+\lambda+\delta}}{N} \right]
+E\left[\sum_i \frac{Z^N_i(t)^{\delta}}{N} \right]}  \nonumber \\ 
& \le & 
\Ch C_{1,\lambda+\delta}
E\left[\frac{|Z^N(0)|}{N} \int_0^t ds  
\sum_{i} \frac{Z^N_i(s)^{1+\lambda+\delta}}{N} \right]
+\lg \psi_{\delta}, c_0\rg .                             \label{4.19}
\end{eqnarray} 
In addition, Lemma 4.3 with 
$\alpha=1, \gamma=\lambda+\delta$ and $\epsilon=1$ reads 
\begin{eqnarray} 
\lefteqn{
E\left[\frac{|Z^N(0)|}{N}\int_0^t ds 
\sum_i\frac{Z^N_i(s)^{1+\lambda+\delta}}{N} \right]}  \nonumber \\ 
& \le &  
\left(t
E\left[\left(\frac{|Z^N(0)|}{N}\right)^{2+\lambda+\delta}\right]
\right)^{\frac{1}{1+\lambda+\delta}}
\left(E\left[\int_0^t ds \sum_i 
\frac{Z^N_i(s)^{2+\lambda+\delta}}{N}\right]
\right)^{\frac{\lambda+\delta}{1+\lambda+\delta}}.   \label{4.20}   
\end{eqnarray} 
Set $p=1+\lambda+\delta$ and $q=1+1/(\lambda+\delta)$. 
By combining the above two inequalities 
\[ 
A_N\ := \ E\left[\int_0^t ds 
\sum_i \frac{Z^N_i(s)^{2+\lambda+\delta}}{N} \right]
\le 
B_N\left(A_N\right)^{1/q} 
+\frac{2\lg \psi_{\delta}, c_0\rg }{\check{C}_{\delta}}, 
\] 
where 
\[ 
B_N 
=\frac{2\Ch C_{1,\lambda+\delta}}{\check{C}_{\delta}}
\left(t
E\left[\left(\frac{|Z^N(0)|}{N}\right)^{2+\lambda+\delta}\right]
\right)^{1/p}. 
\] 
Applying Lemma 4.4 and (\ref{4.15}) and noting that 
$\frac{q}{q-1}=p$, we obtain 
\begin{eqnarray*} 
A_N 
& \le &  
\left(B_N\right) ^{\frac{q}{q-1}}+ \frac{q}{q-1}\cdot
\frac{2\lg \psi_{\delta}, c_0\rg }{\check{C}_{\delta}}             \\ 
& \le & 
\left(\frac{2\Ch C_{1,\lambda+\delta}}{\check{C}_{\delta}}
\right)^{1+\lambda+\delta} 
\overline{m}_{2+\lambda+\delta} t  
+\frac{2(1+\lambda+\delta)\lg \psi_{\delta}, c_0\rg }
{\check{C}_{\delta}} \\ 
& =: & 
\overline{C}(t). 
\end{eqnarray*} 
This proves (\ref{4.16}). The estimate (\ref{4.17}) 
can be deduced from (\ref{4.16}), 
Lemma 4.3 ($\alpha=0, \gamma=1+\lambda$ and $\epsilon=\delta$) and (\ref{4.15}). 
(\ref{4.18}) follows from (\ref{4.19}), (\ref{4.20}), 
(\ref{4.15}) and (\ref{4.16}). 
We complete the proof of Proposition 4.5.  \qed 
\begin{cor}
Under the same assumptions as in Proposition 4.5, 
for each $t>0$, 
there exists a constant $C(t)$ independent of 
$N\in\N$ and $f\in B_{c}$ such that 
\be 
E\left[\lg M_f^N\rg (t)\right]  
\le 
\frac{\Vert f \Vert_{\infty}^{2}}{N}C(t).    \label{4.21} 
\ee 
\end{cor}
{\it Proof.}~This is immediate from (\ref{4.12}), 
(\ref{4.13}), (\ref{4.15}) and (\ref{4.17}).  
\qed 

\medskip 

We end this section with a lemma 
which will be used in the next section.  Recall that 
for $\gamma>0$ the function $\Psi_{\gamma}$ 
on $\Omega$ is defined by  
$\Psi_{\gamma}(\bz)=\sum z_i^{\gamma}$. 
\begin{lm}
Let $\bz\in\Omega$ and put 
$\xi^N=N^{-1}\Xi(\bz)
=N^{-1}\sum \one_{\{z_i>0\}}\delta_{z_i}$. \\ 
(i) For any $f\in B_c$ 
\begin{eqnarray} 
\lefteqn{
\frac{1}{N}\left|L^N\Phi_f(\bz)\right|} \nonumber \\    
& \le &  
\Ch C_{1,\lambda} \Vert f \Vert_{\infty}   
\lg \psi_1,\xi^N \rg    
\lg \psi_{1+\lambda},\xi^N \rg 
+\frac{\Cc}{2} \Vert f \Vert_{\infty}
\lg \psi_{2+\lambda},\xi^N \rg        \label{4.22}     \\ 
& \le &  
\Ch C_{1,\lambda} \Vert f \Vert_{\infty}   
\lg \psi_1,\xi^N \rg^{\frac{2+\lambda}{1+\lambda}}    
\lg \psi_{2+\lambda},\xi^N \rg^{\frac{\lambda}{1+\lambda}}    
+\frac{\Cc}{2} \Vert f \Vert_{\infty}
\lg \psi_{2+\lambda},\xi^N \rg            \label{4.23}  
\end{eqnarray} 
provided that $\lg\psi_{2+\lambda},\xi^N\rg<\infty$. \\ 
(ii) Let $\gamma>1$ and assume that 
$E\left[\left|Z^N(0)\right|
^{2+\lambda+\gamma}\right]<\infty$. Then 
\be  
\wt{M}_{\gamma}^N(t):= 
\frac{1}{N}\Psi_{\gamma}(Z^N(t))
-\frac{1}{N}\Psi_{\gamma}(Z^N(0))
-\frac{1}{N}\int_0^tL^N\Psi_{\gamma}(Z^N(s))ds 
                                                    \label{4.24} 
\ee 
is a martingale. Moreover, 
\be 
\frac{1}{N}\left|L^N\Psi_{\gamma}(\bz)\right| 
\le 
\Ch C_{1,\lambda+\gamma}
\lg \psi_1,\xi^N \rg    
\lg \psi_{1+\lambda+\gamma},\xi^N \rg 
+\frac{\check{C}_{\gamma}}{2} 
\lg \psi_{2+\lambda+\gamma},\xi^N \rg   \label{4.25} 
\ee 
provided that $\lg\psi_{2+\lambda+\gamma},\xi^N\rg<\infty$. 
\end{lm}
{\it Proof.}~(i) (\ref{4.22}) can be shown in just a similar 
way to (\ref{4.10}) and (\ref{4.11}) in view of (\ref{4.9}) and 
\begin{eqnarray*} 
\frac{1}{N}L^N\Phi_f(\bz)    
& = & 
\frac{1}{2N^2}\sum_{i\ne j} z_iz_j \Hh(z_i,z_j) 
\left\{f(z_i+z_j)-f(z_i)-f(z_j)\right\}                \\         
& & 
+\frac{1}{2N}\sum_i z_i^{2} 
\int_0^{1}du \Hc(uz_i,(1-u)z_i)
\left\{f(uz_i)+f((1-u)z_i)-f(z_i)\right\}. 
\end{eqnarray*} 
(\ref{4.23}) is a consequence of H\"older's inequality. \\ 
(ii) The proof is quite analogous to that of Lemma 4.2 
in view of calculations at the beginning of 
the proof of Proposition 4.5 
with $\gamma$ in place of $\delta$. 
So the details are left to the reader.   \qed 

\section{Derivation of the macroscopic equation}
\setcounter{equation}{0} 
\subsection{Tightness arguments}
Before studying the limit of the rescaled processes discussed 
in Section 4, relative compactness of their laws must be argued. 
In fact, it will be convenient to consider, 
rather than $\xi^N(t)$, the measure-valued process 
\[ 
\mu^N(t):= \frac{1}{N}\sum_i Z^N_i(t)\delta_{Z^N_i(t)}, 
\] 
which takes values in $\cM_f$, 
the space of finite measures on $(0,\infty)$, almost surely  
as long as $P(|Z^N(0)|<\infty)=1$. 
Denote by $C_c$  
the set of continuous functions on $(0,\infty)$ 
with compact support and set 
$C_{+,c}=B_+\cap C_{c}$. 
We begin the tightness argument by introducing 
a metric on $\cM_f$ compatible with the weak topology by 
\[ 
d_w(\nu,\nu')
=\sum_{k=0}^{\infty}2^{-k}
\left(|\lg h_k, \nu\rg -\lg h_k, \nu'\rg|\wedge 1\right), 
\qquad \nu,\nu' \in\cM_f 
\] 
where $h_0\equiv 1$ and $\{h_1,h_2,\ldots\}\subset C_{+,c}$ 
is as in the proof of Proposition 3.17 of \cite{Re}.   
(Alternatively, see A 7.7 of \cite{Ka}.)  
In particular, denoting the vague convergence 
by $\stackrel{v}{\to}$, we have, 
for $\eta, \eta_1,\eta_2,\cdots \in\cM$, 
$\eta_n \stackrel{v}{\to}\eta$ iff  
$\lg h_k, \eta_n\rg \to \lg h_k,\eta\rg $ 
for all $k\in\N$. 
$(\cM_f, d_w)$ is a complete, separable metric space. 
Let $\stackrel{w}{\to}$ stand for the weak convergence. 
Given $a,b, \gamma>0$ arbitrarily, let 
\[ 
\cM_{a,b}^{\gamma}
=\{\nu\in\cM_f :~ \lg 1, \nu\rg\le a, 
\lg \psi_{\gamma}, \nu\rg\le b\}. 
\] 
As will be seen in the next two lemmas, 
it is a compact set in $\cM_f$ 
and plays an important role. 
We introduce an auxiliary function $\varphi_R$ 
for each $R>0$ by 
\[ \varphi_R(y)
=\left\{
\begin{array}{ll} 
1 & (y\le R) \\  
R+1-y & (R\le y \le R+1) \\  
0 & (y\ge R+1), 
\end{array} \right. 
\] 
which is bounded and continuous on $(0,\infty)$. 
\begin{lm}
Let $\cM_{a,b}^{\gamma}$ be as above and 
$\nu_1,\nu_2,\ldots\in \cM_{a,b}^{\gamma}$. 
Assume that $\nu_n\stackrel{w}{\to}\nu$ 
for some $\nu\in\cM_f$. 
Then \\ 
(i) $\nu \in \cM_{a,b}^{\gamma}$. \\
(ii) For any $\alpha\in(0,\gamma)$, 
$\lg \psi_{\alpha}, \nu_n\rg \to \lg \psi_{\alpha}, \nu\rg$ 
as $n\to \infty$.  
\end{lm}
{\it Proof }~
This is a special case of Lemma 4.1 in \cite{EW00}. 
But we give a proof for completeness.  \\ 
(i) It is clear that  
$\lg 1, \nu\rg=\lim_{n\to\infty}\lg 1, \nu_n\rg\le a$. 
Also, 
\[
\lg \psi_{\gamma}\varphi_R, \nu\rg
=\lim_{n\to\infty}\lg \psi_{\gamma}\varphi_R, \nu_n\rg 
\le \liminf_{n\to\infty}\lg \psi_{\gamma}, \nu_n\rg\le b. 
\] 
Letting $R\to\infty$, we get 
$\lg \psi_{\gamma}, \nu\rg \le b$. 
Hence $\nu\in \cM_{a,b}^{\gamma}$. \\ 
(ii) Fix $\alpha\in(0,\gamma)$ arbitrarily. 
Observe that for each $R>0$ 
\begin{eqnarray*} 
|\lg \psi_{\alpha},\nu_n \rg- \lg \psi_{\alpha},\nu \rg| 
& \le & 
|\lg \psi_{\alpha}\varphi_R,\nu_n \rg
- \lg \psi_{\alpha}\varphi_R,\nu \rg|  \\ 
& & 
+\lg \psi_{\alpha}(1-\varphi_R),\nu_n \rg
+\lg \psi_{\alpha}(1-\varphi_R),\nu \rg.  
\end{eqnarray*} 
The first term on the right side converges to 0 
as $n\to \infty$. As for the second and 
third terms, we have a uniform bound in $n$ 
\[ 
\lg \psi_{\alpha}(1-\varphi_R),\nu_n \rg 
\le\lg \frac{\psi_{\alpha}}{\psi_{\gamma}} \cdot 
\one_{[R,\infty)} \cdot \psi_{\gamma}, \nu_n \rg 
\le R^{\alpha-\gamma}\lg \psi_{\gamma}, \nu_n\rg 
\le bR^{\alpha-\gamma}
\] 
and similarly 
$\lg \psi_{\alpha}(1-\varphi_R),\nu \rg  \le 
bR^{\alpha-\gamma}$, which vanishes as 
$R\to\infty$. Therefore, 
$\lg \psi_{\alpha}, \nu_n\rg \to \lg \psi_{\alpha}, \nu\rg$ 
as $n\to\infty$. \qed  

\medskip 

\noindent 
Since we know from Lemma 5.1 (i) that 
$\cM_{a,b}^{\gamma}$ are closed subsets of 
$(\cM_f,d_w)$, the next lemma is regarded as 
a slight generalization of Lemma 4.2 in \cite{EW00}, 
which corresponding to the case $\gamma=1$. 
Their proof will be arranged in an obvious manner.   
\begin{lm}
For any $a,b,\gamma>0$,  
$\cM_{a,b}^{\gamma}$ is compact. 
\end{lm}
{\it Proof.}~
It is sufficient to prove that 
any sequence $\{\nu_n\}\subset \cM_{a,b}^{\gamma}$ 
has a weakly convergent subsequence. 
Since $\sup_n\lg 1,\nu_n\rg \le a$, 
we can choose a subsequence $\{\nu_{n_k}\}$ 
for which 
$a_0:=\lim_{k\to\infty}\lg 1,\nu_{n_k}\rg$ exists. 
In case $a_0=0$, $\nu_{n_k}\stackrel{w}{\to}0$. 
In case $a_0>0$, we may assume further that 
$\lg 1,\nu_{n_k}\rg>a_0/2$ for all $k$ and consider 
probability measures 
$\tilde{\nu}_k:=\lg 1,\nu_{n_k}\rg^{-1}\nu_{n_k}$. 
The family $\{\tilde{\nu}_k\}$ is tight because 
\[ 
\tilde{\nu}_k([R,\infty)) 
= \frac{\lg \one_{[R,\infty)},\nu_{n_k}\rg}
{\lg 1,\nu_{n_k}\rg} 
\le 
R^{-\gamma} 
\frac{\lg \psi_{\gamma}\one_{[R,\infty)},\nu_{n_k}\rg}
{\lg 1,\nu_{n_k}\rg} 
< R^{-\gamma} \frac{2b}{a_0}. 
\] 
Taking its subsequence $\{\tilde{\nu}_{k_l}\}$ 
such that $\tilde{\nu}_{k_l}\stackrel{w}{\to}
\tilde{\nu}$ as $l\to \infty$ for some 
$\tilde{\nu}\in\cM_f$, we see that 
$\nu_{n_{k_l}}\stackrel{w}{\to}a_0\tilde{\nu}$ 
as desired. \qed 

\medskip 

We now prove the compact containment property of 
the laws of $\{\mu^N(t):~ t\ge 0\}$. 
Note the triviality 
$\lg \psi_a,\mu^N(t)\rg=\lg \psi_{a+1},\xi^N(t)\rg$. 
More generally, for any $f\in B_c$ 
$\lg f,\mu^N(t)\rg=\lg f^{\star},\xi^N(t)\rg$, 
where $f^{\star}\in B_c$ is defined to be 
$f^{\star}(y)=yf(y)$. 
\begin{pr}
Assume that (\ref{4.2}) holds for some $\delta>1$. 
Suppose that $\sum\delta_{Z^N_i(0)}$ is 
${\rm Po}(Nc_0)$-distributed for each $N=1,2,\ldots$. 
Then for any $T>0$ and $\epsilon\in (0,1)$ 
there exist $a,b>0$ such that 
\be 
\inf_N P\left(\mu^N(t)\in \cM_{a,b}^{\delta-1} 
\quad 
\mbox{for all} \quad t\in [0,T] \right) 
\ge 1-\epsilon.                            \label{5.1}
\ee  
\end{pr}
{\it Proof.}~ Lemma 4.1 with $a=2+\lambda+\delta$ 
ensures the validity of (\ref{4.15}), which enables us to apply 
Proposition 4.5 and Lemma 4.7 with $\gamma=\delta$. 
Since 
\[ 
\lg 1,\mu^N(t)\rg =\lg \psi_1,\xi^N(t)\rg
=N^{-1}|Z^N(t)|=N^{-1}|Z^N(0)|=\lg \psi_1,\xi^N(0)\rg
\] 
and 
\begin{eqnarray*} 
\lg \psi_{\delta-1},\mu^N(t)\rg
& = & \lg \psi_{\delta},\xi^N(t)\rg              \\ 
& = & 
\lg \psi_{\delta},\xi^N(0)\rg  
+\frac{1}{N}\int_0^tL^N\Psi_{\delta}(Z^N(s))ds 
+\wt{M}_{\delta}^N(t), 
\end{eqnarray*} 
Chebyshev's inequality and Doob's  inequality for 
submartingale (e.g. Corollary 2.17, Chapter 2 in \cite{EK}) 
together yield 
\begin{eqnarray*} 
\lefteqn{1-P\left(\mu^N(t)\in \cM_{a,b}^{\delta-1} 
\quad 
\mbox{for all} \quad t\in [0,T] \right)}         \\ 
& \le & 
P\left(\sup_{t\in [0,T]}\lg 1,\mu^N(t)\rg >a\right)
+P\left(
\sup_{t\in [0,T]}\lg \psi_{\delta-1},\mu^N(t)\rg >b \right) \\ 
& \le & 
P\left(\lg \psi_1,\xi^N(0)\rg >a\right)
+P\left(\lg \psi_{\delta},\xi^N(0)\rg >\frac{b}{3}\right)  \\
& & 
+P\left(\frac{1}{N}\int_0^T\left|L^N
\Psi_{\delta}(Z^N(s))\right|ds >\frac{b}{3} \right) 
+P\left(\sup_{t\in [0,T]}|\wt{M}_{\delta}^N(t)|>\frac{b}{3}\right)\\
& = & 
a^{-1}E\left[\lg \psi_1,\xi^N(0)\rg \right]
+3b^{-1}E\left[\lg \psi_{\delta},\xi^N(0)\rg\right]   \\ 
& & 
+3b^{-1}
E\left[\frac{1}{N}\int_0^T\left|L^N
\Psi_{\delta}(Z^N(s))\right|ds\right] 
+3b^{-1}E\left[|\wt{M}_{\delta}^N(T)|\right]. 
\end{eqnarray*} 
Therefore, the proof of (\ref{5.1}) reduces to showing 
that the four expectations in the above 
are bounded in $N$. The first two ones are finite 
and independent of $N$: 
\[ 
E\left[\lg \psi_1,\xi^N(0)\rg \right]
=\lg \psi_1,c_0\rg 
\quad \mbox{and} \quad 
E\left[\lg \psi_{\delta},\xi^N(0)\rg \right]
=\lg \psi_{\delta},c_0\rg. 
\] 
For the third expectation we deduce  
from (\ref{4.25}) and (\ref{4.20}) that 
\begin{eqnarray*} 
\lefteqn{ E_N \ := \   
E\left[\frac{1}{N}\int_0^T\left|L^N\Psi_{\delta}(Z^N(s))
\right|ds\right]}         \\ 
& \le & 
\Ch C_{1,\lambda+\delta}
E\left[\int_0^T \lg \psi_1,\xi^N(0) \rg    
\lg \psi_{1+\lambda+\delta},\xi^N(s) \rg ds \right]   
+\frac{\check{C}_{\delta}}{2}                                
E\left[\int_0^T \lg \psi_{2+\lambda+\delta},\xi^N(s) \rg ds 
\right] \\ 
& \le & 
\Ch C_{1,\lambda+\delta} 
\left(T \cdot 
E\left[\lg \psi_1, \xi^N(0)\rg^{2+\lambda+\delta}\right]
\right)^{\frac{1}{1+\lambda+\delta}}
\left(E\left[\int_0^T  
\lg \psi_{2+\lambda+\delta},\xi^N(s) \rg ds \right] 
\right)^{\frac{\lambda+\delta}{1+\lambda+\delta}}       \\ 
& & 
+\frac{\check{C}_{\delta}}{2}                                
E\left[\int_0^T \lg \psi_{2+\lambda+\delta},\xi^N(s) \rg ds 
\right]. 
\end{eqnarray*} 
Thanking to 
$\overline{m}_{2+\lambda+\delta}<\infty$ and (\ref{4.16}), 
this is bounded in $N$ . 
Lastly, in view of (\ref{4.24}) 
\begin{eqnarray*} 
E\left[|\wt{M}_{\delta}^N(T)|\right] 
& \le & 
E\left[\langle \psi_{\delta},\xi^N(T) \rangle \right]   
+E\left[\langle \psi_{\delta},\xi^N(0) \rangle \right] 
+E_N                                            \\ 
& \le & 
C^*(T) + \lg \psi_{\delta},c_0\rg 
+ E_N, 
\end{eqnarray*} 
where the last inequality follows from (\ref{4.18}). 
We complete the proof of Proposition 5.3 since 
we have already seen that $\sup_N E_N<\infty$.    \qed 

\medskip 

Let $D([0,\infty),\cM_f)$ denote the space of 
right continuous functions $\nu(\cdot)$ 
from $[0,\infty)$ into $\cM_f$ with left limits. 
This space is equipped with a metric 
which induces the Skorohod topology.  
We now state the main result of this subsection, 
which establishes tightness of the rescaled processes. 
Roughly speaking, 
the reason why a stronger assumption on $c_0$ 
than the one in Proposition 5.3 is made here 
is that showing equicontinuity of 
$t\mapsto \lg f,\xi^N(t)\rg$ requires to dominate 
$N^{-1}E[|L^N\Phi_f(Z^N(t))|]$ uniformly in $t>0$. 
\begin{th}
Assume (\ref{4.2}) with $\delta=2+\lambda$, 
i.e., 
\[ 
0< \langle \psi_1, c_0 \rangle <\infty \quad 
\mbox{and} \quad 
\langle \psi_{4+2\lambda}, c_0 \rangle <\infty.   
\] 
Suppose that $\sum \delta_{Z_i^N(0)}$ is 
${\rm Po}(Nc_0)$-distributed for each $N=1,2,\ldots$. 
Then the sequence $\{\cP^N\}_{N=1}^{\infty}$ of 
the laws $\cP^N$ of $\{\mu^N(t):~t\ge 0\}$ 
on $D([0,\infty),\cM_f)$ is relatively compact. 
\end{th}
{\it Proof.}~ Lemma 5.2 and Proposition 5.3 
together imply that for arbitrarily fixed $t\ge 0$ 
the family of the laws of 
$\mu^N(t)$ $(N=1,2,\ldots)$ on $\cM_f$ 
is relatively compact. 
Therefore, for the same reasoning as in 
the proof of Theorem 2.1 in \cite{EW}, 
which exploits Corollary 7.4 in Chapter 3 of \cite{EK}, 
it is sufficient to prove that, for any $T>0$ 
and $\epsilon>0$, there exists $0<\Delta<1$ such that 
\be 
\limsup_{N\to\infty}
P\left(\max_{i:~t_i<T} \sup_{s\in[t_i,t_{i+1})}
d_w(\mu^N(s),\mu^N(t_i))>\frac{\epsilon}{2}\right)
\le \epsilon,                                        \label{5.2}
\ee 
where $t_i=i\Delta \ (i=0,1,2,\ldots)$.  
By (\ref{5.1}) with $\delta=2+\lambda$ 
we can find $a,b>0$ such that   
\[   
\inf_NP\left(W^N_{a,b}(T)\right) 
\ge 1-\frac{\epsilon}{2},    
\]   
where $W^N_{a,b}(T)$ is the event that 
$\mu^N(t)\in \cM_{a,b}^{1+\lambda}$  
for all $t\in [0,T+1]$. So, we only have to prove 
\be 
\limsup_{N\to\infty}
P\left(\left\{\max_{i:~t_i<T} \sup_{s\in[t_i,t_{i+1})}
d_w(\mu^N(s),\mu^N(t_i))>\frac{\epsilon}{2}\right\}
\cap W^N_{a,b}(T)\right)
\le \frac{\epsilon}{2}.                                \label{5.3}
\ee 
It follows from (\ref{4.23}) that, on $W^N_{a,b}(T)$,  
for any $t\in [0,T]$ and $s\in[t,t+\Delta)$ 
\begin{eqnarray*} 
\lefteqn{|\lg h_k, \mu^N(s)\rg-\lg h_k, \mu^N(t)\rg|
\ = \ 
|\lg h_k^{\star},\xi^N(s)\rg-\lg h_k^{\star}, \xi^N(t)\rg|}  \\ 
& \le & 
\left|M^N_{h_k^{\star}}(s)-M^N_{h_k^{\star}}(t)\right| 
+\frac{1}{N}\int_t^s 
\left|L^N\Phi_{h_k^{\star}}(Z^N(u))\right|du  \\ 
& \le & 
2\sup_{u\in[0,T+1]}\left|M^N_{h_k^{\star}}(u)\right|
+ 
\Delta \left(\Ch C_{1,\lambda} 
a^{\frac{2+\lambda}{1+\lambda}}    
b^{\frac{\lambda}{1+\lambda}}    
+\frac{\Cc}{2} b\right)\Vert h_k^{\star} \Vert_{\infty}.   
\end{eqnarray*} 
Hence, taking $k_0\in\N$ such that 
$\sum_{k=k_0+1}^{\infty}2^{-k}<\epsilon/4$ and then 
$\Delta$ sufficiently small so that  
\[ 
\Delta \left(\Ch C_{1,\lambda} 
a^{\frac{2+\lambda}{1+\lambda}}    
b^{\frac{\lambda}{1+\lambda}}    
+\frac{\Cc}{2} b\right)\Vert h_k^{\star} \Vert_{\infty}
<\frac{\epsilon}{8k_0} 
\quad \mbox{for all $k\in\{1,\ldots,k_0\}$},  
\]  
we get 
\begin{eqnarray*} 
\lefteqn{P\left(\left\{\max_{i:~t_i<T} \sup_{s\in[t_i,t_{i+1})}
d_w(\mu^N(s),\mu^N(t_i))>\frac{\epsilon}{2}\right\}
\cap W^N_{a,b}(T)\right)}                         \\ 
& \le & 
P\left(\sum_{k=1}^{k_0}2^{-k} 
\left\{2\sup_{u\in[0,T+1]}
\left|M^N_{h_k^{\star}}(u)\right|
+\frac{\epsilon}{8k_0}\right\}
+\sum_{k=k_0+1}^{\infty}2^{-k}
>\frac{\epsilon}{2}\right)             \\ 
& \le & 
\sum_{k=1}^{k_0}P\left(
2\sup_{u\in[0,T+1]}\left|M^N_{h_k^{\star}}(u)\right|
+\frac{\epsilon}{8k_0}>\frac{\epsilon}{4k_0}\right)    \\ 
& \le & 
\frac{16k_0}{\epsilon}\sum_{k=1}^{k_0}E\left[
\left|M^N_{h_k^{\star}}(T+1)\right|\right], 
\end{eqnarray*} 
in which the last inequality is implied by 
Doob's inequality. Since each expectation 
in the above sum converges to 0 as $N\to\infty$ 
by Corollary 4.6, (\ref{5.3}) and hence (\ref{5.2}) 
are obtained. 
The proof of Theorem 5.4 is complete.  \qed

\subsection{Studying the limit laws}
This section is devoted to the study of 
the weak limit of an arbitrary convergent subsequence  $\{\cP^{N_l}\}_{l=1}^{\infty}$, say, of 
the laws $\cP^N$ of $\{\mu^N(t):~t\ge 0\}$ 
on $D([0,\infty),\cM_f)$. 
Although our main concern will be 
proving that under the limit law the weak form of 
(\ref{4.5}) is satisfied almost surely, 
some properties on the limit are shown in advance. 
Let $C([0,\infty),\cM_f)$ be the space of 
continuous functions from $[0,\infty)$ to $\cM_f$. 
According to \S 10 of Chapter 3 of \cite{EK} 
it is equipped with the metric 
\[ 
d_U(\nu_1(\cdot),\nu_2(\cdot))
=\int_0^{\infty}e^{-u}\sup_{t\in[0,u]}
\{d_w(\nu_1(t),\nu_2(t)) \wedge 1\}du, 
\] 
which gives the topology of uniform convergence
 on compact subsets of $[0,\infty)$. 
\begin{lm}
If $\{\mu^{N_l}(t):~t\ge 0\}$ converges to 
a process $\{\mu(t):~t\ge 0\}$ 
in distribution on $D([0,\infty),\cM_f)$ as $l\to\infty$, 
then $\mu(\cdot)\in C([0,\infty),\cM_f)$ a.s. 
\end{lm}
{\it Proof.}~As in the proof of Lemma 5.1 in \cite{EW}, 
we employ Theorem 10.2 in Chapter 3 of \cite{EK}. 
To this end, observe that if 
$\mu^N(t)\ne\mu^N(t-)=N^{-1}\sum z_i \delta_{z_i}$ 
for some $\bz\in\Omega$, 
the signed measure $\mu^N(t)-\mu^N(t-)$ equals either 
\[ 
\frac{1}{N}(z_i+z_j)\delta_{z_i+z_j}
-\frac{1}{N}\left(z_i\delta_{z_i}+z_j\delta_{z_j}\right) 
\quad \mbox{for some $z_i,z_j>0$ with $i\ne j$} 
\] 
or 
\[ 
\frac{1}{N}\left(y\delta_y+(z_i-y)\delta_{z_i-y}\right) 
-\frac{1}{N}z_i\delta_{z_i}  \quad 
\mbox{for some $z_i>0$ and $y\in(0,z_i)$}. 
\] 
This implies that 
\begin{eqnarray*} 
\lefteqn{d_w(\mu^N(t),\mu^N(t-))}                   \\ 
& \le & 
\sup_{i\ne j}\sum_{k=0}^{\infty}2^{-k}
\min\left\{1,
N^{-1}\left|h_k^{\star}(z_i+z_j)
-h_k^{\star}(z_i)-h_k^{\star}(z_j)\right|
\right\} \\ 
& & 
\vee \ \sup_{i}\sup_{y\in(0,z_i)}\sum_{k=0}^{\infty}2^{-k}
\min\left\{1,
N^{-1}\left|h_k^{\star}(y)+h_k^{\star}(z_i-y)
-h_k^{\star}(z_i)\right|
\right\}   \\ 
& \le  & 
\sum_{k=1}^{\infty}2^{-k}\min\{1,
 3N^{-1}\Vert h_k^{\star}\Vert_{\infty}\}
\to 0 \quad \mbox{as} \quad N\to \infty,  
\end{eqnarray*} 
in which all the terms for $k=0$ in the sums vanish  
because of $h_0^{\star}(y)=y$. So, the above mentioned 
theorem proves the assertion. \qed 

\medskip 

\noindent 
Given $\nu\in\cM$, we define two measures 
$\nu^{\star}$ and $\nu_{\star}$ on $(0,\infty)$ 
by $\nu^{\star}(dy)=y\nu(dy)$ 
and $\nu_{\star}(dy)=y^{-1}\nu(dy)$, respectively. 
For instance, $\mu^N(t)=\xi^N(t)^{\star}$ 
and conversely $\xi^N(t)=\mu^N(t)_{\star}$
\begin{lm}
Suppose that the same assumptions 
as in  Proposition 5.3 hold. 
If $\{\mu^{N_l}(t):~t\ge 0\}$ converges to 
a process $\{\mu(t):~t\ge 0\}$ 
in distribution on $D([0,\infty),\cM_f)$ as $l\to\infty$, 
then 
\be 
P\left(\mu(0)=c_0^{\star}, \ \lg 1, \mu(t)\rg 
=\lg 1, \mu(0)\rg \ \mbox{for all} 
\ t\ge 0 \right)=1                                   \label{5.4}
\ee  
and for each $T>0$ there exist a constant $\wt{C}(T)$
such that 
\be 
E\left[\lg \psi_{\delta-1},\mu(T) \rg
+\int_0^T \left\{\lg 1,\mu(0) \rg    
\lg \psi_{\lambda+\delta},\mu(s) \rg 
+\lg \psi_{1+\lambda+\delta},\mu(s) \rg\right\} ds 
\right]    \le \wt{C}(T).                              \label{5.5} 
\ee 
\end{lm}
{\it Proof.}~ 
By the assumption on initial distributions  
and Lemma 4.1 together 
\[ 
P(\mu(0)_{\star}=c_0)=1 
\quad \mbox{or equivalently}  \quad 
P(\mu(0)=c_0^{\star})=1. 
\]  
For each $N$, 
$P(\lg 1,\mu(t)^N\rg =\lg 1,\mu^N(0)\rg$ for all $t\ge 0)=1$  
and as is seen easily 
\[ 
\{\nu(\cdot)\in D([0,\infty),\cM_f):~ 
\ \lg 1, \nu(t)\rg =\lg 1, \nu(0)\rg 
\ \mbox{for all} \ t\ge 0\} 
\] 
is a closed subset of $D([0,\infty),\cM_f)$. 
So, by the assumed convergence in distribution 
$P\left(\lg 1, \mu(t)\rg 
=\lg 1, \mu(0)\rg \ \mbox{for all} \ t\ge 0 \right)=1$, 
and summarizing, we have shown (\ref{5.4}). 

We proceed to verification of (\ref{5.5}). 
In estimating $E_N$ in the proof of Proposition 5.3 
we have shown the existence of  
a constant $\wt{C}_1(T)$ independent of $N$ such that 
\[ 
E\left[\int_0^T\left\{\lg 1,\mu^N(0) \rg    
\lg \psi_{\lambda+\delta},\mu^N(s) \rg 
+\lg \psi_{1+\lambda+\delta},\mu^N(s) \rg\right\} ds 
\right] \le \wt{C}_1(T).                                
\] 
Combining this with (\ref{4.18}), we obtain 
\begin{eqnarray} 
\lefteqn{ 
E\left[\lg \psi_{\delta-1},\mu^N(T) \rg\right]} 
                                            \label{5.6}     \\ 
& + & 
E\left[\int_0^T\left\{\lg 1,\mu^N(0) \rg    
\lg \psi_{\lambda+\delta},\mu^N(s) \rg 
+\lg \psi_{1+\lambda+\delta},\mu^N(s) \rg\right\} ds 
\right] \le \wt{C}(T),                  \nonumber  
\end{eqnarray} 
where $\wt{C}(T) =\wt{C}_1(T)+C^{*}(T)$.  
Therefore, for any $a,R>0$ and $l\in\N$ 
\begin{eqnarray} 
\lefteqn{
E\left[a \wedge \lg \psi_{\delta-1}
\varphi_{R},\mu^{N_l}(T) \rg \right]} 
                                            \label{5.7}     \\ 
& + & 
E\left[\int_0^T\min\left\{a, \lg 1,\mu^{N_l}(0) \rg    
\lg \psi_{\lambda+\delta}\varphi_{R},\mu^{N_l}(s) \rg 
+\lg \psi_{1+\lambda+\delta}\varphi_{R},\mu^{N_l}(s) \rg
\right\} ds \right] \le \wt{C}(T).               \nonumber  
\end{eqnarray} 
It is not difficult to check that 
the function on $D([0,\infty),\cM_f)$ 
\[ 
\nu(\cdot)\mapsto 
\int_0^T\min\left\{a, \lg 1,\nu(0) \rg    
\lg \psi_{\lambda+\delta}\varphi_{R},\nu(s) \rg 
+\lg \psi_{1+\lambda+\delta}\varphi_{R},\nu(s) \rg
\right\} ds 
\] 
is bounded and continuous. 
Hence, letting $l\to\infty$ in (\ref{5.7}) yields 
(\ref{5.7}) with $\mu$ in placed of $\mu^{N_l}$. 
Since $a$ and $R$ are arbitrary, (\ref{5.5}) holds true 
by virtue of the monotone convergence theorem. \qed 

\medskip 

We are in a position to state the main result of 
this section. It should be emphasized that 
the continuity of $K$ and $F$ is required 
only in the next theorem. 
\begin{th}
In addition to the assumptions of Theorem 5.4, 
assume that $\Hh$ and $\Hc$ are continuous. 
If $\{\mu^{N_l}(t):~t\ge 0\}$ converges 
to a process $\{\mu(t):~t\ge 0\}$ 
in distribution on $D([0,\infty),\cM_f)$ as $l\to\infty$, 
then $\{\xi^{N_l}(t):~t\ge 0\}$ converges in distribution 
on $D([0,\infty),\cM)$ to $\{\xi(t):~t\ge 0\}$ defined by 
$\xi(t)=\mu(t)_{\star}$. Moreover, 
with probability 1, it holds that 
for any $f\in B_c$ and $t\ge 0$ 
\begin{eqnarray} 
\lg f,\xi(t)\rg  -\lg f,c_0 \rg     
& = & 
\frac{1}{2}\int_0^tds 
\int \xi(s)^{\otimes 2}(dxdy)xy \Hh(x,y) (\Box f)(x,y)  \label{5.8} \\ 
& & 
-\frac{1}{2}\int_0^tds 
\int \xi(s)(dx) x\int_0^{x}dy\Hc(y,x-y) (\Box f)(y,x-y)  \nonumber
\end{eqnarray} 
with the integrals on the right side being 
absolutely convergent. 
\end{th}
From the analytic view point, 
this theorem particularly implies 
the existence of a $\cM$-valued weak solution 
to (\ref{4.5}) with symmetric continuous homogeneous 
functions $\Hh$ and $\Hc$ of degree $\lambda\ge 0$  satisfying (H1) and (H2), respectively, 
and with initial measure $c_0$ such that 
$0< \langle \psi_1, c_0 \rangle <\infty$ and 
$\langle \psi_{4+2\lambda}, c_0 \rangle <\infty$. 
Unfortunately, the uniqueness of the solution 
has not been proved and accordingly 
the convergence of the laws of $\{\mu^N(t):~t\ge 0\}$ 
($N=1,2,\ldots)$ has not been obtained. 
Concerning this point, some comments will be given 
at the end of this section. 

\medskip 

\noindent 
{\it Proof of Theorem 5.7.}~ 
In order to prove the first half 
it is sufficient to verify continuity of the map 
$\{\nu(t):~ t\ge 0\} \mapsto 
\{\nu(t)_{\star}:~ t\ge 0\}$ 
from $D([0,\infty),\cM_f)$ to $D([0,\infty),\cM)$. 
But, that continuity follows from an 
equivalent condition to the convergence 
with respect to the Skorohod topology 
(e.g. Proposition 5.3 of Chapter 3 in \cite{EK}) 
in terms of the metric on the state space 
$\cM_f$ or $\cM$ together with continuity of the map 
$\nu\mapsto \nu_{\star}$ from $\cM_f$ to $\cM$.  
(cf. Problem 13 of Chapter 3 in \cite{EK}.) \par 
The proof of the last half is 
divided into five steps.  \par 
{\it Step 1.}~ 
As was sketched roughly in \S 4.1, 
the argument will be based on (\ref{4.6}).  
So, for any $f\in C_c$ and $t>0$, set 
for $\nu(\cdot)\in D([0,\infty),\cM_f)$ 
\begin{eqnarray*} 
I_{t,f}(\nu(\cdot)) 
& = &
\frac{1}{2}\int_0^tds 
\int \nu(s)^{\otimes 2}(dxdy)\Hh(x,y) (\Box f)(x,y)  \\ 
& & 
-\frac{1}{2}\int_0^tds 
\int \nu(s)(dx) \int_0^{x}dy\Hc(y,x-y) (\Box f)(y,x-y)  
\end{eqnarray*} 
provided that both integrals converge absolutely.  
Then, by the assumption and (\ref{4.23}) 
we have almost surely 
\[ 
\lg f,\xi^N(t)\rg -\lg f,\xi^N(0)\rg-I_{t,f}(\mu^N(\cdot)) 
=
M_f^N(t)+D_f^N(t) \quad \mbox{for all $t\ge 0$},  
\] 
where 
\begin{eqnarray*} 
D_f^N(t) 
& = & 
\frac{1}{2}\int_0^t ds\int \left(\xi^N(s)^{[2]}
-\xi^N(s)^{\otimes 2}\right)(dxdy)xy \Hh(x,y) 
(\Box f)(x,y)                                            \\ 
& = & 
\frac{1}{2N^2}\int_0^t ds\sum_{i} Z_i^N(s)^2 
\Hh(Z_i^N(s),Z_i^N(s))\left\{f(2Z_i^N(s))-2f(Z_i^N(s))\right\}.  
\end{eqnarray*} 
By the same calculations as in \S 4.1 
(see the observation after (\ref{4.6}))  
it is readily shown that for some constant $C_f'$ 
independent of $N\in \N$ and $T>0$ 
\[ 
E\left[\sup_{t\in[0,T]} |D_f^N(t)|\right] 
\le \frac{C_f'T}{N}. 
\] 
By combining this with (\ref{4.21}) 
and using Doob's inequality we get for any $\epsilon>0$ 
\begin{eqnarray} 
\lefteqn{P\left(\sup_{t\in[0,T]} 
\left|\lg f,\xi^N(t)\rg -\lg f,\xi^N(0)\rg 
-I_{t,f}(\mu^N(\cdot))\right| > \epsilon \right)} 
                                                   \nonumber   \\ 
& \le & 
P\left(\sup_{t\in[0,T]} 
\left|M_f^N(t)\right| > \frac{\epsilon}{2} \right) 
+ P\left(\sup_{t\in[0,T]} 
\left|D_f^N(t)\right| > \frac{\epsilon}{2} \right)  
                                                   \nonumber   \\ 
& \le & 
\left(\frac{2}{\epsilon}\right)^2
\frac{4\Vert f\Vert_{\infty}^2C(T)}{N}
+\frac{2}{\epsilon}\cdot\frac{C_f'T}{N}.     \label{5.9} 
\end{eqnarray} 
Therefore, the main task in the rest of the proof is 
to show, in a suitable sense, convergence of 
$I_{t,f}(\mu^{N_l}(\cdot))$ to $I_{t,f}(\mu(\cdot))$ 
as $l \to \infty$. But $I_{t,f}$ cannot be defined 
as a function on $D([0,\infty),\cM_f)$ 
and we need to handle by a cut-off argument. \par 
{\it Step 2.}~ 
For each $R>0$ decompose $I_{t,f}$ in the form   
\be  
I_{t,f}(\nu(\cdot))
=
I_{t,f,R}^{(2)}(\nu(\cdot))
+\widetilde{I}_{t,f,R}^{(2)}(\nu(\cdot))
-I_{t,f,R}^{(1)}(\nu(\cdot))
-\widetilde{I}_{t,f,R}^{(1)}(\nu(\cdot)), 
                                                       \label{5.10} 
\ee 
where 
\[ 
I_{t,f,R}^{(2)}(\nu(\cdot))
= \frac{1}{2}\int_0^tds 
\int \nu(s)^{\otimes 2}(dxdy)
\Hh(x,y) (\Box f)(x,y)\varphi_R(x+y), 
\] 
\[ 
\widetilde{I}_{t,f,R}^{(2)}(\nu(\cdot))
= \frac{1}{2}\int_0^tds 
\int \nu(s)^{\otimes 2}(dxdy)
\Hh(x,y) (\Box f)(x,y)(1-\varphi_R(x+y)), 
\] 
\[ 
I_{t,f,R}^{(1)}(\nu(\cdot))
= \frac{1}{2}\int_0^tds \int \nu(s)(dx) 
\int_0^{x}dy\Hc(y,x-y) (\Box f)(y,x-y)\varphi_R(x) 
\] 
and 
\[ 
\widetilde{I}_{t,f,R}^{(1)}(\nu(\cdot))  
= 
\frac{1}{2}\int_0^tds \int \nu(s)(dx) 
\int_0^{x}dy\Hc(y,x-y) (\Box f)(y,x-y)(1-\varphi_R(x)). 
\] 
Of course $I_{t,f,R}^{(2)}$ and $I_{t,f,R}^{(1)}$ 
should be dominant for $R$ large. 
Putting $\delta=\lambda+2$, we actually claim that 
\be 
\sup_{t\in[0,T]}
\left|\widetilde{I}_{t,f,R}^{(2)}(\nu(\cdot))\right|  
\le 
\frac{\Ch C_{1,\lambda+\delta}}{R^{\delta}} 
\Vert \Box f \Vert_{\infty} \int_0^T  \lg 1,\nu(s) \rg \lg 
\psi_{\lambda+\delta},\nu(s) \rg ds      \label{5.11} 
\ee 
and 
\be \sup_{t\in[0,T]}
\left|\widetilde{I}_{t,f,R}^{(1)}(\nu(\cdot))\right|  
\le 
\frac{\Cc}{2R^{\delta}}\Vert \Box f \Vert_{\infty} \int_0^T 
\lg \psi_{1+\lambda+\delta},\nu(s) \rg ds   \label{5.12} 
\ee 
whenever each integral on the right side is finite.  
Indeed, by (\ref{3.4}) and (\ref{4.7}) 
\begin{eqnarray*} 
\left|\widetilde{I}_{t,f,R}^{(2)}(\nu(\cdot))\right|     
& \le & 
\frac{\Ch}{2}
\int_0^tds \int \nu(s)^{\otimes 2}(dxdy)
(x+y)^{\lambda} |(\Box f)(x,y)|\one_{\{x+y>R\}} \\ 
& \le & 
\frac{\Ch}{2}\Vert \Box f \Vert_{\infty}
\int_0^tds \int \nu(s)^{\otimes 2}(dxdy)
(x+y)^{\lambda+\delta}
\frac{\one_{\{x+y>R\}}}{R^{\delta}}               \\ 
& \le  & 
\frac{\Ch C_{1,\lambda+\delta}}{2R^{\delta}} 
\Vert \Box f \Vert_{\infty}
\int_0^tds \int \nu(s)^{\otimes 2}(dxdy)
(x^{\lambda+\delta}+y^{\lambda+\delta}), 
\end{eqnarray*} 
from which (\ref{5.11}) is immediate. 
Similarly, by (\ref{3.6})  
\begin{eqnarray*} 
\left|\widetilde{I}_{t,f,R}^{(1)}(\nu(\cdot))\right|     
& \le & 
\frac{\Cc}{2}\Vert \Box f \Vert_{\infty} 
\int_0^tds \int \nu(s)(dx) 
x^{1+\lambda}\one_{\{x>R\}} \\ 
& \le & 
\frac{\Ch}{2}\Vert \Box f \Vert_{\infty}
\int_0^tds \int \nu(s)(dx)x^{1+\lambda+\delta}
\frac{\one_{\{x>R\}}}{R^{\delta}}              \\ 
& \le  & 
\frac{\Cc}{2R^{\delta}}\Vert \Box f \Vert_{\infty} 
\int_0^tds \int \nu(s)(dx)x^{1+\lambda+\delta} 
\end{eqnarray*} 
and thus (\ref{5.12}) is valid. 
Now consider (\ref{5.11}) and (\ref{5.12}) 
with $\nu(\cdot)$ being a random element 
$\mu^N(\cdot)$ or $\mu(\cdot)$. 
Thanking to Chebyshev's inequality, 
taking expectations and then 
using (\ref{5.6}), (\ref{5.4}) and (\ref{5.5}) lead to 
\be 
P\left(\sup_{t\in[0,T]}\left\{
\left|\widetilde{I}_{t,f,R}^{(2)}(\mu^N(\cdot))\right|
+\left|\widetilde{I}_{t,f,R}^{(1)}(\mu^N(\cdot))\right|\right\}
> \epsilon\right) 
\le \frac{2\Ch C_{1,\lambda+\delta}+\Cc}{2\epsilon R^{\delta}}
\Vert \Box f \Vert_{\infty} \widetilde{C}(T)      \label{5.13} 
\ee 
and 
\be 
P\left(\sup_{t\in[0,T]}\left\{
\left|\widetilde{I}_{t,f,R}^{(2)}(\mu(\cdot))\right|
+\left|\widetilde{I}_{t,f,R}^{(1)}(\mu(\cdot))\right|\right\}
> \epsilon\right) 
\le \frac{2\Ch C_{1,\lambda+\delta}+\Cc}{2\epsilon R^{\delta}}
\Vert \Box f \Vert_{\infty} \widetilde{C}(T).     \label{5.14} 
\ee 

{\it Step 3.}~ Let $f\in C_c$ be arbitrary. 
Clearly proving that 
\be 
P\left(\lg f,\xi(t) \rg -\lg f,\xi(0) \rg 
-I_{t,f}(\mu (\cdot))
=0 \quad \mbox{for all $t\ge 0$}\right)=1         \label{5.15} 
\ee 
is equivalent to showing that for any 
$T>0$ and $\epsilon\in (0,1)$   
\be 
P\left(\sup_{t\in[0,T]}\left|\lg f,\xi(t) \rg -\lg f,\xi(0) \rg
-I_{t,f}(\mu (\cdot))\right|>2\epsilon\right)=0.     \label{5.16} 
\ee 
We claim here that the latter can be reduced to 
establishing the inequality 
\begin{eqnarray} 
\lefteqn{
P\left(\sup_{t\in[0,T]}\left|\lg f,\xi(t) \rg -\lg f,\xi(0) \rg
-I_{t,f,R}^{(2)}(\mu (\cdot))+I_{t,f,R}^{(1)}(\mu (\cdot))
\right|>\epsilon\right)}                            \label{5.17} \\  
& \le & 
\liminf_{l \to \infty}
P\left(\sup_{t\in[0,T]}\left|\lg f,\xi^{N_l}(t) \rg 
-\lg f,\xi^{N_l} (0) \rg
-I_{t,f,R}^{(2)}(\mu^{N_l} (\cdot))
+I_{t,f,R}^{(1)}(\mu^{N_l} (\cdot))
\right|>\epsilon\right)    \nonumber 
\end{eqnarray} 
for each $R>0$. Indeed, by (\ref{5.10}) and (\ref{5.14}) 
\begin{eqnarray} 
\lefteqn{ 
P\left(\sup_{t\in[0,T]}\left|\lg f,\xi(t) \rg -\lg f,\xi(0) \rg
-I_{t,f}(\mu (\cdot))
\right|>2\epsilon\right)}                            \label{5.18} \\  
& \le & 
P\left(\sup_{t\in[0,T]}\left|\lg f,\xi(t) \rg -\lg f,\xi(0) \rg
-I_{t,f,R}^{(2)}(\mu (\cdot))+I_{t,f,R}^{(1)}(\mu (\cdot))
\right|>\epsilon\right) \nonumber  \\ 
& &  
+\frac{2\Ch C_{1,\lambda+\delta}+\Cc}{2\epsilon R^{\delta}}
\Vert \Box f \Vert_{\infty} \widetilde{C}(T). \nonumber 
\end{eqnarray}  
On the other hand, by (\ref{5.9}) and (\ref{5.13}) 
\begin{eqnarray*} 
\lefteqn{ \liminf_{l \to \infty} 
P\left(\sup_{t\in[0,T]}\left|\lg f,\xi^{N_l}(t) \rg 
-\lg f,\xi^{N_l}(0) \rg
-I_{t,f,R}^{(2)}(\mu^{N_l} (\cdot))
+I_{t,f,R}^{(1)}(\mu^{N_l} (\cdot))
\right|>\epsilon\right) }                                    \\  
& \le & 
\liminf_{l \to \infty} 
P\left(\sup_{t\in[0,T]}\left|\lg f,\xi^{N_l}(t) \rg 
-\lg f,\xi^{N_l}(0) \rg-I_{t,f}(\mu^{N_l} (\cdot))
\right|>\frac{\epsilon}{2}\right)
\nonumber  \\ 
& &  
+\frac{2\Ch C_{1,\lambda+\delta}+\Cc}{2\epsilon R^{\delta}}
\Vert \Box f \Vert_{\infty} \widetilde{C}(T)                \\ 
& \le & 
\frac{2\Ch C_{1,\lambda+\delta}+\Cc}{2\epsilon R^{\delta}}
\Vert \Box f \Vert_{\infty} \widetilde{C}(T). 
\end{eqnarray*}  
Therefore, 
this combined with  (\ref{5.17}) and (\ref{5.18}) yields  
\[ 
P\left(\sup_{t\in[0,T]}\left|\lg f,\xi(t) \rg 
-\lg f,\xi(0) \rg -I_{t,f}(\mu (\cdot))
\right|>2\epsilon\right) 
\le  
\frac{2\Ch C_{1,\lambda+\delta}+\Cc}{\epsilon R^{\delta}}
\Vert \Box f \Vert_{\infty} \widetilde{C}(T) 
\] 
and the aforementioned claim 
follows since $R>0$ is arbitrary. \par 
{\it Step 4.}~ 
We shall prove (\ref{5.17}) for arbitrarily fixed 
$T, R>0$, $\epsilon\in (0,1)$ and 
$f\in \{h_1^{\star}, h_2^{\star}, \ldots\}(\subset C_c)$ 
at least. 
But, because of the triviality that 
$x > \epsilon$ if and only if $1 \wedge x > \epsilon $,  
(\ref{5.17}) can be rewritten into 
\be 
P\left(1 \wedge \Upsilon (\mu(\cdot)) > \epsilon\right)
\le \liminf_{l\to\infty} 
P\left(1 \wedge \Upsilon(\mu^{N_l}(\cdot))>\epsilon\right), 
                                                         \label{5.19} 
\ee 
where $\Upsilon=\Upsilon_{T,f,R}$ is a Borel measurable 
function on $D([0,\infty),\cM_f)$ defined by 
\[  
\Upsilon (\nu(\cdot)) 
=\sup_{t\in[0,T]}\left|\lg f_{\star},\nu(t) \rg 
-\lg f_{\star},\nu(0) \rg 
-I_{t,f,R}^{(2)}(\nu (\cdot))
+I_{t,f,R}^{(1)}(\nu (\cdot))\right|.    
\]  
Furthermore, (\ref{5.19}) can be reduced to 
showing that for any $n\in\N$ 
\begin{eqnarray} 
\lefteqn{
P\left(1 \wedge \Upsilon (\mu(\cdot)) > \epsilon \ 
|~\lg 1,\mu(0)\rg \le a_n\right)}         \nonumber    \\ 
& \le & 
\liminf_{l\to\infty} 
P\left(1 \wedge \Upsilon (\mu^{N_l}(\cdot))>\epsilon \ 
|~\lg 1,\mu^{N_l}(0)\rg \le a_n\right),          \label{5.20} 
\end{eqnarray} 
where $a_1,a_2,\ldots$ are such that 
$\lim_{n\to\infty}a_n=\infty$, 
$a_n>\lg \psi_1,c_0\rg(=E[\lg 1,\mu^{N_l}(0)\rg])$ 
and  $\lim_{l\to\infty}P(\lg 1,\mu^{N_l}(0)\rg \le a_n) 
=P(\lg 1,\mu(0)\rg \le a_n)$ for each $n\in\N$. 
Indeed, assuming (\ref{5.20}), we get 
\begin{eqnarray*}
\lefteqn{
P\left(1 \wedge \Upsilon (\mu(\cdot)) > 
\epsilon\right)}                                \\  
& = &  
P\left(1 \wedge \Upsilon (\mu(\cdot)) > \epsilon 
\ |~\lg 1,\mu(0)\rg \le a_n\right)   \             
P\left(\lg 1,\mu(0)\rg \le a_n\right)            \\ 
& \le & 
\liminf_{l\to\infty}
P\left(1 \wedge \Upsilon (\mu^{N_l}(\cdot))>\epsilon \ 
|~\lg 1,\mu^{N_l}(0)\rg \le a_n\right)   \ 
P\left(\lg 1,\mu(0)\rg \le a_n\right)               \\ 
& \le & 
\liminf_{l\to\infty}
\frac{P\left(1 \wedge \Upsilon (\mu^{N_l}(\cdot))>\epsilon\right)}
{P\left(\lg 1,\mu^{N_l}(0)\rg \le a_n\right)}  \cdot 
P\left(\lg 1,\mu(0)\rg \le a_n\right)              \\ 
& \le & 
\frac{
\liminf_{l\to\infty}
P\left(1 \wedge \Upsilon (\mu^{N_l}(\cdot))>\epsilon\right)}
{1-a_n^{-1}\lg \psi_1,c_0\rg}  \cdot 
P\left(\lg 1,\mu(0)\rg \le a_n\right),  
\end{eqnarray*} 
which tends to the right side of (\ref{5.19}) as $n\to\infty$. 
For any $a>0$, 
set $\cM_{\le a}=\{\nu\in\cM_f:~\lg 1,\nu\rg \le a\}$, 
which is regarded as a closed subspace of $\cM_f$. 
Accordingly $D([0,\infty),\cM_{\le a})$ is 
a closed subspace of $D([0,\infty),\cM_{f})$. 
Note that by (\ref{5.4}) 
\[ 
P\left(\mu(\cdot)\in D([0,\infty),\cM_{\le a}) \  
|~\lg 1,\mu(0)\rg \le a\right)= 1 
\] 
and similarly 
\[ 
P\left(\mu^{N_l}(\cdot)\in D([0,\infty),\cM_{\le a}) \  
|~\lg 1,\mu^{N_l}(0)\rg \le a\right)= 1.  
\] 
By the assumption of 
convergence of $\{\mu^{N_l}(t):~t\ge 0\}$ 
to $\{\mu(t):~t\ge 0\}$ in distribution as $l\to\infty$ 
together with 
$\lim_{l\to\infty}P(\lg 1,\mu^{N_l}(0)\rg \le a_n) 
=P(\lg 1,\mu(0)\rg \le a_n)$ $(n=1,2,\ldots)$ 
it is not difficult to verify that for each $n\in\N$ 
the sequence of the conditional laws 
of $\{\mu^{N_l}(t):~t\ge 0\}$ 
on $D([0,\infty),\cM_{\le a_n})$ 
given $\lg 1,\mu^{N_l}(0)\rg \le a_n$
converges to the conditional law of $\{\mu(t):~t\ge 0\}$ 
on $D([0,\infty),\cM_{\le a_n})$ 
given $\lg 1,\mu(0)\rg \le a_n$. Therefore, 
(\ref{5.20}) is naturally expected to follow 
as a consequence of certain continuity of $\Upsilon$ 
restricted on $\cM_{\le a_n}$.  
In fact, by virtue of 
Theorem 10.2 in Chapter 3 of \cite{EK} 
and Lemma 5.5 together, we can conclude (\ref{5.20}) 
as soon as the continuity of $\Upsilon$ 
restricted on $\cM_{\le a_n}$ 
with respect to the metric $d_U$ 
(defined at the beginning of this subsection) 
is checked to hold. We show below more generally 
that continuity of $\Upsilon$ on $\cM_{\le a}$ 
for any $a>0$.  
\par 
For this purpose, take an arbitrary sequence  
$\{\nu_n(\cdot)\}_{n=1}^{\infty}$ of 
$D([0,\infty),\cM_{\le a})$ 
and $\nu(\cdot)\in D([0,\infty),\cM_{\le a})$ such that 
$d_{U}(\nu_n(\cdot), \nu(\cdot))\to 0$ as $n\to\infty$. 
Then our task here is to show that 
$\Upsilon (\nu_n(\cdot)) \to \Upsilon (\nu(\cdot))$ 
for any $f\in \{h_1^{\star}, h_2^{\star}, \ldots\}$. 
From general inequalities of the form 
\[ 
\left|\sup_t |\phi_1(t)| -\sup _t|\phi_2(t)|\right| 
\le \sup_t|\phi_1(t)-\phi_2(t)| 
\le \sup_t|\phi_1(t)|+\sup_t|\phi_2(t)| 
\] 
we deduce 
\begin{eqnarray*} 
\lefteqn{ 
\left|\Upsilon (\nu_n(\cdot))
-\Upsilon (\nu(\cdot))\right|}                        \\ 
& \le & 
2\sup_{t\in[0,T]}
\left|\lg f_{\star},\nu_n(t) \rg 
-\lg f_{\star},\nu(t) \rg  \right|                     \\ 
& & 
+ \sup_{t\in[0,T]}
\left|I_{t,f,R}^{(1)}(\nu_n (\cdot))
-I_{t,f,R}^{(1)}(\nu(\cdot))\right|    
+\sup_{t\in[0,T]}
\left|I_{t,f,R}^{(2)}(\nu_n (\cdot))
-I_{t,f,R}^{(2)}(\nu(\cdot))\right|                \\ 
&=: & 
2s_n + s^{(1)}_n + s^{(2)}_n. 
\end{eqnarray*} 
Letting $i\in \N$ be such that $f=h_i^{\star}$ 
or $f_{\star}=h_i$, we have 
\[ 
s_n= \sup_{t\in[0,T]}\left|\lg h_i,\nu_n(t) \rg 
-\lg h_i,\nu(t) \rg  \right|,
\]   
which converges to 0 as $n\to \infty$  
by $\sup_{t\in[0,T]} d_w(\nu_n(t), \nu(t))\to 0$. 
As for $s_n^{(1)}$ observe that 
\[ 
2s_n^{(1)} 
\le 
\int_0^Tdt\left| \lg g,\nu_n(t)\rg-\lg g,\nu(t)\rg \right|, 
\] 
where 
$\ds{g(x)=\int_0^{x}dy\Hc(y,x-y) 
(\Box h_i^{\star})(y,x-y)\varphi_R(x)}$. 
Since the assumed continuity of $\Hc$ and (\ref{3.6}) 
together assure that 
$g$ is a bounded continuous function on $(0,\infty)$, 
$\lg g,\nu_n(t)\rg \to \lg g,\nu(t)\rg $ 
for each $t\ge0$. Furthermore, 
\[ 
\sup_{t\in[0,T]} 
\left| \lg g, \nu_n(t)\rg -\lg g, \nu(t)\rg  \right|
\le 
2a\Vert g \Vert_{\infty}.   
\] 
So, the dominated convergence theorem 
proves that $s_n^{(1)}\to 0$ as $n\to \infty$.  
Basically a similar strategy can be adopted 
to $s_n^{(2)}$: 
\[ 
2s_n^{(2)} 
\le 
\int_0^Tdt\left| \lg G,\nu_n(t)^{\otimes 2}\rg 
-\lg G,\nu(t)^{\otimes 2}\rg  \right|, 
\] 
where 
$G(x,y):=\Hh(x,y) (\Box h_i^{\star})(x,y)\varphi_R(x+y)$ 
is verified to be bounded and continuous 
on $(0,\infty)^2$ thanks to (\ref{3.4}) 
as well as the assumption that $\Hh$ is continuous. 
Here, we claim that 
$\nu_n(t)^{\otimes 2} \stackrel{w}{\to}  
\nu(t)^{\otimes 2}$. In this respect,  
we rely on a slight generalization of 
Theorem 2.8 in \cite{Bi}, which implies 
in particular that if a sequence $\{p_n\}$ 
of probability measures on $(0,\infty)$ 
converges weakly to $p$, then 
$p_n^{\otimes 2} \stackrel{w}{\to} p^{\otimes 2}$. 
Generalizing this assertion to finite measures 
is easy by considering the normalized measures 
and it follows that $\lg G,\nu_n(t)^{\otimes 2}\rg 
\to \lg G,\nu(t)^{\otimes 2}\rg $ 
as $n\to \infty$ for each $t\ge0$. 
In addition,  
\[ 
\sup_{t\in[0,T]} 
\left| \lg G, \nu_n(t)^{\otimes 2}\rg 
-\lg G, \nu(t)^{\otimes 2}\rg  \right|
\le 
2a^2\Vert G \Vert_{\infty}. 
\] 
Hence again by the dominated convergence 
theorem $s_n^{(2)}\to 0$. 
Consequently we have proved 
the continuity of $\Upsilon=\Upsilon_{T,f,R}$ 
on each $\cM_{\le a}$ 
with respect to $d_U$ for any $T,R>0$ and 
$f\in \{h_1^{\star}, h_2^{\star}, \ldots\}$. 
As was already discussed 
this implies (\ref{5.15}) for those $f$'s.  \par 
{\it Step 5.}~
The remaining task is derivation of 
the weak form (\ref{5.8}).  
By combining (\ref{5.15}) for $f=h_i^{\star}$ 
$(i=1,2,\ldots)$ with 
Lemma 5.6 (implying in particular (\ref{5.15}) 
for $f=h_0^{\star}$) and then 
recalling the relation $\xi(t)^{\star}=\mu(t)$,  
we have proved so far that, with probability 1, 
for any $i\in\Zp$ and $t \ge 0$ 
\begin{eqnarray*} 
\lg h_i,\mu(t)\rg  -\lg h_i,c_0^{\star} \rg  
& = & 
\frac{1}{2}\int_0^tds \int \mu(s)^{\otimes 2}(dxdy) 
\Hh(x,y)(\Box h_i^{\star})(x,y)  \\ 
& & 
-\frac{1}{2}\int_0^tds \int \mu(s)(dx)
\int_0^{x}dy\Hc(y,x-y) (\Box h_i^{\star})(y,x-y) 
\end{eqnarray*} 
with the integrals on the right side being 
absolutely convergent. Since $\{h_0,h_1,\ldots\}$ is 
measure-determining, the above equalities 
are regarded as an equality among finite measures. 
In other words, one can replace $h_i$ by 
arbitrary bounded Borel functions $f$. 
In particular, replacing $h_i$ by 
$f_{\star}$ with $f\in B_c$ being arbitrary, 
we obtain (\ref{5.8}).  
The proof of Theorem 5.7 is complete.      \qed 

\medskip 

\noindent 
{\it Remark.}~
As mentioned earlier the uniqueness of weak 
solutions of  (\ref{4.5}) has not been proved. 
There is quite an extensive literature concerning 
the existence and/or uniqueness of solutions to 
coagulation-fragmentation equations. 
Well-posedness in the sense of 
measure-valued solutions was studied in e.g. \cite{FL} 
(for coagulation equations) 
and \cite{C2} (for coagulation 
multiple-fragmentation equations). 
Below we take up a result in \cite{BLL}. 
While it considers classical solutions, 
the setting of that paper is well adapted to 
a special case of our models. 
Given $\lambda\ge 0$, let 
$\Hc(x+y)=2 (x+y)^{\lambda}$ 
or $F(x+y)=2 (x+y)^{\lambda+1}$. 
Then (\ref{1.1}) 
is rewritten into 
\begin{eqnarray} 
\frac{\partial}{\partial t}c(t,x)
& = & 
-x^{2+\lambda} c(t,x) 
+2\int_x^{\infty}y^{\lambda+1}c(t,y)dy    \label{5.21}  \\ 
& & 
+\frac{1}{2}\int_0^x K(y,x-y)c(t,y)c(t,x-y)dy  
-c(t,x)\int_0^{\infty}K(x,y)c(t,y)dy.      \nonumber  
\end{eqnarray} 
This coincides with the equation (1.3) in \cite{BLL} 
with $\alpha=\lambda+2$ and $\nu=0$. 
One assumption made on $K$ in \cite{BLL} (cf. (1.5) there) 
is the following: \\ 
for some $C>0$ and 
$0\le \sigma \le \rho<\alpha(=\lambda+2)$ 
\[ 
0 \le K(x,y) \le C 
\left[(1+x)^{\rho}(1+y)^{\sigma}+
(1+x)^{\sigma}(1+y)^{\rho}\right], 
\quad x,y>0. 
\] 
Under our assumptions (\ref{1.4}) and (H1) on $K$, 
this condition is fulfilled with $C=\Ch C_{1,\lambda}$, 
$\rho=\lambda+1$ and $\sigma=1$. Therefore, 
Theorem 2.2 (i) in \cite{BLL} implies, among other things, 
the existence of a unique nonnegative classical solution 
to (\ref{5.21}) which is local in time. As argued  
in a closely related article \cite{BL}, showing 
the existence of a global solution requires 
a priori bound for the moments of solutions. 
For such a purpose, an analogue to  
our calculations in \S 4.2 could be useful, 
although we will not pursue this point here. 

\bigskip 
\bigskip 

\noindent 
{\Large \bf Appendix} 

\medskip 
\medskip 
\noindent 
In this appendix we prove Lemma 4.1. \\ 
{\it Proof of Lemma 4.1.}~
By the assumption for any $f\in B_{+,c}$ 
\[ 
E\left[\exp(-\lg f,\eta^N\rg )\right] 
=E\left[\exp\left(-\frac{1}{N}\sum_if(Y_i^N)\right)\right] 
=\exp\left(-N\lg 1-e^{-f/N}, \zeta\rg \right). 
\] 
As $N\to \infty$ 
the most right side converges to 
$\exp(-\lg f, \zeta\rg)$ by Lebesgue's convergence theorem. 
This proves the first assertion. To prove 
the convergence (\ref{4.4}) in case $a(=:n)\in \N$ 
it suffices to give the moment formula of the form  
\[ 
E\left[\left(\sum_iY_i^N\right)^n\right]
=
n!\sum_{k=1}^n\frac{N^k}{k!}\sum 
\frac{\lg\psi_{n_1}, \zeta\rg\cdots \lg\psi_{n_k}, \zeta\rg}
{n_1!\cdots n_k!},            
\] 
where the inner summation is taken over 
$k$-tuples $(n_1,\ldots,n_k)$ of positive integers 
such that $n_1+\cdots+n_k=n$.  
But, the proof of the above formula clearly reduces to 
showing that for each $R>0$ 
\[ 
E\left[\left(\sum_iY_i^N\one_{(0,R]}(Y_i^N)
\right)^n\right]
=n!\sum_{k=1}^n\frac{N^k}{k!}\sum 
\frac{\lg\psi_{n_1}^{(R)}, \zeta\rg
\cdots \lg\psi_{n_k}^{(R)}, \zeta\rg}{n_1!\cdots n_k!},  
\] 
where $\psi_{l}^{(R)}(y)=y^l\one_{(0,R]}(y)$, and  
this version  is derived by comparing the coefficients 
of $t^n$ after expanding in $t$ each side of 
\[ 
E\left[\exp\left(t\sum_iY_i^N\one_{(0,R]}(Y_i^N)
\right)\right] 
=\exp\left(N\int_{(0,R]} \zeta(dy)(e^{ty}-1) \right). 
\]  

It remains to prove (\ref{4.4}) in case  
$a\in (1,\infty)\setminus \N$. 
For such an $a$, put $n=[a]$ and $\alpha=a-[a]\in (0,1)$. 
Let $C_{3,\alpha}=\alpha/\Gamma(1-\alpha)$. 
Combining the aforementioned moment formula 
with the identity 
$y^{\alpha}=C_{3,\alpha} 
\int ds s^{-(1+\alpha)}(1-e^{-sy})$ 
for $y>0$ we deduce 
\begin{eqnarray*} 
\lefteqn{
E\left[\left(\sum_iY_i^N\right)^a\right]}            \\  
&=& 
C_{3,\alpha}\int \frac{ds}{s^{1+\alpha}}
\left\{E\left[\left(\sum_iY_i^N\right)^n\right] 
-E\left[\left(\sum_iY_i^N\right)^n
\exp\left(-s\sum_iY_i^N\right)\right] \right\}            \\ 
&=& 
C_{3,\alpha}n!\int \frac{ds}{s^{1+\alpha}}
\left\{\sum_{k=1}^n\frac{N^k}{k!} \sum 
\frac{\int y^{n_1}\zeta(dy)\cdots \int y^{n_k}\zeta(dy)}
{n_1!\cdots n_k!}\right.                                                  \\
& & 
-\left.\sum_{k=1}^n\frac{N^k}{k!}\sum 
\frac{\int y^{n_1}e^{-sy}\zeta(dy)\cdots \int y^{n_k}e^{-sy}\zeta(dy)}
{n_1!\cdots n_k!}
E\left[\exp\left(-s\sum_iY_i^N\right)\right] \right\}. 
\end{eqnarray*} 
Here, the second equality follows from the fact that 
under the transformed measure 
\[ 
\tilde{P}_s\left(\sum \delta_{Y^N_i}\in\bullet \right) 
:= E\left[\exp\left(-s\sum_iY_i^N\right); 
\ \sum \delta_{Y^N_i}\in\bullet
 \right] 
E\left[\exp\left(-s\sum_iY_i^N\right)\right]^{-1}
\] 
$\sum \delta_{Y^N_i}$ is 
${\rm Po}(e^{-sy}\zeta(dy))$-distributed. 
(See e.g. \cite{B}, p.80, Lemma 2.4.) 
For $s>0$ set 
$\Lambda(s) =\int(1-e^{-sy})\zeta(dy)$ so that 
$\Lambda(s)\le s\int y\zeta(dy)$ and  
\[ 
E\left[\exp\left(-s\sum_iY_i^N\right)\right]
= \exp\left(-N\Lambda(s)\right). 
\] 
By the change of variable $u:=Ns$ in the integral 
with respect to $ds$ 
\begin{eqnarray*} 
\lefteqn{
E\left[\left(\frac{\sum_iY_i^N}{N}\right)^a\right]}            \\  
&=& 
C_{3,\alpha}n!\int \frac{du}{u^{1+\alpha}}
\left\{\sum_{k=1}^n\frac{N^{k-n}}{k!}\sum 
\frac{\int y^{n_1}\zeta(dy)\cdots \int y^{n_k}\zeta(dy)}
{n_1!\cdots n_k!}\right.                                                  \\
& & 
-\left.\sum_{k=1}^n\frac{N^{k-n}}{k!}\sum 
\frac{\int y^{n_1}e^{-\frac{uy}{N}}\zeta(dy)\cdots 
\int y^{n_k}e^{-\frac{uy}{N}}\zeta(dy)}{n_1!\cdots n_k!} 
\exp\left(-N\Lambda\left(\frac{u}{N}\right)\right)\right\}  \\ 
&=& 
C_{3,\alpha}n!\sum_{k=1}^n\frac{N^{k-n}}{k!}\sum 
\int \frac{du}{u^{1+\alpha}}
\left\{\frac{\int y^{n_1}\zeta(dy)\cdots \int y^{n_k}\zeta(dy)}
{n_1!\cdots n_k!}\right.                                                  \\
& & 
\hspace*{2cm}\left. -\frac{\int y^{n_1}e^{-\frac{uy}{N}}\zeta(dy)
\cdots \int y^{n_k}e^{-\frac{uy}{N}}\zeta(dy)}{n_1!\cdots n_k!}
\exp\left(-N\Lambda\left(\frac{u}{N}\right)\right)\right\}.  
\end{eqnarray*} 
Here, we claim that each integral with respect to $du$ 
in the last expression is not only finite 
but also bounded in $N$.  
Indeed, for such an integral corresponding to 
$k=1$ (or $n_1=n$), we can observe by Fubini's theorem 
\begin{eqnarray*} 
\lefteqn{
C_{3,\alpha}\int \frac{du}{u^{1+\alpha}} 
\left|\int y^{n}\zeta(dy)
- \int y^{n}e^{-\frac{uy}{N}}\zeta(dy) 
\exp\left(-N\Lambda\left(\frac{u}{N}\right)\right) \right|} \\ 
& = & 
C_{3,\alpha}\int \frac{du}{u^{1+\alpha}} 
\left(\int y^{n}\zeta(dy)- \int y^{n}e^{-\frac{uy}{N}}\zeta(dy)\right) \\ 
& & 
+C_{3,\alpha}\int \frac{du}{u^{1+\alpha}} 
\int y^{n}e^{-\frac{uy}{N}}\zeta(dy)
\left\{1-\exp\left(-N\Lambda\left(\frac{u}{N}\right)\right)\right\}  \\ 
& \le & 
\int y^n\left(C_{3,\alpha}\int \frac{du}{u^{1+\alpha}}
(1-e^{-\frac{uy}{N}})\right)\zeta(dy)                     \\ 
& & 
+\int y^{n}\zeta(dy)C_{3,\alpha}
\int \frac{du}{u^{1+\alpha}} 
\left\{1-\exp\left(-u\int y\zeta(dy)\right)\right\}  \\ 
& = & 
\int y^n \left(\frac{y}{N}\right)^{\alpha}\zeta(dy)
+\int y^{n}\zeta(dy)\left(\int y \zeta(dy)\right)^{\alpha}. 
\end{eqnarray*} 
It is not difficult to obtain 
analogous estimates for any integral 
corresponding to $k\in\{2,\ldots,n-1\}$. Thus, 
we have identified the main term as $N\to\infty$: 
\begin{eqnarray*} 
\lefteqn{
E\left[\left(\frac{\sum_iY_i^N}{N}\right)^a\right]+o(1)}            \\  
& = & 
C_{3,\alpha}\int \frac{du}{u^{1+\alpha}}
\left\{\left(\int y\zeta(dy)\right)^n 
-\left(\int ye^{-\frac{uy}{N}}\zeta(dy)\right)^n
\exp\left(-N\Lambda\left(\frac{u}{N}\right)\right)\right\}.  
\end{eqnarray*} 
Lastly, this integral converges to 
$\lg \psi_1,\zeta\rg^{a}$ since as before 
\begin{eqnarray*} 
\lefteqn{ 
C_{3,\alpha}\int \frac{du}{u^{1+\alpha}}
\left|\left(\int y\zeta(dy)\right)^n 
-\left(\int ye^{-\frac{uy}{N}}\zeta(dy)\right)^n\right|} \\ 
& \le & 
\int y \left(
C_{3,\alpha}\int \frac{du}{u^{1+\alpha}} 
\left(1-e^{-\frac{uy}{N}}\right)\right)\zeta(dy)  
\cdot n\left(\int y\zeta(dy)\right)^{n-1}                   \\ 
& = & 
\int y\left(\frac{y}{N}\right)^{\alpha}\zeta(dy)
\cdot n\left(\int y\zeta(dy)\right)^{n-1} 
\ \to \ 0 
\end{eqnarray*} 
and by Lebesgue's convergence theorem 
\begin{eqnarray*} 
\lefteqn{
C_{3,\alpha}\int \frac{du}{u^{1+\alpha}} 
\left(\int ye^{-\frac{uy}{N}}\zeta(dy)\right)^n
\left\{1-\exp\left(-N\Lambda\left(\frac{u}{N}\right)\right)\right\}}  \\ 
&\to & 
C_{3,\alpha}\int \frac{du}{u^{1+\alpha}}
\left(\int y\zeta(dy)\right)^n
\left\{1-\exp\left(-u\int y \zeta(dy)\right)\right\}  \\ 
&=& 
\left(\int y\zeta(dy)\right)^{n+\alpha} 
\ = \ 
\lg \psi_1,\zeta\rg^{a}. 
\end{eqnarray*} 
Consequently (\ref{4.4}) holds 
and the proof of Lemma 4.1 is complete. 
\qed

\end{document}